\documentclass[10pt,a4paper]{article}

\usepackage[numbers,sort&compress]{natbib}
\usepackage{doi}
\usepackage{graphicx}
\usepackage{amsmath,amssymb,amsfonts}
\usepackage{amsthm}
\usepackage{mathrsfs}
\usepackage[title]{appendix}
\usepackage[dvipsnames]{xcolor} 
\usepackage{textcomp}
\usepackage{manyfoot}
\usepackage{algorithm}
\usepackage{algorithmicx}
\usepackage{algpseudocode}
\usepackage{listings}

\usepackage{cleveref}
\usepackage{float}

\usepackage{threeparttable}
\usepackage{booktabs}
\usepackage{caption}

\usepackage{fix-cm}
\newcommand{\arxivtablefont}{\normalfont\fontsize{8bp}{9.5bp}\selectfont}
\let\arxivoriginaltable\table
\let\endarxivoriginaltable\endtable
\renewenvironment{table}[1][tbp]{%
    \arxivoriginaltable[#1]%
    \begin{center}%
    \begin{threeparttable}%
    \arxivtablefont
    \frenchspacing
    \renewcommand{\raggedright}{\rightskip=0pt plus 1fil\parfillskip=0pt\relax}%
    \sloppy\adjdemerits=100\linepenalty=100
    \renewcommand{\footnotesize}{\fontsize{7pt}{8pt}\selectfont}%
}{%
    \end{threeparttable}%
    \end{center}%
    \endarxivoriginaltable
}
\DeclareCaptionFont{tablefont}{\normalfont\fontsize{10bp}{12bp}\selectfont}
\makeatletter
\def\arxivbottomrule{\noalign{\ifnum0=`}\fi
    \hrule\@height3pt\@width0pt
    \hrule\@height0.75pt
    \hrule\@height3pt\@width0pt
    \futurelet\@tempa\@xhline}
\makeatother
\let\bottomrule\arxivbottomrule
\AddToHook{cmd/appendix/after}{%
    \setcounter{table}{0}%
}

\usepackage{multirow}
\usepackage{subcaption}
\usepackage{array}
\usepackage{makecell}
\usepackage{enumitem}
\usepackage{pifont}

\newcommand{\xmark}{\ding{55}}

\usepackage{xurl}

\makeatletter
\def\paragraphfont{%
  \reset@font
  \fontsize{10bp}{12bp}%
  \bfseries\itshape
  \selectfont
  \raggedright
}
\renewcommand\paragraph{%
  \@startsection{paragraph}{4}{\z@}%
    {-12pt \@plus -4pt \@minus -2pt}%
    {3pt}%
    {\paragraphfont}%
}
\makeatother

\hypersetup{
  colorlinks,
  breaklinks=true,
  plainpages=false,
  citecolor=blue,
  linkcolor=blue,
  urlcolor=blue,
  bookmarksopen=true,
  bookmarksnumbered=false,
  bookmarksdepth=5
}

\newcommand{\eg}{\textit{e}.\textit{g}.}
\newcommand{\ie}{\textit{i}.\textit{e}.}
\newcommand{\FNN}{\mathrm{FNN}}
\newcommand{\CNN}{\mathrm{CNN}}
\newcommand{\UNet}{\mathrm{U}\text{-}\mathrm{Net}}
\newcommand{\GNN}{\mathrm{GNN}}

\newcommand{\IC}{\mathrm{IC}}
\newcommand{\ILU}{\mathrm{ILU}}

\newcommand{\secref}[1]{Sec.~\ref{#1}}

\usepackage{todonotes}

\usepackage{easyReview}

\usepackage[normalem]{ulem}

\AtBeginDocument{%
    \setlength{\abovedisplayskip}{4pt}%
    \setlength{\belowdisplayskip}{4pt}%
    \setlength{\abovedisplayshortskip}{4pt}%
    \setlength{\belowdisplayshortskip}{4pt}%
}%

\usepackage{xparse}
\NewDocumentEnvironment{ignore}{+b}
    {\ignorespaces}
    {\ignorespacesafterend}

\usepackage{geometry}
\usepackage{authblk}

\begin{document}

\title{An overview of machine learning-enhanced iterative methods for systems of linear and nonlinear equations}

\author[a]{Yuhuang Meng}
\author[a,b]{Jing Zhao}
\author[a,*]{Alexander Heinlein}

\affil[a]{Delft Institute of Applied Mathematics, Delft University of Technology, Mekelweg 4, Delft, 2628 CD, the Netherlands}
\affil[b]{Department of Hydrodynamics and Forecasting, Deltares, Boussinesqweg 1, Delft, 2629 HV, the Netherlands}

\date{}

\begin{center}
    \begin{minipage}{0.85\textwidth}
        \maketitle
    \end{minipage}
\end{center}

\begingroup
    \renewcommand{\thefootnote}{\fnsymbol{footnote}}
    \footnotetext[1]{\normalsize Corresponding author}
    \renewcommand{\thefootnote}{}
    \footnotetext[2]{\normalsize \textit{Email addresses:}
    \texttt{y.meng@tudelft.nl} (Yuhuang Meng),
    \texttt{jing.zhao@deltares.nl} (Jing Zhao),
    \texttt{a.heinlein@tudelft.nl} (Alexander Heinlein)}
\endgroup

\begin{abstract}
    Systems of equations arise in a wide range of scientific and engineering applications.
    The present work focuses on solvers for general systems of equations, including but not limited to those arising from partial differential equations.
    These systems can be broadly categorized into linear and nonlinear problems.
    For large linear systems, iterative solvers are generally preferred over direct methods due to the latter's superlinear growth of computational costs.
    Although convergence theory is well-developed under certain assumptions on the coefficient matrix, many classes of systems still pose open challenges.
    These difficulties become even more severe for systems of nonlinear equations, where nonlinear solvers typically rely on repeated linearization.
    For example, Newton's method may even converge quadratically near the solution; it can also converge slowly or diverge when the initial guess is not chosen appropriately.
    A wide range of solvers with diverse variants and hyperparameter settings exists, and the development of efficient and robust iterative methods remains an active area of research.
    Recently, machine learning (ML) techniques have been applied to enhance the efficiency of classical iterative methods while preserving their interpretability and reliability.
    We refer to these ML-enhanced iterative methods as \emph{hybrid iterative methods}, in the sense that they combine classical iterative methods with ML.
    This paper provides a comprehensive overview of state-of-the-art approaches to constructing \emph{hybrid iterative methods} for systems of both linear and nonlinear equations, while also discussing open challenges and outlining potential directions for future research.
    \\

\end{abstract}

\tableofcontents

\section{Introduction}\label{sec:introduction}

This section is organized as follows.
\Cref{sec:introduction-background} presents the background and motivation, \Cref{sec:introduction-concepts} provides a high-level overview of classical iterative methods and introduces the terminology used throughout the paper, and \Cref{sec:introduction-scope} defines the scope and organization of this review.

\subsection{Background and motivation}\label{sec:introduction-background}

Systems of equations are prevalent in various scientific and engineering applications, appearing in fields such as image processing, circuit analysis, chemical kinetics, structural mechanics, and fluid dynamics.
In scientific computing, they often arise from the discretization of partial differential equations (PDEs).
Accurately solving these systems is crucial for modeling, analyzing, and understanding real-world phenomena.
Notably, we focus on, but do not restrict ourselves to, systems arising from PDE discretizations.

Depending on the type of systems, iterative methods can be broadly classified into two categories: methods for systems of linear and nonlinear equations.
Linear systems can be solved by direct methods, such as Cholesky, LU, and QR factorizations; however, these approaches are often impractical for large-scale problems due to their high computational demands.
Iterative methods, by contrast, offer a more efficient and practical alternative, as each iteration is typically inexpensive and usually requires storing only one or a few vectors.
We use the term linear iterative methods to refer to iterative algorithms specifically designed for solving systems of linear equations.
Classical iterative methods for linear systems include both stationary methods, such as Jacobi and Gauss-Seidel, and non-stationary methods, such as Krylov subspace methods.
These iterative methods are known for their robustness and theoretical guarantees, supported by rigorous convergence and error analysis.
Nevertheless, their efficiency heavily depends on the choice of components, such as initial guesses, preconditioners, update strategies, and algorithmic parameters.
Selecting or constructing effective components remains challenging, as these choices are often problem-dependent and may require prior knowledge of the underlying system, such as sparsity pattern or conditioning.
Moreover, constructing such components can be computationally expensive, particularly for large-scale problems, which may limit the overall efficiency and robustness of classical iterative methods in practice.

Many nonlinear problems can be addressed using methods from nonlinear optimization; see, for instance,~\cite{nocedal2006numerical,ulbrich2012nichtlineare,boyd2004convex} for detailed overviews.
Here, we focus on methods commonly used in scientific computing, including fixed-point iterations, Newton-type methods, and gradient-based minimization of residual merit functions.
The convergence behavior of nonlinear methods is generally more delicate than in the linear setting: convergence results are often local, while global guarantees typically require additional assumptions such as convexity.
Global convergence guarantees require additional problem- and algorithm-specific assumptions.
For example, Newton's method may exhibit quadratic convergence near the solution, but it can also converge slowly or even diverge when the initial guess is poor.
These challenges motivate the search for more robust and efficient solution strategies.

In recent years, scientific machine learning (SciML)~\cite{baker2019workshop} has emerged as a broad research field at the intersection of machine learning (ML) and scientific computing, with the aim of improving the efficiency and accuracy of scientific and engineering computations.
Representative research directions include physics-informed machine learning (PIML)~\cite{raissi2019physicsinformed, karniadakis2021physicsinformed}, operator learning~\cite{lu2021learninga, li2021fourier}, data-driven model discovery~\cite{brunton2016discovering, raissi2018deep, chen2018neural}, and surrogate modeling~\cite{alizadeh2020managing, sturek2025surrogate}, among others.
In this review, we focus on a specific class of methods within SciML: ML-enhanced iterative methods for systems of linear and nonlinear equations.
These methods incorporate ML into classical iterative frameworks and are referred to as \emph{hybrid iterative methods (HIMs)} throughout the paper.
The term \emph{hybrid} emphasizes that ML enhances and complements the underlying classical method rather than replacing it entirely.

\subsection{Iterative methods and terminology}\label{sec:introduction-concepts}

To facilitate the discussion and categorization of \emph{HIMs}, we first present general iterative frameworks and introduce the terminology used throughout this review.

\paragraph{\textbf{General iterative framework for linear systems}}

Consider the linear system
\begin{equation*}
    A\mathbf{x} = \mathbf{b},
\end{equation*}
where $A\in\mathbb{R}^{N\times N}$ is the coefficient matrix, $\mathbf{b} \in \mathbb{R}^N$ denotes the right-hand side vector, and $\mathbf{x} \in \mathbb{R}^N$ is the unknown solution vector.
Throughout this paper, for simplicity, we focus on square systems, as considered in the vast majority of the reviewed studies.
The objective is to compute $\mathbf{x}$ for given $A$ and $\mathbf{b}$.
Starting from an initial guess $\mathbf{x}_0$, an iterative method generates a sequence of approximations according to
\begin{equation}\label{equ:prelims-general-iter}
    \mathbf{x}_{k+1} = \Phi_k(\mathbf{x}_k, \mathbf{x}_{k-1}, \cdots; \boldsymbol{\gamma}_k),
\end{equation}
where $\Phi_k$ denotes the update (or iteration) function, $(\mathbf{x}_k, \mathbf{x}_{k-1}, \cdots)$ represents the iterates accumulated up to iteration $k$, and $\boldsymbol{\gamma}_k$ collects the associated algorithmic parameters.
For many methods, the iteration can be expressed more specifically as
\begin{equation}\label{equ:prelims-update}
    \begin{aligned}
        \mathbf{x}_{k+1}
         & = \mathbf{x}_k + \alpha_k \mathbf{d}_k                                                            \\
         & = \mathbf{x}_k + \alpha_k \Psi_k(\mathbf{x}_k, \mathbf{x}_{k-1}, \cdots; \boldsymbol{\gamma}'_k),
    \end{aligned}
\end{equation}
where $\mathbf{d}_k:=\Psi_k(\mathbf{x}_k, \mathbf{x}_{k-1}, \cdots; \boldsymbol{\gamma}'_k) \in \mathbb{R}^N$ is the update direction, $\alpha_k \in \mathbb{R}$ is a step-size or relaxation parameter, and $\boldsymbol{\gamma}'_k$ denotes the parameters involved in constructing the update direction.
More detailed discussions of these iterative methods are provided in later sections; see~\Cref{sec:linear-core_components-prelims,sec:linear-precond-prelims,sec:linear-mgdd-prelims}.

Based on this representation, we introduce the following terminology:
\begin{itemize}
    \item Under this representation, we refer to the initial guess $\mathbf{x}_0$, the update function $\Phi_k$ (or $\Psi_k$), and the parameters $\boldsymbol{\gamma}_k$ as the three \emph{core components} of an iterative method.
          We use this terminology to organize different ML-enhanced approaches.
    \item 
          Preconditioning techniques accelerate convergence by transforming the original linear system into an equivalent system that is easier to solve.
          For example, a left-preconditioned system takes the form $M^{-1}A\mathbf{x} = M^{-1}\mathbf{b}$.
          Ideally, $M^{-1}$ approximates $A^{-1}$ sufficiently well that $M^{-1}A$ has favorable spectral properties, while remaining inexpensive to construct and apply.
          In this formulation, the preconditioner $M^{-1}$ becomes part of the update function $\Phi_k$ (or $\Psi_k$) and may therefore be regarded as one of its components.
          For example, the preconditioned Richardson iteration takes the form
          \[
              \mathbf{x}_{k+1} = \mathbf{x}_{k} + \alpha_k M^{-1}(\mathbf{b} - A \mathbf{x}_{k}),
          \]
          where the update direction is given by $\mathbf{d}_k = \Psi_k(\mathbf{x}_k; M^{-1})=M^{-1}(\mathbf{b} - A \mathbf{x}_{k})$, and $\alpha_k$ is the relaxation parameter.
          However, constructing an effective and computationally inexpensive preconditioner is often challenging and highly problem-dependent.
    \item 
          \emph{Subspace correction methods} accelerate convergence by reducing difficult error components through corrections on suitably chosen subspaces.
          This general framework encompasses multigrid, domain decomposition, deflation, and augmentation methods, which may be used either as standalone iterative solvers or as preconditioners.
          Their performance depends strongly on the construction of components such as smoothers, coarsening strategies, coarse spaces, and transfer operators, providing several opportunities for ML enhancement.
    \item 
          For the same problem, different iterative solvers, that is, different combinations of iterative schemes $\Phi_k$ (or $\Psi_k$) and preconditioners $M^{-1}$, may exhibit significantly different convergence behavior and computational costs.
          Since the resulting configuration space can be large and problem-dependent, manual selection often relies on expert knowledge and empirical tuning.
          The task of automatically choosing a suitable solver configuration using ML is referred to as \emph{automated solver selection}.
\end{itemize}

\paragraph{\textbf{General iterative framework for nonlinear systems}}

Consider the nonlinear system
\begin{equation*}
    F(\mathbf{x}) = 0,
\end{equation*}
where $F: \mathbb{R}^N \rightarrow \mathbb{R}^N$ is a nonlinear function and $\mathbf{x} \in \mathbb{R}^{N}$ is the unknown solution vector.
Nonlinear systems are typically solved numerically using iterative methods that seek an approximate solution $\mathbf{x}^* \in \mathbb{R}^N$ satisfying $F(\mathbf{x}^*) \approx 0$.
Analogously to~\Cref{equ:prelims-general-iter}, a general iterative method for nonlinear systems can be written as
\begin{equation}\label{equ:prelims-general-iter-nonlinear}
    \mathbf{x}_{k+1} = \Phi_k(\mathbf{x}_k, \mathbf{x}_{k-1}, \cdots; \boldsymbol{\gamma}_k).
\end{equation}
Commonly used nonlinear solvers include fixed-point iterations and Newton-type methods; gradient-based methods can also be applied to residual merit functions.
Further details on these nonlinear iterative methods are provided in~\Cref{sec:nonlinear-prelims}.

Much of the terminology introduced for linear solvers is also applicable to nonlinear solvers.
First, the \emph{core components} remain the initial guess, the update function, and the algorithmic parameters.
More importantly, nonlinear solvers are generally more sensitive to the choice of the initial guess than linear solvers.
Second, \emph{nonlinear preconditioning} may be employed to improve convergence.
Third, multigrid and domain decomposition methods can be extended to nonlinear problems.
Finally, \emph{solver selection} has received relatively limited attention for nonlinear systems in the literature.
This is partly because the set of commonly used nonlinear solvers is smaller and partly because most linearization-based nonlinear solvers rely on linear solvers at each iteration.

\subsection{Scope and organization of this paper}\label{sec:introduction-scope}

\begin{figure}[!htbp] 
    \centering
    \includegraphics[width=\textwidth]{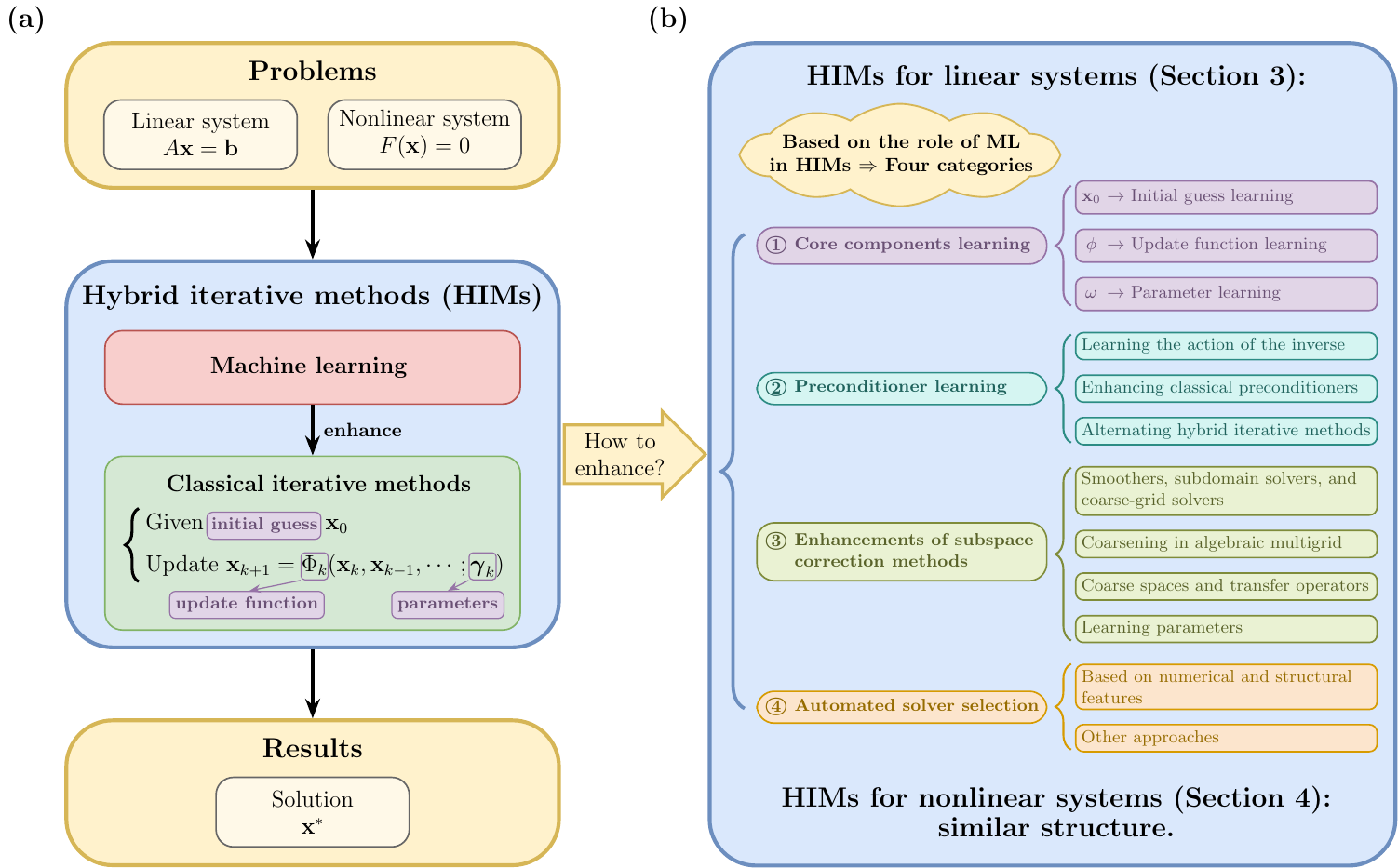}
    \caption{
        Overview of ML-enhanced iterative methods, referred to in this paper as \emph{hybrid iterative methods (HIMs)}.
        \textbf{(a)} General workflow of \emph{HIMs} for solving linear and nonlinear systems of equations.
        \textbf{(b)} Classification of \emph{HIMs} for linear systems into four categories, as discussed in~\Cref{sec:linear}.
        The classification of \emph{HIMs} for nonlinear systems, discussed in~\Cref{sec:nonlinear}, follows a similar structure.
    }
    \label{fig:overview-HIMs}
\end{figure}

The scope of this review and the main categories of \emph{HIMs} considered are illustrated in~\Cref{fig:overview-HIMs}(a).
Rather than entirely replacing classical numerical solvers, most existing \emph{HIMs} leverage ML techniques to construct or select components of classical iterative methods, where the construction or choice of these components may depend strongly on the problem instance.
In this way, \emph{HIMs} combine the robustness and, when available, convergence guarantees of classical numerical methods with the flexibility offered by ML.
For example, in contrast to fully learned approaches for solving PDEs, such as physics-informed neural networks (PINNs)~\cite{raissi2019physicsinformed}, classical iterative methods equipped with ML-based preconditioners may retain the accuracy and reliability of the underlying numerical solver.
Moreover, ML-based preconditioners can be trained over a range of parameter values, making them particularly attractive for parametric linear systems.
This hybridization makes \emph{HIMs} a promising class of methods for solving systems of equations arising from challenging scientific and engineering applications.

As shown in~\Cref{fig:overview-HIMs}(b), \emph{HIMs} can be broadly categorized into methods for linear systems and methods for nonlinear systems, which are discussed in~\Cref{sec:linear,sec:nonlinear}, respectively.
Following the discussion in~\Cref{sec:introduction-concepts}, we classify \emph{HIMs} for linear systems into four categories:
\begin{enumerate}
    \item \emph{Core components learning}: using ML to construct the initial guess, update function, or algorithmic parameters.
    \item \emph{Preconditioner learning}: using ML to construct or improve preconditioners from data, including two main approaches, namely \emph{learning the action of the inverse} and \emph{enhancing classical preconditioners}, as well as one usage strategy referred to as \emph{alternating hybrid iterative methods}.
    \item \emph{Enhancements of subspace correction methods}: using ML to design components of multigrid, domain decomposition, deflation, and augmentation methods, such as smoothers, subdomain solvers, coarsening strategies, coarse spaces, and transfer operators.
    \item \emph{Automated solver selection}: using ML to select a suitable solver, including an iterative method, a preconditioner, or both, from a prescribed search space based on problem features.
\end{enumerate}

We note that this categorization is not unique and that the different categories are not mutually exclusive.
First, \emph{preconditioner learning} may be regarded as a subcategory of \emph{update function learning}, since ML models can be leveraged to construct update functions through learned preconditioners.
Nevertheless, \emph{preconditioner learning} has developed into a broad and established research area and is therefore discussed separately.
Second, the \emph{enhancement of subspace correction methods} may overlap with \emph{preconditioner learning} when the resulting methods are employed as preconditioners for classical iterative solvers.
Third, although both \emph{update function learning} and \emph{automated solver selection} concern update functions, they operate at different levels.
The former aims to approximate or enhance classical iterative methods, whereas the latter formulates the task as a classification problem whose goal is to select an appropriate iterative method, preconditioner, or combination of the two.

Existing reviews focus on specific topics, such as the combination of domain decomposition methods and ML~\cite{heinlein2021combining, klawonn2024machine}, or on broader SciML topics, such as PIML~\cite{karniadakis2021physicsinformed, cuomo2022scientific, meng2025when} and neural operator learning~\cite{lu2022comprehensive, viswanath2026neural}.
However, a systematic review of ML-enhanced classical iterative methods for systems of equations is still lacking.
This review fills this gap by organizing existing approaches according to how ML is incorporated into classical iterative frameworks, covering both linear and nonlinear systems.
We also highlight their benefits, limitations, and open challenges.

The paper is organized as follows.
\Cref{sec:prelims-ml} introduces relevant ML techniques, with an emphasis on their use in hybrid iterative methods for linear and nonlinear systems of equations.
In~\Cref{sec:linear}, we review the literature on \emph{HIMs} for linear systems according to four categories: core components learning, preconditioner learning, enhancement of subspace correction methods, and automated solver selection.
\Cref{sec:nonlinear} focuses on \emph{HIMs} for nonlinear systems and follows a classification similar to that used for linear systems, except for automated solver selection, for which relatively little literature is available.
Other topics related to \emph{HIMs} but outside the main scope of this review, including ML-enhanced time integration methods, ML approaches inspired by classical iterative methods, and ML-enhanced direct solvers for linear systems, are briefly discussed in~\Cref{sec:related}.
Finally, \Cref{sec:conclusions} concludes the review and presents broader perspectives for future research.

\section{Machine learning preliminaries for hybrid iterative methods}\label{sec:prelims-ml}

This section briefly introduces the machine learning tools commonly used in \emph{HIMs}, including proper orthogonal decomposition (POD), feedforward neural networks (FNNs), convolutional neural networks (CNNs), graph neural networks (GNNs), recurrent neural networks (RNNs), and meta-learning.
The discussion is not intended to be exhaustive, but focuses on the techniques most relevant to the \emph{HIMs} reviewed in this paper.
These techniques are presented from the perspective of scientific computing, with an emphasis on their key characteristics and their connections to classical iterative methods.
We begin by introducing some basic concepts in machine learning.

For the purposes of this review, ML-based approaches are characterized by three closely related ingredients: \emph{data}, \emph{models}, and \emph{learning}~\cite[p.~11]{deisenroth2020mathematics}.
\emph{Data} can be viewed as samples drawn from an underlying (unknown) distribution.
ML methods are data-driven approaches that aim to extract patterns from data.
\emph{Models} define a space of possible functions that can represent the data, \eg, a parameterized function.
\emph{Learning} (or training) refers to adjusting the parameters of the model so that it captures patterns in the observed data (training data), where model performance is typically evaluated using a predefined metric.
From a mathematical perspective, this corresponds to optimizing the model parameters with respect to an objective function.
In this sense, ML aims to learn a function or representation from data that captures the underlying structure of the data.

To illustrate these concepts, we consider a simple supervised regression problem~\cite{abu-mostafa2012learning}.
Let the input be denoted by $\mathbf{x} \in \mathcal{X}$, and let $f:\mathcal{X} \to \mathcal{Y}$ be an unknown target function, where $\mathcal{X} \subseteq \mathbb{R}^m$ and $\mathcal{Y} \subseteq \mathbb{R}^n$ denote the input and output spaces, respectively.
Throughout this subsection, $\mathbf{x}$ denotes an ML input, except in expressions involving a linear system such as $A\mathbf{x}=\mathbf{b}$, where it denotes the unknown solution vector.
We are given a training dataset $\mathcal{D}=\{(\mathbf{x}_i, \mathbf{y}_i)\}_{i=1}^M$, where $\mathbf{y}_i = f(\mathbf{x}_i) \in \mathbb{R}^n$ denotes the corresponding output (or label), and $M$ is the number of training samples.
The goal is to find a model $g \in \mathcal{G}$ that approximates $f$, where $\mathcal{G}$ represents a hypothesis space of candidate functions.
In practice, the model is parameterized as $g_\theta$, where $\theta$ denotes the model parameters.
Training the model is achieved by solving an optimization problem over $\theta$ with respect to a predefined loss function.
A common choice is the mean squared error (MSE), defined as
\[
    \mathcal{L}(\theta) = \frac{1}{M}\sum_{i=1}^M \| g_\theta(\mathbf{x}_i) - \mathbf{y}_i \|_2^2,
\]
which measures the discrepancy between the model predictions and the target outputs.
Once trained, the model can be used to approximate $f$ at a new input $\mathbf{x}' \in \mathcal{X}$ via $f(\mathbf{x}') \approx g_\theta(\mathbf{x}')$.

ML techniques are commonly categorized into three paradigms: supervised learning, unsupervised learning, and reinforcement learning.
In supervised learning, models are trained using labeled data, where corresponding labels $\mathbf{y}_i$ are provided.
Typical supervised learning tasks are regression and classification.
Unsupervised learning aims to extract latent structures or patterns from unlabeled data.
Representative examples include clustering and dimensionality reduction techniques, such as the $k$-means algorithm and principal component analysis (PCA).
Reinforcement learning (RL) focuses on learning decision-making strategies through interactions of an agent with an environment.
Instead of relying on labeled outputs, an RL agent selects an action $\mathbf{a}_t$ based on the current state $\mathbf{s}_t$ and receives a reward $r_t$ that evaluates the selected action, where $t$ denotes the interaction time step.
For further technical details on ML principles, we refer to~\cite{abu-mostafa2012learning, deisenroth2020mathematics} and, in particular, to~\cite{goodfellow2016deep, bishop2024deep} for deep learning.
Next, we introduce several ML methods that are commonly used in the construction of \emph{HIMs}.

\paragraph{\textbf{Proper orthogonal decomposition}}

POD~\cite{lumley1967structure} is a popular linear dimensionality reduction technique in unsupervised learning.
In the context of statistics, it is also known as principal component analysis (PCA)~\cite{hotelling1933analysis, pearson1901liii}.
Consider a high-dimensional dataset $\mathcal{D}=\{\mathbf{x}_i\}_{i=1}^M$,
where $\mathbf{x}_i \in \mathbb{R}^m$ and $m \gg 1$.
The sample mean is defined as $\boldsymbol{\mu} = \frac{1}{M}\sum_{i=1}^M \mathbf{x}_i$.
The goal of POD is to extract dominant features by projecting the centered data $\mathbf{x}_i - \boldsymbol{\mu} \in \mathbb{R}^m$ onto a lower-dimensional subspace $\mathbb{R}^r$ with $r\ll m$.

Specifically, we seek a set of orthonormal basis vectors $\{\boldsymbol{\phi}_j\}_{j=1}^r$, $\boldsymbol{\phi}_j \in \mathbb{R}^m$, such that any data sample $\mathbf{x}$ can be approximated as
\[
    \mathbf{x} \approx \boldsymbol{\mu} + \sum_{j=1}^r z_j \boldsymbol{\phi}_j.
\]
The coefficients are given by the orthogonal projection $z_j = \boldsymbol{\phi}_j^T (\mathbf{x}-\boldsymbol{\mu})$.
Compared with the original high-dimensional data, this reduced representation significantly decreases both memory requirements and computational cost by working in a low-dimensional space.
Let $\Phi = [\boldsymbol{\phi}_1, \cdots, \boldsymbol{\phi}_r] \in \mathbb{R}^{m\times r}$ denote the POD basis matrix.
The low-dimensional representation (also called latent variables) is $\mathbf{z} = \Phi^T (\mathbf{x}-\boldsymbol{\mu}) \in \mathbb{R}^r$, and the reconstruction is $\hat{\mathbf{x}} = \boldsymbol{\mu} + \Phi \mathbf{z} $.
Define the centered data matrix $X=[\mathbf{x}_1-\boldsymbol{\mu},\dots,\mathbf{x}_M-\boldsymbol{\mu}] \in \mathbb{R}^{m\times M}$.
The POD basis $\Phi$ is obtained by minimizing the reconstruction error, given by
\[
    \begin{aligned}
        \min_{\Phi}\quad & \frac{1}{M} \left\|X - \Phi\Phi^T X \right\|_F^2 \\
        \text{s.t.}\quad & \Phi^T\Phi = I_r,
    \end{aligned}
\]
where $I_r \in \mathbb{R}^{r\times r}$ denotes the identity matrix.
Let $X = U\Sigma V^T$ be the singular value decomposition (SVD) of $X$.
An optimal POD basis can be chosen as $\Phi = U_r$, where $U_r$ contains the first $r$ left singular vectors.
Typical applications of POD include dimensionality reduction, noise reduction, and feature extraction.
Within the scope of this paper, POD has been used to construct reduced-order-model (ROM)-based preconditioners for Krylov subspace methods~\cite{pasetto2017reduced} and to develop nonlinear preconditioning strategies~\cite{luo2023pin^la-settled}.

\paragraph{\textbf{Feedforward neural networks}}
FNNs, also known as multilayer perceptrons (MLPs) with fully connected architectures, are widely used nonlinear function approximators~\cite{goodfellow2016deep}.
The term \textit{feedforward} means that information flows from the input layer to the output layer without feedback connections.
We consider an FNN $\mathcal{N}_\theta(\mathbf{x})$ that approximates a target function $f(\mathbf{x})$.

\begin{figure}[!htbp] 
    \centering
    \includegraphics[width=0.35\textwidth]{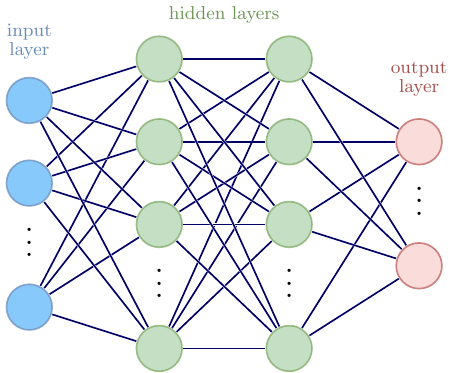}
    \caption{
        Schematic of a fully connected feedforward neural network (FNN) with two hidden layers.
        Nonlinearity is introduced through neuron activation functions, while the output layer is often linear.
    }
    \label{fig:prelims-ml-fcns}
\end{figure}

As illustrated in~\Cref{fig:prelims-ml-fcns}, an FNN is composed of an input layer, multiple hidden layers, and an output layer.
It can be represented as a composition of layerwise mappings,
\begin{equation*}
    \mathcal{N}_\theta (\mathbf{x}) = F_L \circ F_{L-1} \circ \cdots \circ F_1(\mathbf{x}),
\end{equation*}
where
\begin{equation*}
    \begin{cases}
        F_l(\mathbf{x})=\sigma(W_{l}\mathbf{x} + b_{l}), & l=1, \cdots, L-1, \\
        F_L(\mathbf{x}) = W_{L}\mathbf{x} + b_{L}.
    \end{cases}
\end{equation*}
Here, $W_{l}$ and $b_{l}$ denote the weight matrix and bias vector of the $l$-th layer, respectively, and the network parameters are given by $\theta = \{W_l, b_l\}_{l=1}^L$.
The function $\sigma(\cdot)$ is called the activation function, which introduces nonlinearity and allows the network to approximate complex functions.
Common activation functions include the sigmoid function $\sigma(\mathbf{x}) = \frac{1}{1+e^{-\mathbf{x}}}$, the hyperbolic tangent function $\sigma(\mathbf{x}) = \tanh(\mathbf{x})$, and the rectified linear unit (ReLU) function $\mathrm{ReLU}(\mathbf{x}) = \max(0, \mathbf{x})$.
The network output is denoted by $\hat{\mathbf{y}} = \mathcal{N}_\theta (\mathbf{x})$.
Networks with multiple hidden layers are commonly referred to as deep neural networks (DNNs), which are capable of learning complex representations.

Training a neural network involves optimizing the parameters $\theta$ by minimizing a predefined loss function over a training dataset $\mathcal{D}=\{(\mathbf{x}_i, \mathbf{y}_i)\}_{i=1}^M$.
This optimization problem is typically solved using gradient-based methods, most notably stochastic gradient descent (SGD), and higher-order methods, such as limited-memory Broyden-Fletcher-Goldfarb-Shanno (L-BFGS)~\cite{liu1989limited}.
Various first-order optimization methods have been developed to improve convergence and training stability, including momentum-based methods and adaptive learning rate algorithms such as AdaGrad~\cite{duchi2011adaptive} and Adam~\cite{kingma2017adam}; see also~\cite{kopanicakova2026introduction} for an overview of optimization methods for training SciML models.

An important characteristic of neural networks' training dynamics is their spectral bias, also known as the frequency principle: during training, they tend to learn low-frequency components of a target function more rapidly than high-frequency ones~\cite{rahaman2019spectral, xu2020frequency, xu2025overview, xu2025understanding}.
Although this behavior may limit their ability to resolve fine-scale features, it can be exploited in \emph{HIMs}.
For many PDE discretizations, suitable stationary relaxation schemes efficiently damp high-frequency error components while reducing low-frequency components more slowly.
Neural networks are widely observed to exhibit a spectral bias toward learning low-frequency components first, which can complement the smoothing properties of these relaxation schemes.
This idea is exploited, for example, in \emph{initial guess learning} (\Cref{sec:linear-core_components-initial}), where neural-network predictions are refined by stationary iterations~\cite{xu2020frequency}, and in the \emph{alternating HIM} (\Cref{sec:linear-precond-alternating}) proposed in~\cite{zhang2024blending}, which is based on deep operator networks (DeepONets)~\cite{lu2021learninga}.
This is conceptually reminiscent of multigrid methods, in which coarse-grid correction addresses smooth, low-frequency errors, while relaxation on the fine grid damps oscillatory, high-frequency errors.

Despite the spectral bias, FNNs are widely used as nonlinear function approximators.
However, FNNs operate on vectorized inputs and outputs and do not explicitly exploit the structural relationships among data elements, \ie, the entries of an input $\mathbf{x}_i$.
This limitation may reduce their flexibility and effectiveness when handling structured data.
For instance, images and graphs exhibit strong locality properties: images are defined on grid-like domains, whereas graphs are characterized by irregular relational structures.
In addition, FNNs are not well-suited for modeling sequential data.
These data structures frequently arise in \emph{HIMs}, for example, through image or graph representations of coefficient matrices and iterative sequences.
To better capture such different types of structure, specialized architectures such as CNNs, GNNs, and RNNs have been developed for image-like, graph-structured, and sequential data, respectively.
These architectures are briefly introduced below.

\paragraph{\textbf{Convolutional neural networks}}
CNNs~\cite{lecun1989backpropagation, lecun1998gradient} are a class of neural network models originally designed for processing data defined on regular grids, such as images.
They have been widely applied to image processing tasks, including image classification and segmentation.

\begin{figure}[!htbp]
    \centering
    \includegraphics[width=0.6\textwidth]{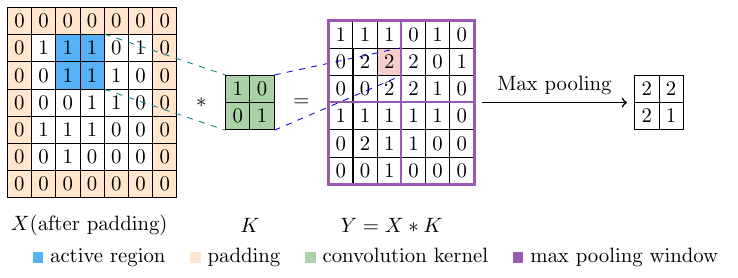}
    \caption{
        Illustration of a 2D convolution and max pooling.
        The input $X$ has size $5\times5$ and is extended to $7\times7$ by zero padding, analogous to homogeneous Dirichlet boundary conditions, whereas padding with repeated boundary values resembles homogeneous Neumann boundary conditions.
        A convolution with a $2\times2$ kernel $K$ produces a $6\times6$ feature map, which is then downsampled by taking the maximum over each $3\times3$ subblock.
    }
    \label{fig:prelims-ml-conv}
\end{figure}

CNNs consist of convolutional layers, pooling layers, and activation functions.
An illustration of the 2D convolution operation, padding, and pooling is shown in~\Cref{fig:prelims-ml-conv}.
Convolutional layers extract local features by applying convolutional kernels to the input.
For a 2D input image $X \in \mathbb{R}^{H \times W}$ and a kernel $K \in \mathbb{R}^{k_h \times k_w}$, the convolution operation without padding and with unit stride produces a feature map $Y$ given by
\[
    Y_{i,j} = (X*K)_{i,j} = \sum_{m=0}^{k_h-1} \sum_{n=0}^{k_w-1} X_{i+m, j+n} K_{m,n},
\]
where $0\le i\le H-k_h$, $0\le j\le W-k_w$, and * denotes the cross-correlation operator commonly used in CNNs.
This operation enables CNNs to effectively capture spatial locality.
The padding layer extends the input with additional boundary values (typically zeros).
It is commonly used to preserve the spatial resolution of feature maps, as the convolution operation without padding reduces the spatial dimensions of the output.
The pooling layer performs local aggregation operations; for example, max pooling selects the maximum value within each local region.
It is used to downsample feature maps and reduce computational cost.
This operation is analogous to the restriction step in multigrid or domain decomposition methods, where fine-grid information is transferred to a coarser representation.

\begin{figure}[!htbp]
    \centering
    \includegraphics[width=0.75\textwidth]{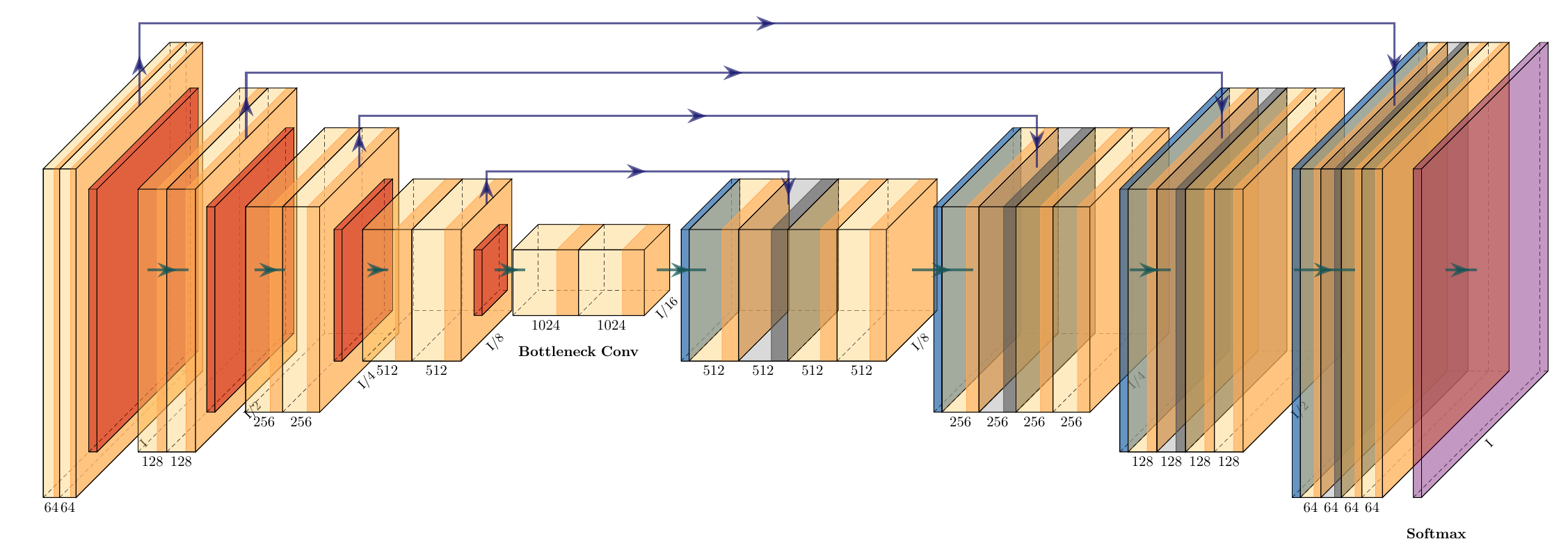}
    \caption{
        Illustration of the U-Net architecture, taken from PlotNeuralNet~\cite{plotneuralnet}.
        The architecture consists of an encoder path and a decoder path, where corresponding refinement levels are connected by skip connections. This structure resembles the process in a multigrid V-cycle.
    }
    \label{fig:prelims-unet}
\end{figure}

U-Net~\cite{ronneberger2015unet} is a specialized CNN architecture designed for image segmentation.
Similar to autoencoders, a U-Net has a U-shaped structure consisting of a contracting path (encoder) and an expanding path (decoder).
Skip connections are utilized to combine high-resolution features from the encoder with upsampled features in the decoder, which improves segmentation accuracy.
This encoder-decoder structure is reminiscent of a multigrid V-cycle, as illustrated in~\Cref{fig:prelims-unet}.
The connection between CNN architectures and multigrid methods has been further explored in MgNet~\cite{he2019mgnet}.

\begin{figure}[!htbp]
    \centering
    \includegraphics[width=1.0\textwidth]{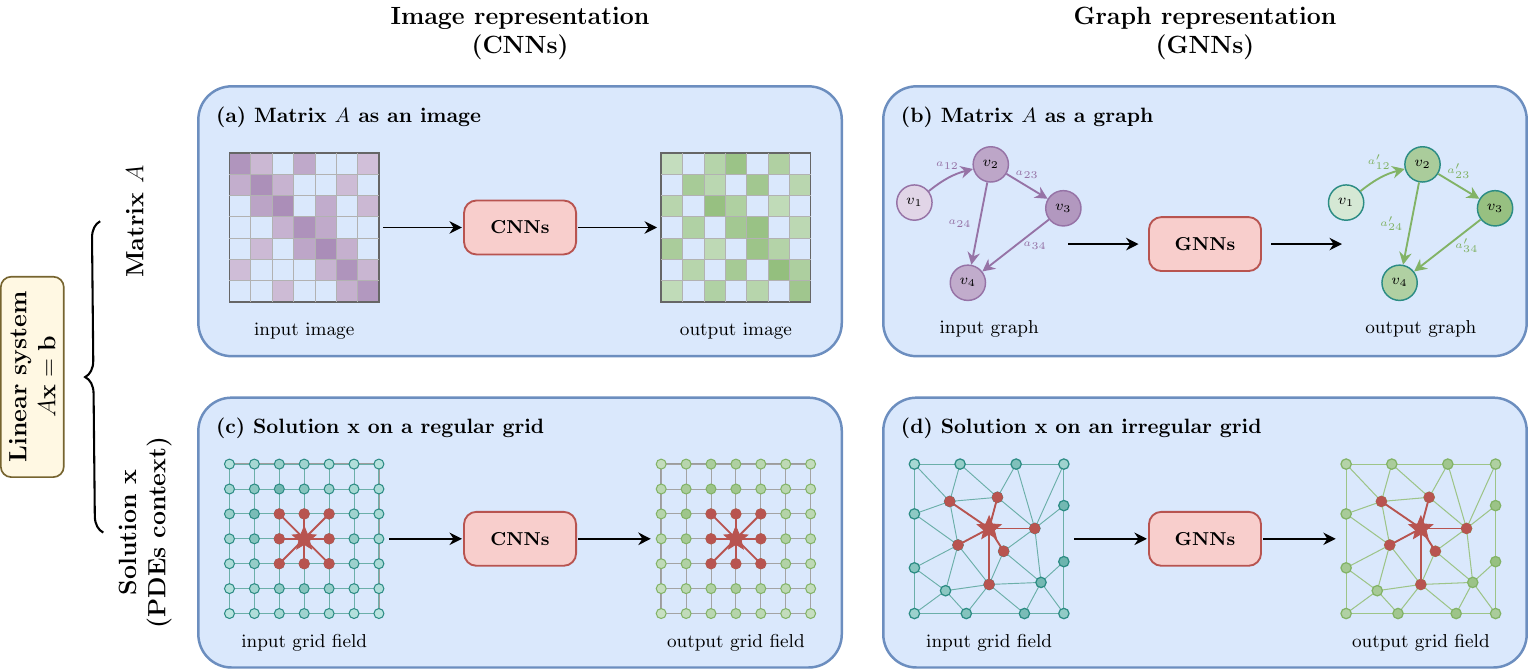}
    \captionsetup{skip=8pt}
    \caption{
        Applications of CNNs and GNNs to the linear system $A\mathbf{x}=\mathbf{b}$.
        \textbf{(a)} The coefficient matrix $A$ is represented as an image and processed by a CNN, where blank cells denote zero entries of $A$.
        \textbf{(b)} The matrix $A$ is represented as a graph, in which matrix entries define edge connections and attributes, and is processed by a GNN.
        \textbf{(c)} For a PDE discretized on a uniform Cartesian grid, the solution $\mathbf{x}$ is represented as a grid-based field and can be naturally processed by CNNs.
        \textbf{(d)} For a PDE discretized on an irregular (also called unstructured) mesh, the solution is represented on a graph and can be handled more flexibly by GNNs.
        The red star denotes a representative degree of freedom, and the red nodes and connecting edges indicate its local neighborhood.
        Thus, CNNs exploit the fixed spatial arrangement of regular grids, whereas GNNs aggregate information over general connectivity patterns.
    }
    \label{fig:prelims-struct_unstruct_data}
\end{figure}

The application of CNNs to solving systems of linear equations $A\mathbf{x}=\mathbf{b}$ can be viewed from two perspectives:
\begin{enumerate}
    \item 
          The coefficient matrix $A$ can be interpreted as an image, enabling CNNs to approximate (potentially nonlinear) mappings between matrices and vectors (or matrices), as shown in~\Cref{fig:prelims-struct_unstruct_data}(a).
          For example, in~\cite{sappl2019deep}, CNNs take the coefficient matrix as input and output a corresponding preconditioner in matrix form.
    \item 
          For linear systems arising from PDE discretizations, the solution $\mathbf{x}$ defined on a regular (uniform) grid can also be treated as an image, as illustrated in~\Cref{fig:prelims-struct_unstruct_data}(c).
          In this setting, the local and sparse connectivity of CNNs is well suited to exploiting the local spatial structure commonly induced by PDE discretizations.
          In~\cite{stanziola2021helmholtz}, a modified U-Net architecture is employed to learn the mapping from the current solution and residual to its corresponding update.
\end{enumerate}

\paragraph{\textbf{Graph neural networks}}

Building upon CNNs, GNNs~\cite{scarselli2008graph, zhou2020graph} extend the idea of localized information aggregation to data defined on irregular domains that can be represented as graphs.
From this perspective, image data can be viewed as a special case of graph-structured data, where pixels correspond to nodes connected by a fixed and regular adjacency pattern.
Graph-structured data are prevalent in real-world applications, including social networks, transportation systems, and scientific computing.

Despite differences in data regularity, CNNs and GNNs share a common underlying principle: both learn representations by aggregating information from local neighborhoods.
In CNNs, this is achieved through convolution operations on structured grids, where the relative ordering of neighbors is fixed and explicitly considered.
In contrast, GNNs generalize this idea to irregular domains via a message-passing mechanism, where each node updates its representation by aggregating information from its neighbors in a permutation-invariant manner, \ie, the aggregation does not depend on the ordering of neighbors.

Formally, a graph is denoted by $\mathcal{G} = (\mathcal{V}, \mathcal{E})$, where $\mathcal{V}$ and $\mathcal{E}$ represent the sets of nodes and edges, respectively.
Each node $n \in \mathcal{V}$ is associated with an initial feature vector $\mathbf{h}_n^{(0)}$.
A GNN updates node representations through a sequence of message-passing layers.
At the $k$-th layer, the representation of node $n$ is updated as
\[
    \mathbf{h}_n^{(k+1)} = \operatorname{Update}\left(
    \mathbf{h}_n^{(k)},
    \operatorname{Aggregate}\left(\left\{ \mathbf{h}_j^{(k)}: j \in \operatorname{Neighbors}(n) \right\}\right)
    \right),
\]
where $\operatorname{Neighbors}(n)$ denotes the set of neighboring nodes of node $n$.
The $\operatorname{Aggregate}(\cdot)$ function aggregates features from neighboring nodes, for example via summation, mean, or max pooling, and the $\operatorname{Update}(\cdot, \cdot)$ function combines the current node representation with the aggregated message, typically using a learnable function such as an FNN, which introduces nonlinearity into the model.
By stacking multiple layers, each node incorporates information from higher-order neighborhoods, enabling the model to capture both local and global structural patterns.
Depending on the task, GNNs can perform node-, edge-, or graph-level predictions, where the inputs and outputs correspond to nodes, edges, or entire graphs, respectively~\cite{zhou2020graph}.
Various extensions introduce specialized aggregation and update schemes to enhance expressive power or adapt to different graph types.
For instance, graph convolutional networks (GCNs)~\cite{kipf2017semisupervised} generalize convolution operations to graph-structured data.
Other widely used architectures include graph sample and aggregate (GraphSAGE)~\cite{hamilton2017inductive}, graph attention networks (GATs)~\cite{velivckovic2018graph}, and graph isomorphism networks (GINs)~\cite{xu2019how}.

Analogous to CNNs, the application of GNNs to linear systems $A\mathbf{x}=\mathbf{b}$ can also be viewed from two perspectives:
\begin{enumerate}
    \item 
          A graph can be represented by its adjacency matrix, and conversely, a matrix can be expressed as a directed graph.
          As illustrated in~\Cref{fig:prelims-struct_unstruct_data}(b), the coefficient matrix $A$ is represented as a graph and processed by a GNN.
          In particular, the correspondence between matrix sparsity and graph structure is widely exploited when approximating sparse matrices using graph-based representations~\cite{chen2022graph}.
          For example, in~\cite{hausner2024neural}, GNNs are employed to construct sparse incomplete Cholesky factor matrices.
    \item 
          GNNs are well-suited for modeling solutions arising from the discretization of PDEs on irregular grids, where the resulting solution data $\mathbf{x}$ are inherently connected irregularly, as shown in~\Cref{fig:prelims-struct_unstruct_data}(d).
\end{enumerate}
We refer the reader to~\cite{saad_iterative_2003, moore2025graph} for further discussion of GNNs and linear algebra.

\paragraph{\textbf{Recurrent neural networks}}

RNNs~\cite{bengio1994learning, medsker1999recurrent} are specifically designed to handle sequential data, such as time series and natural language.
RNNs employ recurrent connections to capture temporal dependencies, in which information from previous time steps is propagated through the hidden state.

\begin{figure}[!htbp]
    \centering
    \includegraphics[width=0.6\textwidth]{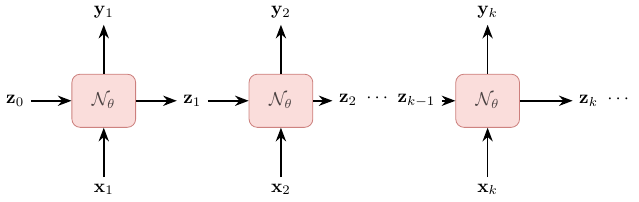}
    \caption{
        Illustration of a recurrent neural network (RNN).
        At each step $k$, the shared network $\mathcal{N}_\theta$ updates the hidden state $\mathbf{z}_k$ from the previous hidden state $\mathbf{z}_{k-1}$ and the current input $\mathbf{x}_k$, and produces the output $\mathbf{y}_k$.
        Here, $\mathcal{N}_\theta$ represents the mappings $g_1$ and $g_2$ in~\Cref{equ:rnn}.
    }
    \label{fig:prelims-ml-rnns}
\end{figure}

In RNNs, a hidden variable, denoted by $\mathbf{z}$, is introduced to encode information from previous steps.
Given a sequence of inputs $\{\mathbf{x}_k\}_{k=1}^T$, the hidden state is updated recursively.
At each step, the hidden state and output are given by
\begin{equation}\label{equ:rnn}
    \mathbf{z}_k = g_1(\mathbf{z}_{k-1}, \mathbf{x}_k), \quad
    \mathbf{y}_k = g_2(\mathbf{z}_k),
\end{equation}
where $\mathbf{z}_k$ denotes the hidden state at time step $k$, and the initial hidden state $\mathbf{z}_0$ is typically initialized to zero or treated as a learnable parameter.
The function $g_1(\cdot, \cdot)$ is a nonlinear mapping that updates the hidden state $\mathbf{z}_k$ based on the previous state $\mathbf{z}_{k-1}$ and the current input $\mathbf{x}_k$, and $g_2(\cdot)$ maps the hidden state $\mathbf{z}_k$ to the output $\mathbf{y}_k$.
Both $g_1(\cdot, \cdot)$ and $g_2(\cdot)$ are typically parameterized by learnable weights.
The overall architecture of the RNN is illustrated in~\Cref{fig:prelims-ml-rnns}.

The recurrent structure of RNNs naturally resembles an iterative process.
Hence, RNN architectures have been employed to learn update functions of iterative methods~\cite{stanziola2021helmholtz}.
However, standard RNNs may struggle to capture long-range dependencies due to the vanishing gradient problem~\cite{bengio1994learning}.
To address this issue, more advanced architectures, such as long short-term memory (LSTM) networks~\cite{hochreiter1997long}, gated recurrent units (GRUs)~\cite{cho2014properties}, and transformers~\cite{vaswani2017attention}, have been developed to improve stability and capture long-range dependencies.

\paragraph{\textbf{Meta-learning}}

Meta-learning, also known as ``learning to learn''~\cite{thrun1998learning}, is a framework that aims to design models and algorithms capable of rapidly adapting to new tasks.
Rather than learning a model for a single fixed problem, meta-learning focuses on improving the learning process itself~\cite{hospedales2022metalearning}.
Typical applications include few-shot learning~\cite{finn2017modelagnostic} and neural architecture search (NAS)~\cite{liu2019darts}.
We refer the reader to~\cite{hospedales2022metalearning} for a comprehensive overview.
In the context of ML-enhanced iterative methods, meta-learning has been employed to select parameters for iterative solvers~\cite{arisaka2023principled}.

\section{Machine learning-enhanced linear iterative methods}\label{sec:linear}

This section reviews ML-enhanced iterative methods (\emph{HIMs}) for linear systems of equations, organized into four categories: \emph{core components learning}~(\Cref{sec:linear-core_components}), \emph{preconditioner learning}~(\Cref{sec:linear-precond}), \emph{enhancement of subspace correction methods}~(\Cref{sec:linear-mgdd}), and \emph{automated solver selection}~(\Cref{sec:linear-solver_selection}).
In each subsection, we first provide the relevant preliminaries on classical methods and discuss the motivation for incorporating ML.
Readers familiar with the classical methods may proceed directly to the corresponding review and refer back to the preliminaries for notation and definitions as needed.
We then review and compare the existing literature from both algorithmic and numerical perspectives, highlighting the key ideas, benefits, and limitations of the proposed approaches.
Throughout the discussion of each topic, we also identify open challenges and potential directions for future research.

\subsection{Core components learning}\label{sec:linear-core_components}

As shown in~\Cref{equ:prelims-general-iter}, a general iterative method starts from an initial guess and generates a sequence of approximations by repeatedly applying an update rule.
Such an iterative method is determined by three \emph{core components}: the initial guess, the update function, and the associated algorithmic parameters.
These components provide natural interfaces through which ML can be incorporated into a classical iterative solver.
In~\Cref{sec:linear-core_components-prelims}, we first recap the preliminaries on classical stationary iterative methods and Krylov subspace methods, and then discuss the motivation for \emph{core components learning}.
We then review three directions: initial guess learning, update function learning, and parameter learning in~\Cref{sec:linear-core_components-initial,sec:linear-core_components-update,sec:linear-core_components-param}, respectively.

\subsubsection{Preliminaries and motivation}\label{sec:linear-core_components-prelims}

For the linear system $A\mathbf{x}=\mathbf{b}$, when $A$ is nonsingular, this system has a unique solution, \ie, $\mathbf{x}_\text{exact} = A^{-1}\mathbf{b}$.
Direct methods obtain the solution in a finite number of steps in exact arithmetic, whereas iterative methods generate successive approximations and are typically terminated when a prescribed tolerance is reached~\cite[p.~58]{quarteroni2007numerical}.
Certain iterative methods, including some Krylov subspace methods, also terminate after a finite number of iterations in exact arithmetic.
For dense matrices, standard dense direct factorizations require $\mathcal{O}(N^3)$ arithmetic operations and $\mathcal{O}(N^2)$ storage.
For sparse matrices, the computational and memory costs depend strongly on the sparsity structure and ordering, since fill-in may occur in the factors.
For large-scale sparse systems, factorization costs and memory growth can therefore become substantial, motivating the use of iterative solvers.
Iterative methods for linear systems are commonly classified into two categories: stationary and non-stationary methods.
An iterative method is called stationary if its update function and algorithmic parameters are independent of the iteration index; otherwise, it is non-stationary, as in Krylov subspace methods.
Although iterative methods are generally viewed as generating convergent sequences of approximations, Krylov subspace methods may converge in a finite number of steps in exact arithmetic.
Both classes are briefly reviewed below.

\paragraph{Stationary iterative methods}\label{sec:prelims-linear-stationary}
To illustrate the basic structure of stationary iterative methods, we begin with the Richardson iteration~\cite{lewisfry1911approximate},
\begin{equation*}
    \begin{aligned}
        \mathbf{x}_{k+1}
         & = \mathbf{x}_k + \alpha_k \mathbf{r}_{k}              \\
         & = (I - \alpha_k A)\mathbf{x}_k + \alpha_k \mathbf{b},
    \end{aligned}
\end{equation*}
where $\mathbf{r}_{k}=\mathbf{b} - A\mathbf{x}_k$ denotes the residual at the $k$-th iteration and $\alpha_k \in \mathbb{R}$ is the relaxation parameter.
The method is stationary when $\alpha_k=\alpha$ is fixed throughout the iterations; otherwise, it is referred to as the non-stationary Richardson iteration.
For the stationary Richardson iteration, the error satisfies $\mathbf{e}_{k+1}=(I-\alpha A)\mathbf{e}_k$, and the method converges for every initial guess if and only if $\rho(I-\alpha A)<1$.
The convergence of Richardson iteration can often be significantly improved through preconditioning, leading to the preconditioned Richardson iteration
\begin{equation}\label{equ:prelims-preconditioned-richardson}
    \begin{aligned}
        \mathbf{x}_{k+1}
         & = \mathbf{x}_k + \alpha_k M^{-1}(\mathbf{b}-A\mathbf{x}_k)      \\
         & = (I-\alpha_k M^{-1}A)\mathbf{x}_k + \alpha_k M^{-1}\mathbf{b},
    \end{aligned}
\end{equation}
where $M\in\mathbb{R}^{N\times N}$ is a nonsingular matrix referred to as the preconditioner.
When $\alpha_k=1$, if $M^{-1}A$ is close to the identity matrix, then the spectral radius of the iteration matrix $B= I- M^{-1}A$
is typically small, resulting in faster convergence.

Many classical stationary iterative methods, including Jacobi~\cite{jacobi_ueber_1845}, Gauss-Seidel (GS)~\cite{gauss_letter_1903}, and successive over-relaxation (SOR) methods~\cite{young1954iterative}, can be derived from matrix splittings of the form $A=M-N$, where $M$ is nonsingular and $N=M-A$.
This decomposition yields the iteration
\begin{equation*}
    \begin{aligned}
        \mathbf{x}_{k+1}
         & = M^{-1}N\mathbf{x}_k + M^{-1}\mathbf{b}      \\
         & = (I-M^{-1}A)\mathbf{x}_k + M^{-1}\mathbf{b},
    \end{aligned}
\end{equation*}
which coincides with~\Cref{equ:prelims-preconditioned-richardson} when $\alpha_k=1$.
Therefore, these methods can be interpreted as special cases of preconditioned Richardson iteration.
A classical splitting is obtained from $A=D+L+U$, where $D$ is the diagonal part of $A$, $L$ is the strictly lower triangular part, and $U$ is the strictly upper triangular part.
Based on this decomposition, several classical stationary methods are obtained:
\begin{enumerate}
    \item Jacobi method: $M = D$ and $N = -(L + U)$.
    \item Gauss-Seidel method: $M = D + L$ and $N = - U$.
    \item SOR method: $M = \frac{1}{\omega} D + L$ and $N = (\frac{1}{\omega}-1)D - U$, where $\omega \in (0, 2)$ is the relaxation parameter.
\end{enumerate}

These stationary iterative methods are also commonly referred to as relaxation methods.
In practice, these methods often converge slowly because the spectral radius of the iteration matrix is typically close to one.
Nevertheless, they can be used as a component in more sophisticated methods.
For example, they are often used as smoothers in multigrid methods and as preconditioners for Krylov subspace methods.

\paragraph{Krylov subspace methods}
Krylov subspace methods are a class of non-stationary iterative methods for solving large sparse linear systems, as they do not rely on a fixed iteration matrix.
The main idea is to approximate the exact solution in a sequence of low-dimensional subspaces, known as Krylov subspaces.
Given an initial guess $\mathbf{x}_0$ and the corresponding residual $\mathbf{r}_{0} = \mathbf{b} - A \mathbf{x}_0$, Krylov subspace methods find an approximate solution $\mathbf{x}_k \in \mathbf{x}_0 + \mathcal{K}_k(A, \mathbf{r}_{0})$, where $\mathcal{K}_k(A, \mathbf{r}_{0}) = \text{span}\left\{\mathbf{r}_{0}, A\mathbf{r}_{0}, A^2\mathbf{r}_{0}, \cdots, A^{k-1}\mathbf{r}_{0}\right\}$, $k=1, 2, \cdots$ denotes the $k$-th Krylov subspace associated with $A$ and $\mathbf{r}_0$~\cite{vuik2018krylov}.

A variety of Krylov subspace methods have been developed for coefficient matrices with different structural properties.
The conjugate gradient (CG) method~\cite{Hestenes1952MethodsOC,lanczos_solution_1952} is designed for symmetric positive definite (SPD) systems, whereas the minimal residual (MINRES) method~\cite{paige1975solution} applies to symmetric, possibly indefinite systems.
For general nonsymmetric systems, widely used methods include the generalized minimal residual (GMRES) method~\cite{saad_gmres_1986}, the biconjugate gradient stabilized (BiCGSTAB) method~\cite{van_der_vorst_bi-cgstab_1992}, the induced dimension reduction (IDR) method~\cite{wesseling_numerical_1980}, and the IDR($s$) method~\cite{sonneveld_idr_2009}.
In the following, we briefly review the CG and GMRES methods.

For any vector $\mathbf{x} \in \mathbb{R}^N$, the $A$-norm is defined by $\|\mathbf{x}\|_A^2 = \mathbf{x}^T A\mathbf{x}$.
Given an SPD matrix $A$, the CG method computes $\mathbf{x}_k \in \mathbf{x}_0 + \mathcal{K}_k(A, \mathbf{r}_{0})$ such that the error $\|\mathbf{x}_k - \mathbf{x}_\text{exact}\|_A$ is minimized.
The CG iterate $\mathbf{x}_k$ satisfies the following inequality:
\begin{equation*}
    \| \mathbf{x}_k - \mathbf{x}_\text{exact} \|_A \leq 2 \left( \frac{\sqrt{\kappa_2(A)} - 1}{\sqrt{\kappa_2(A)} + 1} \right)^k \| \mathbf{x}_0 - \mathbf{x}_\text{exact} \|_A,
\end{equation*}
where $\kappa_2(A) = \frac{\lambda_{\max}(A)}{\lambda_{\min}(A)}$ denotes the spectral condition number of $A$, representing the ratio of its largest to smallest eigenvalues.
The condition number thus provides a convergence bound for CG. A more precise assessment of its convergence speed requires consideration of the full eigenvalue distribution, including spectral clustering, and the initial error.
For a nonsymmetric matrix $A$, the GMRES method constructs an orthonormal basis $\left\{\mathbf{v}_1, \mathbf{v}_2, \cdots, \mathbf{v}_k\right\}$ of the $k$-th Krylov subspace $\mathcal{K}_k(A, \mathbf{r}_{0})$ via the Arnoldi iteration.
It then computes $\mathbf{x}_k \in \mathbf{x}_0 + \mathcal{K}_k(A, \mathbf{r}_0) = \mathbf{x}_0 + \text{span} \left\{\mathbf{v}_1, \mathbf{v}_2, \cdots, \mathbf{v}_k\right\}$ such that the residual norm $\|\mathbf{b} - A \mathbf{x}_k\|_2$ is minimized.
We refer readers to~\cite{vuik2018krylov, greenbaum_iterative_1997, golub2013matrix} for more details on CG and GMRES methods.

\paragraph{\textbf{Motivation for core components learning}}

\begin{figure}[!htbp] 
    \centering
    \includegraphics[width=1.0\textwidth]{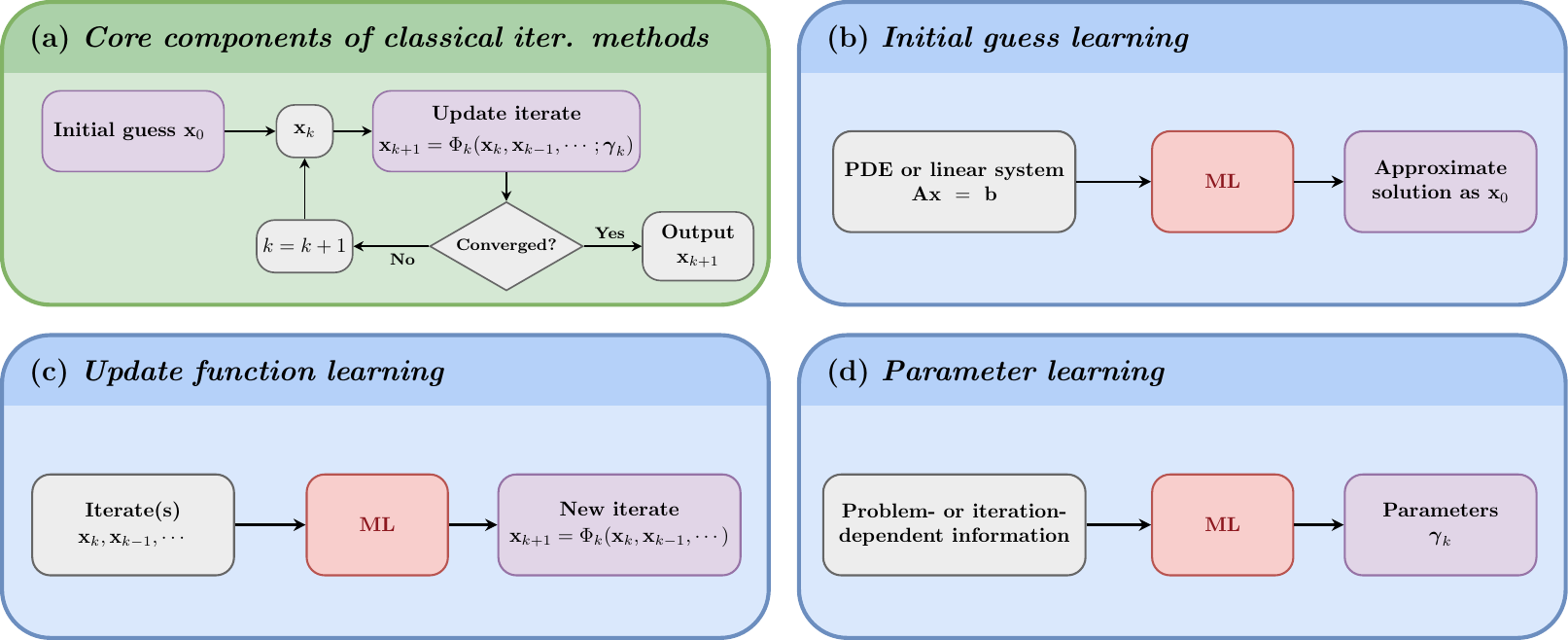}
    \caption{
        Overview of \emph{core components learning}~(\Cref{sec:linear-core_components}).
        \textbf{(a)} Core components of a general iterative method: the initial guess $\mathbf{x}_0$, the update function $\Phi_k$, and parameters $\gamma_k$.
        \textbf{(b)} \emph{Initial guess learning}: ML predicts an approximate solution that serves as the initial guess.
        \textbf{(c)} \emph{Update function learning}: ML replaces or augments the classical update function.
        \textbf{(d)} \emph{Parameter learning}: ML learns algorithmic parameters from problem- or iteration-dependent information.
    }
    \label{fig:core_components}
\end{figure}

Despite their substantial differences, stationary and Krylov subspace methods can be described within a unified three-component framework, as summarized in~\Cref{tab:linear-core_components-prelims}.
A schematic illustration of this framework is provided in~\Cref{fig:core_components}(a).
This unified perspective naturally motivates three corresponding directions for ML enhancement:
\begin{enumerate}
    \item 
          The \textbf{initial guess} $\mathbf{x}_0$ determines the initial error and residual and can therefore strongly affect the convergence of stationary iterative methods. Unlike stationary iterations, which repeatedly apply a fixed error-propagation operator, Krylov methods build approximations in expanding subspaces and can exploit spectral clustering.
          In Krylov subspace methods, it additionally determines the Krylov subspace $\mathcal{K}_k(A,\mathbf{r}_0)$ generated during the iteration.
          A more accurate initial guess may therefore reduce the number of iterations required to reach a prescribed tolerance.
          \emph{Initial guess learning} aims to obtain such initial guesses by using ML models to predict approximate solutions.
    \item 
          The \textbf{update function} determines how information generated during the iteration, such as previous iterates and residuals, is used to construct the next iterate.
          Classical iterative methods typically possess clear mathematical interpretations and, in many cases, theoretical convergence properties.
          \emph{Update function learning} augments the iterative update with ML techniques to accelerate convergence or reduce computational cost.
    \item 
          \textbf{Parameters} such as relaxation parameters, GMRES restart lengths, and stopping tolerances also influence the practical performance of iterative solvers.
          \emph{Parameter learning} selects such parameters based on features of the linear system or information generated during the iteration.
          Similar ideas can also be applied to parameters associated with classical preconditioners and subspace correction methods.
\end{enumerate}

\begin{table}[!ht]
    \renewcommand{\arraystretch}{1.3}
    \centering
    \caption{
        Core components of stationary iteration and Krylov subspace methods.
    }
    \label{tab:linear-core_components-prelims}
    \resizebox{\linewidth}{!}{%
        \begin{minipage}{\linewidth}
            \centering
            \begin{tabular}{
                    >{\raggedright\arraybackslash}p{2.5cm}
                    >{\raggedright\arraybackslash}p{3.2cm}
                    >{\raggedright\arraybackslash}p{4.0cm}
                    >{\raggedright\arraybackslash}p{6.0cm}
                }
                \toprule
                \textbf{method class} & \textbf{initial guess}                                                      & \textbf{update function} & \textbf{parameters}
                \\
                \midrule
                \textbf{stationary methods}
                                      & determines the initial error
                                      & $\mathbf{x}_{k+1} = \mathbf{x}_{k} + M^{-1}(\mathbf{b} - A\mathbf{x}_{k})$
                                      & relaxation/acceleration parameters, stopping tolerance
                \\
                \midrule
                \textbf{Krylov subspace methods}
                                      & determines $\mathbf{r}_0$ and $\mathcal{K}_k(A,\mathbf{r}_0)$
                                      & \eg, CG for SPD systems, GMRES for non-SPD systems
                                      & restart length, truncation or recycling parameters, and stopping tolerances
                \\
                \bottomrule
            \end{tabular}
        \end{minipage}
    }
\end{table}

\subsubsection{Initial guess learning}\label{sec:linear-core_components-initial}

Common choices of an initial guess include a zero vector, a random vector, or the solution from the previous time step in time-dependent problems.
In full multigrid methods, the initial guess on the fine grid is typically obtained by first solving a coarse-grid problem and then interpolating the solution back to the fine grid~\cite{trottenberg2001multigrid}.
We refer the reader to~\cite{tromeur-dervout2006choicea, ye2020improving, austin2021initial} for further discussion of classical choices of initial guesses.
\emph{Initial guess learning} employs an ML model to generate an approximate solution that serves as the initial guess for a classical iterative method, as illustrated in~\Cref{fig:core_components}(b).
Although the convergence robustness of linear iterative methods is generally less sensitive to the choice of the initial guess, there are two main motivations for \emph{initial guess learning}:
\begin{enumerate}
    \item From the perspective of classical iterative solvers, the goal is to use ML to predict an improved initial guess that reduces the initial error and thereby accelerates convergence.
          This is particularly attractive when repeatedly solving a sequence of related linear systems, such as parametric systems, where reducing the iteration count can lead to substantial savings in overall computational cost.
    \item From the perspective of ML, solutions predicted by a learned surrogate model often have limited accuracy, partly due to generalization issues.
          A classical iterative solver can subsequently be employed to refine the predicted solution until a prescribed error tolerance is reached.
\end{enumerate}
In this subsection, we primarily interpret \emph{initial guess learning} from the first perspective, that is, that of classical iterative methods.

When neural networks are employed for this purpose, the resulting hybrid solver can be understood from the perspective of spectral complementarity.
Due to their spectral bias~\cite{rahaman2019spectral, xu2020frequency, xu2025overview, xu2025understanding}, neural networks tend to capture the low-frequency components of the solution more rapidly, whereas classical iterative methods can subsequently reduce the remaining error components, which are often dominated by high-frequency modes.
This complementary behavior provides a natural basis for combining the two approaches.
In the following, we review several representative methods, whose algorithmic differences are summarized in~\Cref{tab:linear-core_components-initial-algorithms}.

\begin{table}[!ht]
    \renewcommand{\arraystretch}{1.3}
    \centering
    \caption{
        Comparison of existing approaches for \emph{initial guess learning}~(\Cref{sec:linear-core_components-initial}).
        The last column summarizes our interpretation of the main motivation behind each approach.
    }
    \label{tab:linear-core_components-initial-algorithms}
    \resizebox{\linewidth}{!}{%
        \begin{minipage}{\linewidth}
            \centering
            \begin{tabular}{
                    lllll
                }
                \toprule
                \textbf{ref.}
                 & \textbf{tested problem}
                 & \textbf{ML}
                 & \textbf{iter. method}
                 & \textbf{main motivation}
                \\
                \midrule
                \cite{xu2020frequency}
                 & Poisson eq.
                 & FNN
                 & Jacobi
                 & from ML: spectral bias of neural networks
                \\
                \cite{ozbay2021poisson}
                 & Poisson eq.
                 & CNN
                 & multigrid/BiCGSTAB
                 & from ML: improve CNN predictions
                \\
                \cite{illarramendi2020hybrid}
                 & Poisson eq. from incompressible flow
                 & CNN
                 & Jacobi
                 & from ML: generalization issue of the CNN
                \\
                \cite{luna2021accelerating}
                 & Poisson eq.
                 & CNN
                 & GMRES
                 & from classical iter.: reduce iterations
                \\
                \cite{vigot2024solving}
                 & Poisson eq. from plasma simulations
                 & GraphSAGE
                 & GMRES
                 & from classical iter.: reduce iterations
                \\
                \cite{nikolopoulos2024aienhanced}
                 & Parametrized elasticity \& Biot problems
                 & CAE, FNN
                 & POD-2G
                 & from classical iter.: reduce iterations
                \\
                \cite{eshaghi2026nows}
                 & Poisson, Darcy, elasticity, Burgers, smoke plume
                 & neural operator
                 & CG/PCG
                 & from classical iter.: reduce iterations
                \\
                \cite{wang2026pretrainingfinetuning}
                 & 3D elasticity homogenization
                 & FNO
                 & PCG
                 & from classical iter.: reduce iterations
                \\
                \cite{wang2026pretrain}
                 & Linear elasticity, nonlinear hyperelasticity, homogenization
                 & Transolver
                 & CG/PCG
                 & from classical iter.: reduce iterations
                \\
                \bottomrule
            \end{tabular}
        \end{minipage}
    }
\end{table}

Building on the spectral bias of neural networks, \citet{xu2020frequency} proposed a hybrid solver for PDEs in which approximate solutions generated by FNNs are refined by the Jacobi method to efficiently reduce the remaining high-frequency errors.
\citet{ozbay2021poisson} propose a CNN-based solver for the Poisson equation on 2D Cartesian grids, designed to handle varying grid resolutions, right-hand sides, and boundary conditions without retraining.
To incorporate boundary conditions, the Poisson problem is decomposed into a homogeneous Poisson problem and several Laplace subproblems with inhomogeneous boundary conditions, which are approximated by two CNN-based subnetworks.
The networks are trained in a supervised manner by minimizing the solution error with respect to reference solutions.
The predicted solution is further employed as an initial guess for multigrid and BiCGSTAB methods to accelerate convergence.

\citet{tompson2017accelerating} introduce CNNs for accelerating pressure projection in incompressible flow simulations.
Building on this idea, \citet{illarramendi2020hybrid} propose a hybrid method combining a CNN and the Jacobi method to improve robustness and accuracy.
The CNN takes the velocity divergence, which forms the right-hand side vector $\mathbf{b}$ of the discretized Poisson system, along with geometric information about boundaries and obstacles, and predicts the pressure field.
The CNN prediction is then refined by Jacobi iterations to satisfy a prescribed tolerance.

\citet{luna2021accelerating} consider a sequence of linear systems arising from discretized 2D Poisson equations, where the coefficient matrix $A$ is fixed and the right-hand side $\mathbf{b}^{(n)}$ varies with the simulation time step.
A CNN is trained to predict an approximate solution $\mathbf{x}_0^{(n)}$ from $\mathbf{b}^{(n)}$, which is used as an initial guess for GMRES to accelerate convergence.
To avoid expensive offline data generation, the CNN is trained online using samples collected during the simulation.

A similar problem setting is considered in \cite{vigot2024solving}, where Poisson equations on unstructured grids in plasma simulations also give rise to a sequence of linear systems with a fixed $A$ and time-varying $\mathbf{b}^{(n)}$.
Rather than directly predicting the solution from $\mathbf{b}^{(n)}$, a GraphSAGE model~\cite{hamilton2017inductive}, denoted by $\mathcal{N}$, is trained to predict the solution increment $\Delta \mathbf{x}^{(n)}=\mathcal{N}(\mathbf{x}^{(n-1)},\mathbf{b}^{(n)})$.
The predicted solution $\mathbf{x}^{(n)}_0=\mathbf{x}^{(n-1)}+\Delta \mathbf{x}^{(n)}$ is then used as the initial guess for GMRES.
The model is trained in a supervised manner using reference solution increments.

\citet{nikolopoulos2024aienhanced} propose an ML-enhanced framework for accelerating iterative solvers for large-scale parameterized linear systems.
A convolutional autoencoder (CAE) and an FNN are trained to predict approximate solutions from the problem parameters via a low-dimensional latent representation.
For a new parameter instance, the predicted solution is used as an initial guess for a two-grid iterative method that combines multigrid with POD-based transfer operators.

\citet{eshaghi2026nows} propose a hybrid framework, called neural operator warm starts (NOWS), for accelerating iterative solvers for linear systems arising from the discretization of parametric PDEs.
In this approach, neural operators learn mappings from PDE parameters to solutions.
The predicted solution is then used as an initial guess for a Krylov subspace method to accelerate its convergence.
In their experiments, neural operators including FNO, variational physics-informed neural operator (VINO)~\cite{eshaghi2025variational} and multi-head neural operator (MHNO)~\cite{eshaghi2026multihead} are exploited.
This method has been extended to the homogenization method in~\cite{wang2026pretrainingfinetuning}.
In~\cite{wang2026pretrain}, the authors propose a ``pretrained and warm-start'' framework in which the transformer-based PDE solvers called Transolver~\cite{wu2024transolver} are employed to generate initial guesses for classical iterative solvers.

\paragraph{\textbf{Discussion}}

In \emph{initial guess learning}, ML models are employed to generate improved initial guesses for classical iterative methods, although the approach can also be viewed from an ML perspective.
By providing improved initial guesses, ML models reduce the initial residual $\mathbf{r}_0$, which can accelerate the convergence of iterative methods.

We have the following general remarks on the current status and future directions of \emph{initial guess learning}:
\begin{enumerate}
    \item \textbf{Low-risk hybridization:}
          The resulting hybrid iterative method is conceptually non-intrusive, since the ML model is called only once at the beginning of the solution process, and the learned initial guess does not require modification of the underlying classical iterative solver.
          This feature makes the approach particularly attractive for integration into existing solvers and software for computational science and engineering, such as large-scale fluid simulations.
    \item \textbf{Iteration reduction and convergence rate:}
          As shown in~\cite{eshaghi2026nows}, compared with a zero initial guess, a learned initial guess can reduce the number of iterations and therefore the wall-clock time required to reach the same convergence tolerance.
          However, it generally does not change the asymptotic convergence rate of the underlying iterative method, which is mainly governed by the spectral properties of the coefficient matrix and can be fundamentally improved, for example, by preconditioning.
    \item \textbf{Is improvement guaranteed?}
          Compared with classical choices of initial guesses, it remains unclear whether a learned initial guess necessarily reduces the iteration number or wall-clock time.
          Although iteration reductions have been widely observed empirically, general theoretical guarantees are still lacking.
          Moreover, in most existing studies, the ML training objective is not directly tied to the actual performance of the iterative solver.
          A general training loss for learning initial guesses and solver parameters is proposed in~\cite{arisaka2023principled} and will be discussed in~\Cref{sec:linear-core_components-param}.
    \item \textbf{Effect of the target tolerance}:
          As reported in~\cite{wang2026pretrainingfinetuning}, the benefit of a learned initial guess also depends on the prescribed convergence tolerance.
          A substantial reduction in iterations may be observed for relatively moderate tolerances.
          However, the relative gain can diminish as the tolerance becomes stricter, since in this case the influence of the initial guess on the iteration number becomes less significant.
    \item \textbf{Offline training cost:}
          \begin{enumerate}
              \item 
                    Training ML models usually introduces an offline computational cost.
                    Therefore, the overall benefit depends on whether this cost can be amortized over many online solves, as in multi-query problems or sequences of related linear systems.
                    Methods that can be trained on coarse-grid discretizations and transferred to finer grids, such as certain neural-operator-based approaches, may help reduce this burden.
                    At the same time, sequences of linear systems have been extensively studied using classical techniques, particularly Krylov subspace recycling methods.
                    A systematic comparison between \emph{initial guess learning} strategies and such classical approaches, in terms of both computational cost and convergence improvement, deserves further investigation.
              \item 
                    The benefit of a learned initial guess strongly depends on its quality, and more accurate initial guesses are generally expected to provide larger reductions in subsequent classical iterations, as also observed in~\cite{wang2026pretrainingfinetuning}.
                    However, obtaining a more accurate prediction typically requires greater model capacity, training effort, or inference cost.
                    This leads to a trade-off between the cost of producing the initial guess and the computational savings achieved by the subsequent iterative solver.
          \end{enumerate}
    \item \textbf{Extensions}:
          \begin{enumerate}
              \item 
                    \emph{Initial guess learning} is a conceptually general strategy that is not restricted to any particular ML architecture.
                    For PDE problems, any ML-based approach capable of producing a compatible approximate solution can potentially be used to provide an initial guess for a subsequent classical iterative solver.
              \item 
                    Most existing studies focus on PDE problems, especially Poisson equations.
                    Nevertheless, \emph{initial guess learning} is not inherently limited to PDE solvers and can be extended more broadly to iterative algorithms that benefit from high-quality initializations.
                    For instance, learned warm starts for optimization problems have been investigated in~\cite{sambharya2024learning}.
          \end{enumerate}
\end{enumerate}

\subsubsection{Update function learning}\label{sec:linear-core_components-update}

\emph{Update function learning} refers to the use of an ML model, denoted by $\mathcal{N}$, to learn an iterative update mapping, such as
\[
    \mathbf{x}_{k+1} = \mathcal{N}(\mathbf{x}_k, \mathbf{x}_{k-1}, \cdots),
\]
or an incremental update of the form
\[
    \mathbf{x}_{k+1} = \mathbf{x}_k + \mathcal{N}(\mathbf{x}_k, \mathbf{x}_{k-1}, \cdots).
\]
ML models can also be applied to augment an existing solver, yielding an update rule such as
\[
    \mathbf{x}_{k+1} = \Phi_{0}(\mathbf{x}_k, \mathbf{x}_{k-1}, \cdots) + \mathcal{N}(\mathbf{x}_k, \mathbf{x}_{k-1}, \cdots),
\]
where $\Phi_{0}$ represents an existing solver, also referred to as the base solver.
The schematic illustration of \emph{update function learning} is shown in~\Cref{fig:core_components}(c).

\citet{hsieh2019learning} propose a learned iterative PDE solver with convergence guarantees by combining neural networks with classical iterative methods.
Specifically, the neural network is trained to learn a correction to the update step produced by a classical solver.
They focus on PDEs with Dirichlet boundary conditions.
A linear iterative solver $\Phi_0$ is defined as
\[
    \mathbf{x}_{k+1} = \Phi_0 (\mathbf{x}_k) = T_c\mathbf{x}_k + \mathbf{c},
\]
where $T_c$ is a constant matrix and $\mathbf{c}$ is a vector.
Given an existing linear iterative solver (base solver) $\Phi_0$, such as the Jacobi method, the learned iteration scheme is defined as
\begin{equation}\label{equ:update_learning-enhanced}
    \mathbf{x}_{k+1} = \Phi(\mathbf{x}_k) = \Phi_0(\mathbf{x}_k) + G\mathcal{N}_\text{linear}(\Phi_0(\mathbf{x}_k) - \mathbf{x}_k),
\end{equation}
where the operator $G$ restricts the learned correction to the interior, thereby preserving the boundary conditions.
$\mathcal{N}_{\text{linear}}$ is a learned linear operator parameterized by a CNN without activation functions or bias terms.
The linear operator $\mathcal{N}_{\mathrm{linear}}$ takes the update produced by the base solver, \ie, $\Phi_0(\mathbf{x}_k)-\mathbf{x}_k$, as input and predicts an additional correction.
Since the learned correction vanishes at any fixed point of the base solver, the correct solution remains a fixed point of the learned iteration scheme.
The networks are trained by minimizing the solution error after $k$ learned iterations:
\[
    \mathcal{L} = \left\| \Phi^k(\mathbf{x}_0)-\mathbf{x}^* \right\|_2^2,
\]
where $\mathbf{x}_0$ is a randomly sampled initial guess, $\Phi^{k}$ denotes $k$ successive applications of the learned iteration scheme $\Phi$, and $\mathbf{x}^*$ is the reference solution.

\citet{stanziola2021helmholtz} develop a learned iterative solver to efficiently solve 2D heterogeneous Helmholtz equations.
The solver is implemented using a recurrent modified U-Net.
Let $\mathbf{r}_k = \mathbf{b} - A\mathbf{x}_k$ denote the residual at the $k$-th iteration.
The ML model $\mathcal{N}$ and the resulting iterative scheme are given by
\begin{equation}\label{equ:update_learning-fully_learned}
    (\mathbf{d}_k, \mathbf{h}_{k+1}) = \mathcal{N}(\mathbf{x}_k, \mathbf{r}_k, \mathbf{h}_k), \qquad
    \mathbf{x}_{k+1} = \mathbf{x}_k + \mathbf{d}_k,
\end{equation}
respectively.
Here, the hidden states $\mathbf{h}_k$ and $\mathbf{h}_{k+1}$ are introduced to retain information from previous iterations, and $\mathcal{N}$ is implemented as a modified U-Net with hidden states at multiple scales.
The model is trained in an unsupervised manner by minimizing the residuals over the iterative trajectory,
\[
    \mathcal{L} = \sum_{k=1}^{K_\mathrm{{loss}}} \|\mathbf{r}_k\|_2^2 = \sum_{k=1}^{K_\mathrm{{loss}}} \|\mathbf{b} - A\mathbf{x}_k\|_2^2,
\]
where $K_\mathrm{{loss}}$ denotes the maximum number of iterations considered during training.

To make the iteration more stable and improve the neural network training proposed in~\cite{stanziola2021helmholtz}, \citet[Chapter~5.3.1]{xiang2022solution} introduces a minimum residual Richardson (MRR) scheme.
Specifically, they incorporate a relaxation parameter (step size) $\alpha_k$ into the update function,
\[
    \mathbf{x}_{k+1} = \mathbf{x}_k + \alpha_k \mathbf{d}_k,
\]
where $\mathbf{d}_k$ is the correction direction predicted by the neural network.
The relaxation parameter is chosen to minimize the residual norm along this direction:
\[
    \begin{aligned}
        \alpha_k
         & = \arg\min_{\alpha} \left\| \mathbf{b} - A \mathbf{x}_{k+1} \right\|_2^2                     \\
         & = \arg\min_{\alpha} \left\| \mathbf{r}_k - \alpha \left( A \mathbf{d}_k \right) \right\|_2^2
        = \frac{(A \mathbf{d}_k)^H \mathbf{r}_k}{\left\| A \mathbf{d}_k \right\|_2^2}.
    \end{aligned}
\]
Here, $(\cdot)^H$ denotes the Hermitian transpose.

\paragraph{\textbf{Discussion}}

We make the following general remarks on existing methods and potential future directions for \emph{update function learning}:

\begin{enumerate}
    \item \textbf{Connection to preconditioner learning:}
          All of the methods reviewed above can be interpreted from the perspective of \emph{preconditioner learning} (see~\Cref{sec:linear-precond}), where the learned update effectively corresponds to learning a preconditioner.
          \begin{enumerate}
              \item In~\cite{hsieh2019learning}, when a stationary iterative method is used as the base solver, its update can be written as $\Phi_0(\mathbf{x}_k)=\mathbf{x}_k+M_0^{-1}\mathbf{r}_k$, where $M_0^{-1}$ denotes the corresponding inverse preconditioner.
                    Substituting this expression into~\Cref{equ:update_learning-enhanced} gives $\Phi(\mathbf{x}_k) = \mathbf{x}_k + (I+G\mathcal{N}_\text{linear})M_0^{-1}\mathbf{r}_k$.
                    Hence, the learned iterator can be interpreted as applying the ML-based preconditioner $(I+G\mathcal{N}_\text{linear})M_0^{-1}$ to the Richardson iteration.
              \item In~\cite{stanziola2021helmholtz, xiang2022solution}, omitting the latent state in~\Cref{equ:update_learning-fully_learned} for simplicity, the iteration can be written as $\mathbf{x}_{k+1} = \mathbf{x}_k + \mathcal{N}(\mathbf{x}_k, \mathbf{r}_k)$, which gives $\mathbf{r}_{k+1} = \mathbf{r}_k - A \mathcal{N}(\mathbf{x}_k, \mathbf{r}_k)$.
                    Since the network is trained by minimizing the residuals along the iterative trajectory, $\mathcal{L} = \sum_{k=1}^{K_\mathrm{{loss}}} \|\mathbf{r}_k\|_2^2$, the learned correction is encouraged to reduce the residual.
                    In the ideal case, where the next residual is driven close to zero, $\mathbf{r}_k \approx A \mathcal{N}(\mathbf{x}_k, \mathbf{r}_k)$, and hence $\mathcal{N}(\mathbf{x}_k, \mathbf{r}_k) \approx A^{-1}\mathbf{r}_k$.
                    From this perspective, the network can be viewed as a learned nonlinear preconditioner.
                    Building on this interpretation, \citet[Chapter~5.3.2]{xiang2022solution} further employs the trained network as a preconditioner within flexible Krylov subspace methods; see also~\Cref{sec:linear-precond-direct}.
          \end{enumerate}

          In summary, methods in~\cite{hsieh2019learning, stanziola2021helmholtz, xiang2022solution} can therefore be regarded as Richardson iterations equipped with learned preconditioners, differing primarily in whether the learned preconditioner is linear~\cite{hsieh2019learning} or nonlinear~\cite{stanziola2021helmholtz, xiang2022solution}.

          However, \emph{update function learning} is conceptually more general than \emph{preconditioner learning}, since the learned update may depend not only on the current residual $\mathbf{r}_k$, but also on additional iterative information such as previous iterates or momentum terms.
          Such update strategies are conceptually related to classical acceleration techniques such as Anderson acceleration and momentum-based methods.
    \item \textbf{Training data and efficient offline training:}
          In the papers reviewed above, the training data consist of different initial guesses $\mathbf{x}_0$, and the loss is generally accumulated over $K_\mathrm{{loss}}$ iterations.
          $K_\mathrm{{loss}}$ ranges from $1$ to $20$ in~\cite{hsieh2019learning}, and it is set to $1000$ in~\cite{stanziola2021helmholtz}.
          For large $K_\mathrm{{loss}}$, backpropagation through the full iterative trajectory can be computationally prohibitive.
          To address this, the authors of~\cite{stanziola2021helmholtz} employ a replay buffer and truncated backpropagation through time (TBPTT)~\cite{jaeger2002tutorial}.
          Specifically, intermediate solver states are stored in the replay buffer, and only a short window of iterations (10 steps in their experiments) is unrolled for backpropagation at each training step.
    \item \textbf{Fair comparison against classical solvers:}
          When comparing the learned update functions against classical iterative solvers, the number of iterations alone may not provide a fair measure of computational efficiency.
          Although a learned update scheme may require substantially fewer iterations, each iteration can be more expensive due to the evaluation of the learned model.
          Therefore, comparisons should account for the actual computational cost, for example through the number of floating-point operations, as considered in~\cite{hsieh2019learning}, or through wall-clock time.
          Moreover, the offline training cost should also be taken into account.
    \item \textbf{Generalization}:
          In~\cite{hsieh2019learning}, the linear iterative solver is enhanced by a learned linear operator that preserves the fixed point of the underlying solver.
          By contrast, learned nonlinear update functions, such as those in~\cite{stanziola2021helmholtz, xiang2022solution}, generally lack theoretical convergence guarantees and may encounter stagnation on unseen problems, as reported in~\cite{xiang2022solution}.
\end{enumerate}

\subsubsection{Parameter learning}\label{sec:linear-core_components-param}

For classical iterative methods and their acceleration techniques, performance often depends critically on the choice of parameters.
In practice, these parameters are typically selected heuristically, which can be computationally expensive, particularly when solving a sequence of linear systems.
ML methods provide a framework for automatically determining such parameters based on problem- or iteration-dependent information, as shown in~\Cref{fig:core_components}(d).

\citet{khodak2024learning} propose an online bandit learning approach for selecting the SOR relaxation parameter when solving sequences of linear systems.
The selection of the relaxation parameter is formulated as a bandit task, in which the algorithm selects a parameter and receives the resulting iteration count as feedback.
The objective is to minimize the total number of iterations over the entire sequence.
They consider two settings.
First, for a general sequence of systems ${(A_1,\mathbf{b}_1),\ldots,(A_T,\mathbf{b}_T)}$, Tsallis-INF~\cite{zimmert2019optimal} is used to select the relaxation parameter online, aiming to achieve a total iteration count close to that obtained with the best fixed parameter over the sequence.
Second, for related systems of the form $A_t=A+c_t\mathbf{I}$, the shift $c_t \in \mathbb{R}$ is used as context to learn an instance-dependent mapping $c_t \to \omega_t$.
Since the method relies solely on iteration counts as feedback, it can also be extended to the selection of other solver components, such as those arising in preconditioners~\cite{khodak2025pcgbandit}; see also~\Cref{sec:linear-solver_selection} for more details.

\citet{cui2025neural} focus on solving Helmholtz equations with high wavenumber and heterogeneous media using a learning-enhanced multigrid method; see also~\Cref{sec:linear-mgdd-smoother}.
On coarse grid levels, the Chebyshev semi-iteration is applied to the normal equation to damp high-frequency error components:
\[
    \mathbf{x}_{k+1} = \mathbf{x}_k + \alpha_k A^T(\mathbf{b} - A\mathbf{x}_k).
\]
After $k$ steps, the error satisfies $\mathbf{e}_k = p_k(A^TA)\mathbf{e}_0$, where $p_k(A^TA) = \prod_{j=0}^{k-1}(I-\alpha_j A^TA)$.
The Chebyshev coefficients $\{\alpha_j\}_{j=0}^{k-1}$ are selected such that the corresponding polynomial is minimized over the spectral interval $\left[ \frac{\lambda_{\max}}{\nu}, \lambda_{\max} \right]$, where $\lambda_{\max}$ is the largest eigenvalue of $A^TA$.
The parameter $\nu > 1$ determines the target spectral interval such that $\lambda_{\max}/\nu$ approximates the median eigenvalue of $A^TA$.
Since the appropriate $\nu$ varies across multigrid levels and problem parameters, a CNN is trained to predict $\nu$ from the Helmholtz problem parameters, including the angular frequency and the slowness field.
The CNN is trained in an unsupervised manner through the differentiable multigrid solver by minimizing the relative residual, $\|\mathbf{b}-A\mathbf{x}_k\|_2^2 / \|\mathbf{b}\|_2^2$.

A related approach is proposed in~\cite{wang2026momentumaccelerated}, where a multi-step Richardson iteration with varying relaxation parameters is introduced and referred to as Richardson($m$), \ie,
\[
    \mathbf{x}_k=T_m(A,\boldsymbol{\alpha})\mathbf{x}_{k-1}+\mathbf{g}_m,
\]
where $\boldsymbol{\alpha}=(\alpha_1,\ldots,\alpha_m)$,
\[
    T_m(A,\boldsymbol{\alpha}) = (I-\alpha_m A)(I-\alpha_{m-1}A)\cdots(I-\alpha_1A),
\]
and
\[
    \mathbf{g}_m =\left( \alpha_m I + \sum_{j=1}^{m-1} \left[ \prod_{i=j+1}^{m}(I-\alpha_iA) \right]\alpha_j \right)\mathbf{b}.
\]
For SPD systems, the relaxation parameters can be chosen using Chebyshev optimal approximation based on the smallest and largest eigenvalues of $A$~\cite{mason2002chebyshev}.
To avoid explicitly computing or estimating such spectral information for each problem instance, the authors consider parameterized PDEs and employ a neural network to learn the mapping from the PDE parameters to $\boldsymbol{\alpha}$.
The resulting learned Richardson($m$) iteration can further be used as a smoother within the Fourier neural solver (FNS) framework~\cite{cui2022fourier,cui2025hybrid}; see~\Cref{sec:linear-precond-alternating}.

\citet{arisaka2023principled} propose a meta-learning-based framework for learning solver parameters (or initial guesses) for iterative methods over a distribution of related problem instances.
They focus on whether learned solver parameters can effectively accelerate convergence.
A common approach is to train parameters by minimizing the solution error after a fixed number of iterations $m$, denoted by $\mathcal{L}_m$, for example, $\mathcal{L}_m = \|\mathbf{x}_m-A^{-1}\mathbf{b}\|_2^2$.
However, the actual goal of acceleration is to minimize the number of iterations required to reach a prescribed tolerance, denoted by $\mathcal{L}_\epsilon$.
Importantly, minimizing $\mathcal{L}_m$ does not necessarily minimize $\mathcal{L}_\epsilon$.
To align the training objective with the goal of accelerating convergence, they propose gradient-based meta-solving (GBMS), which learns solver parameters by directly targeting $\mathcal{L}_\epsilon$.
Since $\mathcal{L}_\epsilon$ is discrete and non-differentiable, they introduce the differentiable surrogate
\[
    \widetilde{\mathcal{L}}_\epsilon^{(k+1)} = \widetilde{\mathcal{L}}_\epsilon^{(k)} + \sigma_a({\mathcal{L}}_k-\epsilon),
\]
where $\widetilde{\mathcal{L}}_\epsilon^{(0)} = 0$, $\sigma_a(\cdot)$ is a sigmoid function with parameter $a$, and $\mathcal{L}_\epsilon(\tau;\omega) = \displaystyle \lim_{k\to\infty} \displaystyle \lim_{a\to\infty} \widetilde{\mathcal{L}}_\epsilon^{(k)}(\tau;\omega)$.

In~\cite{arisaka2024accelerating}, GBMS is extended to legacy numerical solvers that do not support automatic differentiation, such as those in the Portable, Extensible Toolkit for Scientific Computation (PETSc)~\cite{osti_2998643}.
The resulting non-intrusive GBMS (NI-GBMS) estimates unavailable solver gradients by finite differences and reduces their variance using a control variate constructed from an adaptive surrogate model.

Parameters in classical preconditioners and subspace correction methods can also be learned by ML.
For organizational consistency, they are discussed separately:
\begin{itemize}
    \item \textbf{Parameters in classical preconditioners.}
          These are often learned within the framework of \emph{automated solver selection}~\Cref{sec:linear-solver_selection}.
    \item \textbf{Parameters in subspace correction methods.}
          Learning the strong threshold in algebraic multigrid is discussed in~\Cref{sec:linear-mgdd-coarsening}, while other parameters in multigrid and domain decomposition methods are discussed in~\Cref{sec:linear-mgdd-parameters}.
\end{itemize}

\paragraph{\textbf{Discussion}}

Since the performance of parameter learning methods can be problem- and setup-dependent, we do not report a representative numerical comparison here and refer readers to the original studies for detailed evaluations.
We make the following remarks on the future development of \emph{parameter learning}:
\begin{enumerate}
    \item \textbf{Non-intrusive hybridization:}
          Similar to \emph{initial guess learning}~(\Cref{sec:linear-core_components-initial}), many \emph{parameter learning} approaches can be integrated into existing solvers without modifying their internal implementations or accessing low-level solver information.
          Such non-intrusive integration is particularly attractive for practical use, especially with legacy solver codes.
    \item \textbf{Intrusive hybridization via adaptive parameter learning}: 
          In contrast, learning adaptive rather than globally fixed parameters may offer greater flexibility and potentially improve the convergence behavior.
          Examples include learning spatially varying weights in weighted least-squares problems~\cite{wang2023learningbased} (see~\Cref{sec:linear-mgdd-transfer_operator}) and local pseudo-time steps in~\cite{zandbergen2025improving} (see~\Cref{sec:nonlinear}).
          However, such methods typically require more intrusive modifications to the underlying solver.
          Other adaptive parameters, such as local strong-threshold parameters in coarsening, could also be explored.
\end{enumerate}

\subsection{Preconditioner learning}\label{sec:linear-precond}

Preconditioning is one of the most effective and extensively studied techniques for accelerating the convergence of iterative methods.
However, constructing a suitable preconditioner is generally problem-dependent and may be computationally expensive, particularly when solving a sequence of related linear systems.
Recently, ML has emerged as a promising approach for constructing or enhancing preconditioners.

In~\Cref{sec:linear-precond-prelims}, we recap the fundamentals of preconditioning and discuss the motivation for \emph{preconditioner learning}.
We then organize existing \emph{preconditioner learning} methods into three categories.
First, in~\Cref{sec:linear-precond-direct}, we present methods that employ ML models for \emph{learning the action of the inverse}.
Second, in~\Cref{sec:linear-precond-hybrid}, we discuss \emph{enhancing classical preconditioners}, which combine classical preconditioners with learned components, such as incomplete Cholesky or LU factorizations with learned factors.
Third, in~\Cref{sec:linear-precond-alternating}, we review \emph{alternating hybrid iterative methods}, which alternate between a classical iterative solver and an ML-preconditioned iterative solver.
This categorization is not unique and is intended only to provide a convenient organization of the existing literature.

\subsubsection{Preliminaries and motivation}\label{sec:linear-precond-prelims}
As reviewed in~\Cref{sec:linear-core_components-prelims}, commonly used iterative methods for linear systems can be broadly classified into stationary methods and Krylov subspace methods.
For Krylov subspace methods, convergence can depend strongly on the properties of the coefficient matrix.
Preconditioning aims to transform the original system into an equivalent one for which the iterative method converges more rapidly.

Let $M\in\mathbb{R}^{N\times N}$ be a nonsingular matrix, referred to as the \emph{preconditioner}.
The original linear system $A\mathbf{x}=\mathbf{b}$ may be transformed using left preconditioning,
\[
    M^{-1}A\mathbf{x} = M^{-1}\mathbf{b},
\]
or using right preconditioning,
\[
    AM^{-1}\mathbf{y} = \mathbf{b}, \qquad \mathbf{x} = M^{-1}\mathbf{y}.
\]
An effective preconditioner should improve the convergence of the iterative method while remaining inexpensive to apply.
Making the preconditioned operator close to the identity is one common strategy, but it is not necessary in general.
For SPD systems solved by CG, the condition number provides a convergence bound, while spectral clustering also influences the convergence speed.
For nonsymmetric or nonnormal systems, convergence analysis is more involved: GMRES convergence cannot generally be inferred from eigenvalues or the condition number alone.
Indeed, any prescribed positive nonincreasing residual-norm history through iteration $N-1$ can occur for an $N\times N$ matrix with prescribed nonzero eigenvalues~\cite{greenbaum1996any}.

\paragraph{\textbf{Preconditioned Krylov methods}}

Applying a Krylov subspace method to a preconditioned system yields a preconditioned Krylov method, such as preconditioned CG (PCG) or preconditioned GMRES.
Standard preconditioned Krylov methods typically employ a fixed linear preconditioner throughout the iterative process.
For PCG, the fixed linear preconditioner must additionally be SPD.
If the preconditioner is a nonlinear operator or varies at every step, a flexible Krylov method is generally required, such as flexible CG (FCG)~\cite{notay2000flexible} and flexible GMRES (FGMRES)~\cite{saad_flexible_1993,van_der_vorst_gmresr_1994}.
These methods serve as representative examples of standard and flexible preconditioned Krylov methods, respectively.
For comprehensive treatments of classical iterative methods and preconditioning for linear systems, we refer readers to~\cite{saad_iterative_2003,ciaramella_iterative_2022,greenbaum_iterative_1997,axelsson_iterative_1994}.

\paragraph{\textbf{Classical choice of preconditioners}}

Widely used classes of preconditioners~\cite{quarteroni2007numerical, benzi2002preconditioning, golub2013matrix, wathen2015preconditioning, pearson2020preconditioners} include the following:
\begin{enumerate}
    \item \textbf{Matrix splitting preconditioners} are constructed from matrix splittings and include Jacobi, Gauss-Seidel, SOR, and their block variants.
          For example, choosing $M=\operatorname{diag}(A)$ yields the Jacobi preconditioner, whereas choosing $M$ as the block-diagonal part of $A$ yields the block Jacobi preconditioner; both are well suited for parallel computing.
    \item \textbf{Incomplete factorization preconditioners}, such as incomplete Cholesky (IC)~\cite{manteuffel1980incomplete} and incomplete LU (ILU)~\cite{meijerink1977iterative}, provide sparse approximations to exact matrix factorizations.
          IC is typically used for SPD systems in combination with PCG, whereas ILU is commonly employed for non-SPD systems, often together with GMRES.
    \item \textbf{Sparse approximate inverse (SPAI) preconditioners}~\cite{benzi1999comparative} construct a sparse matrix $R\approx A^{-1}$, typically by minimizing $\|AR-I\|_F$ over a prescribed sparsity pattern~\cite{wathen2015preconditioning}.
          Since their application involves sparse matrix-vector products rather than triangular solves, they are naturally amenable to parallel computation~\cite{benzi2002preconditioning}.
    \item \textbf{Subspace correction preconditioners}~\cite{xu1992iterative, tang2009comparison} include multigrid, domain decomposition, deflation, and augmentation methods.
          These approaches reduce error components through corrections on suitably chosen subspaces and can also be used as standalone iterative solvers.
          Further details are provided in~\Cref{sec:linear-mgdd-prelims}.
    \item \textbf{Other preconditioners} include polynomial preconditioners~\cite{johnson1983polynomial}, Green's function-based preconditioners~\cite{Fraenkel_1972}, hierarchical basis preconditioners~\cite{yserentant1986multi}, hierarchical matrix-based preconditioners~\cite{bebendorf2003existence}, and others.
          In this review, we focus primarily on those that have been combined with ML techniques.
\end{enumerate}

\paragraph{\textbf{Motivation for preconditioner learning}}
A major challenge in preconditioning is that effective preconditioners are inherently problem-dependent: there is no universally optimal preconditioner across different iterative methods and problem classes.
Their construction may therefore need to be repeated when the coefficient matrix varies across a sequence of related linear systems.
Moreover, their effectiveness often comes at the expense of computational cost.
More sophisticated preconditioners can substantially reduce the number of iterations, but may also require considerable computational effort in construction or application.
In this context, ML has emerged as a promising tool for learning problem-adapted preconditioners or improving the construction of classical ones.

\begin{figure}[!htbp] 
    \centering
    \includegraphics[width=1.0\textwidth]{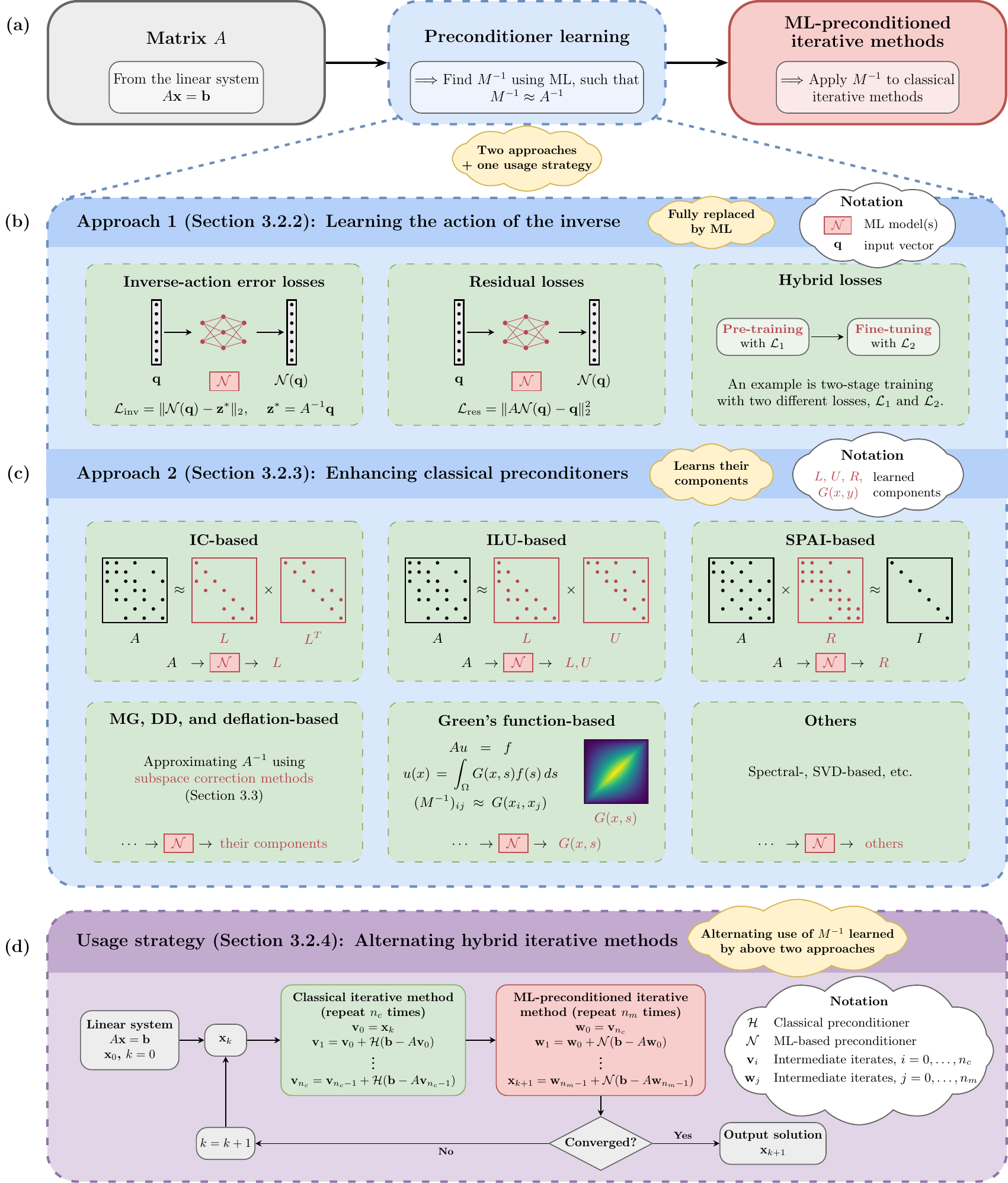}
    \caption{
        Overview of \emph{preconditioner learning}~(\Cref{sec:linear-precond}).
        \textbf{(a)} Given a coefficient matrix $A$, preconditioner learning seeks to construct or enhance a preconditioner $M^{-1}$ that approximates $A^{-1}$ or its action.
        The arrows indicate the workflow of constructing or enhancing $M^{-1}$ from $A$ and subsequently incorporating it into a classical iterative method.
        We summarize two main approaches and one usage strategy.
        \textbf{(b)} \emph{Learning the action of the inverse}~(\Cref{sec:linear-precond-direct}), in which ML models directly approximate the mapping induced by $A^{-1}$.
        \textbf{(c)} \emph{Enhancing classical preconditioners}~(\Cref{sec:linear-precond-hybrid}), in which ML models learn components of classical preconditioners.
        \textbf{(d)} \emph{Alternating hybrid iterative methods}~(\Cref{sec:linear-precond-alternating}), which interleave classical iterative methods with ML-preconditioned steps.
        In \textbf{(d)}, a representative instance is shown in which a classical stationary iteration is alternated with an ML-preconditioned Richardson iteration.
        More generally, both schemes in \textbf{(d)} can be replaced by preconditioned Krylov methods.
    }
    \label{fig:precond}
\end{figure}

\Cref{fig:precond}(a) illustrates the general workflow of \emph{preconditioner learning}, where ML is employed to approximate or enhance the inverse $A^{-1}$ or its action.
We broadly categorize existing approaches into three groups, comprising two approaches and one usage strategy:
\begin{enumerate}
    \item \textbf{Learning the action of the inverse:}
          A straightforward strategy is to use an ML model to approximate the action of $A^{-1}$, typically over a family of related linear systems, such as parameterized systems.
          The resulting learned preconditioner can be incorporated into Richardson iterations or Krylov subspace methods.
          This approach is particularly attractive when one seeks a reusable, problem-dependent preconditioner for a class of related systems.
    \item \textbf{Enhancing classical preconditioners:}
          ML can also be used to enhance classical preconditioners by learning selected components within established frameworks.
          Rather than replacing the classical construction entirely, these approaches aim to reduce setup costs or improve adaptivity while preserving the structure and favorable properties of classical preconditioning techniques.
          Typical examples include learned components in incomplete factorizations, sparse approximate inverses, multigrid, and domain decomposition methods.
          In most cases, the resulting preconditioner remains linear.
    \item \textbf{Alternating hybrid iterative methods:}
          Preconditioners learned by either of the above approaches can be applied at every iteration, yielding ``ML-preconditioned iterative methods''.
          Alternatively, ML-preconditioned and classical iterations can be applied alternately, forming an \emph{alternating hybrid iterative method}.
          This usage strategy can exploit the complementary strengths of the learned and classical components, for example, by reducing low-frequency errors with ML and refining the remaining errors with classical iterations.
\end{enumerate}

Additionally, \emph{update function learning} in~\Cref{sec:linear-core_components-update} and \emph{subspace correction methods} in~\Cref{sec:linear-mgdd} are closely related to \emph{preconditioner learning}.
First, as discussed in~\Cref{sec:linear-core_components-update}, the reviewed approaches to \emph{update function learning} learn a preconditioner, although update functions can also be learned in other ways.
Second, as we mentioned above, \emph{subspace correction methods} can also be used as preconditioners, and their enhancements are more naturally associated with~\Cref{sec:linear-precond-hybrid}.

Neural networks are by far the most widely adopted ML technique, and the resulting preconditioners are therefore often referred to as ``neural preconditioners''~\cite{chen2025graph, trifonov2025learning, li2025neural}.

\subsubsection{Learning the action of the inverse}\label{sec:linear-precond-direct}

Let $\mathcal{N}: \mathbb{R}^n \to \mathbb{R}^n$ denote an ML model, most commonly a neural network.
Rather than constructing an explicit approximation of the matrix $A^{-1}$, the model is trained to approximate its action, \ie,
\[
    \mathcal{N}(\mathbf{q}) \approx A^{-1}\mathbf{q},
    \qquad \mathbf{q}\in\mathbb{R}^n.
\]
On the other hand, learning an explicit matrix representation of $A^{-1}$ instead leads to sparse approximate inverse (SPAI)-based preconditioners, which are discussed in~\Cref{sec:linear-precond-hybrid}.
Although the target operator $A^{-1}$ is linear, the learned mapping $\mathcal{N}$ is generally nonlinear due to the nonlinear architecture used to represent it.
Consequently, its use within iterative solvers may require algorithms that can accommodate nonlinear preconditioning.

Assume that the exact solution is given by $\mathbf{x}^* = A^{-1}\mathbf{b}$.
At the $k$-th iteration, the current approximation is updated as $\mathbf{x}_{k+1} = \mathbf{x}_k + \mathbf{d}_k$.
Ideally, the correction $\mathbf{d}_k$ should be such that the updated iterate satisfies $A\mathbf{x}_{k+1} = \mathbf{b}$.
Substituting the update formula yields $A\mathbf{d}_k = \mathbf{b} - A\mathbf{x}_k = \mathbf{r}_k$, where $\mathbf{r}_k$ denotes the residual.
Hence, the ideal correction is $A^{-1}\mathbf{r}_k$.

Loss functions used for \emph{learning the action of the inverse} can be broadly classified into three categories:
\begin{enumerate}
    \item \textbf{Inverse-action error losses}: These losses directly compare the network output with a reference inverse action.
          For example,
          \[
              \mathcal{L}_{\mathrm{inv}} =
              \left\lVert
              \mathcal{N}(\mathbf{q})-\mathbf{z}^*
              \right\rVert_2^2,
          \]
          where $\mathbf{q}\in\mathbb{R}^n$ denotes the input vector, and $\mathbf{z}^*\in\mathbb{R}^n$ is the corresponding target vector.
          In particular, when $\mathbf{q}=\mathbf{b}$, the target is the exact solution $\mathbf{z}^*=\mathbf{x}^*=A^{-1}\mathbf{b}$; when $\mathbf{q}=\mathbf{r}$, the target is the ideal correction $\mathbf{z}^*=\mathbf{d}^*=A^{-1}\mathbf{r}$.
    \item \textbf{Residual losses}: These losses evaluate the predicted inverse action through the residual of the associated linear system.
          For example,
          \[
              \mathcal{L}_{\mathrm{res}} =
              \left\lVert
              A\mathcal{N}(\mathbf{q})-\mathbf{q}
              \right\rVert_2^2.
          \]
          The term ``residual loss'' refers to the fact that the network prediction is evaluated through the residual of the equation $A\mathbf{z}^*=\mathbf{q}$.
    \item \textbf{Hybrid losses}: Hybrid approaches combine multiple loss functions, either simultaneously within a single training objective or sequentially across different training stages.
          For example, in~\cite{dimola2025neural}, an ML model is first pretrained using the ``residual loss'' above and then fine-tuned using a loss defined in terms of the principal angles between Krylov subspaces.
\end{enumerate}

For simplicity, the formulations above are expressed in terms of the absolute squared $L^2$ norm, although relative errors are also commonly used in practice.
\Cref{tab:linear-precond-direct-comparison} summarizes several representative works on \emph{learning the action of the inverse}.
In the following, we review these methods according to the three categories of loss functions introduced above.

\begin{table}[!htbp]
    \renewcommand{\arraystretch}{1.3}
    \centering
    \caption{
        Comparison of several existing approaches for \emph{learning the action of the inverse}~(\Cref{sec:linear-precond-direct}).
        Here, $\mu$ denotes the problem parameters, and $\mathcal{I}$ denotes the geometric information.
        For~\cite{kaneda2023deep, bai2026hybrid}, $\alpha$ denotes the iteration step size.
        For~\cite{kaneda2023deep, lan2024neuralpreconditioned}, PSDO stands for ``preconditioned steepest descent with orthogonalization.''
        For~\cite{dimola2025neural}, $\sin\left(\angle(\mathbf r,\mathcal W)\right)$ denotes the sine of the principal angle between the residual $\mathbf{r}$ and $\mathcal{W}$, the image of the Krylov subspace.
    }
    \label{tab:linear-precond-direct-comparison}
    \resizebox{\linewidth}{!}{%
        \begin{minipage}{\linewidth}
            \centering
            \begin{tabular}{
                    p{1.0cm}
                    >{\raggedright\arraybackslash}p{6.6cm}
                    >{\raggedright\arraybackslash}p{1.7cm}
                    >{\raggedright\arraybackslash}p{1.3cm}
                    >{\raggedright\arraybackslash}p{3.3cm}
                    >{\raggedright\arraybackslash}p{2.0cm}
                }
                \toprule
                \textbf{ref.}
                 & \textbf{tested problem}
                 & $\mathcal{N}$
                 & \textbf{ML}
                 & \textbf{loss function}
                 & \textbf{iter. method}
                \\
                \midrule
                \cite{ackmann2020machinelearned}
                 & shallow water equations (time-dependent)
                 & $\mathbf{d} = \mathcal{N} (\mathbf{r})$
                 & FNN
                 & $\left\lVert \mathcal{N}(\mathbf{r}) - \mathbf{d} \right\rVert_2$
                 & GCR
                \\
                \cite{rudikov2024neural}
                 & elliptic PDEs (different resolutions)
                 & $\mathbf{d} = \mathcal{N} (\mathbf{r})$
                 & SNO
                 & $\left\lVert \mathcal{N}(\mathbf{r}) - \mathbf{d} \right\rVert_A / \left\lVert \mathbf{d} \right\rVert_A$
                 & FCG
                \\
                \multirow{2}{*}{\cite{giraud2025neural, xiang2025unsupervised}}
                 & parametric Helmholtz equations~\cite{giraud2025neural}, heterogeneous fluid PDEs~\cite{xiang2025unsupervised}
                 & \multirow{2}{*}{$\mathbf{x} = \mathcal{N} (\mathbf{b}, \mu)$}
                 & \multirow{2}{*}{U-Net}
                 & \multirow{2}{*}{$\left\lVert A\mathcal{N}(\mathbf{b}, \mu) - \mathbf{b}\right\rVert_2^2 / \left\lVert\mathbf{b}\right\rVert_2^2$}
                 & \multirow{2}{*}{FGMRES}
                \\
                \cite{kaneda2023deep}
                 & Poisson equations from incompressible flow
                 & $\mathbf{d} = \mathcal{N} (\mathbf{r})$
                 & CNN
                 & $\left\lVert \alpha A\mathcal{N}(\mathbf{r}) - \mathbf{r}\right\rVert_2$
                 & PSDO
                \\
                \multirow{1}{*}{\cite{lan2024neuralpreconditioned}}
                 & Poisson equations with mixed boundary conditions
                 & \multirow{1}{*}{$\mathbf{d} = \mathcal{N}(\mathbf{r}, \mathcal{I})$}
                 & CNN
                 & \multirow{1}{*}{$\left\lVert A\mathcal{N}(\mathbf{r}, \mathcal{I}) - \mathbf{r}\right\rVert_2$}
                 & \multirow{1}{*}{PSDO}
                \\
                \multirow{1}{*}{\cite{bai2026hybrid}}
                 & Poisson equations from incompressible flow
                 & \multirow{1}{*}{$\mathbf{d} = \mathcal{N} (\mathbf{r})$}
                 & DeepONet
                 & \multirow{1}{*}{$\left\lVert \alpha A\mathcal{N}(\mathbf{r}) - \mathbf{r}\right\rVert_2$}
                 & \multirow{1}{*}{Richardson}
                \\
                \cite{chen2025graph}
                 & general sparse matrix (not just PDEs)
                 & $\mathbf{x} = \mathcal{N} (\mathbf{b})$
                 & GCN
                 & $\left\lVert A\mathcal{N}(\mathbf{b}) - A \mathbf{x}\right\rVert_1$
                 & FGMRES
                \\
                \multirow{2}{*}{\cite{dimola2026numerical,dimola2025neural}}
                 & \multirow{2}{*}{parametric 3D-1D mixed-dimensional PDEs}
                 & \multirow{2}{*}{$\mathbf{d} = \mathcal{N} (\mathbf{r}, \mu)$}
                 & \multirow{2}{*}{U-Net}
                 & $\left\lVert A\mathcal{N}(\mathbf{r}) - \mathbf{r}\right\rVert_2$ \cite{dimola2025neural, dimola2026numerical} \& $\left|\sin\!\left(\angle(\mathbf r,\mathcal W)\right)\right|$ \cite{dimola2025neural}
                 & \multirow{2}{*}{FGMRES}
                \\
                \bottomrule
            \end{tabular}
        \end{minipage}
    }
\end{table}

\paragraph{\textbf{Inverse-action error losses}}

An early work by~\citet{ackmann2020machinelearned} considers linear systems arising from the numerical solution of the shallow-water equations with semi-implicit timestepping, where a linear system must be solved at each time step.
The authors employ an FNN to construct a preconditioner for the generalized conjugate residual (GCR) method~\cite{eisenstat1983variational}.
The resulting preconditioner is local, as each grid-point correction is predicted solely from information in a small neighboring stencil.
For each grid point, the FNN takes local $3\times 3$ or $5\times 5$ grid stencils of the residual and the six coefficient fields defining the discretized linear operator as input, and predicts the increment to the fluid thickness at the central grid point.

\citet{rudikov2024neural} propose FCG-NO for linear systems arising from elliptic PDEs, where discretization-invariant neural operators, such as spectral neural operators (SNOs)~\cite{fanaskov2023spectral} and Fourier neural operators (FNOs)~\cite{li2021fourier}, are trained as nonlinear preconditioners for FCG.
The discretization invariance enables preconditioners trained on coarse grids to be transferred to finer grids.
For training, the authors sample $A$ and $\mathbf{x}^*$, set $\mathbf{b}=A\mathbf{x}^*$, and run unpreconditioned CG to collect residual-error pairs $\{\mathbf{r}_k, \mathbf{d}_k\}$, where $\mathbf{r}_k$ is in the Krylov subspace, and $\mathbf{d}_k=\mathbf{x}^*-\mathbf{x}_k$.

\paragraph{\textbf{Residual losses}}

As discussed in~\Cref{sec:linear-core_components-update}, \citet[Chapter~5.3.2]{xiang2022solution} interprets the trained DNN, which maps the current solution $\mathbf{x}_k$ and residual $\mathbf{r}_k$ to an update $\mathbf{d}_k$, as approximating the action of $A^{-1}$.
The network is therefore used as a nonlinear preconditioner for FGMRES and flexible full orthogonalization method (FFOM) in solving Helmholtz equations.
Building on this idea, \citet{giraud2025neural} employ a U-Net to approximate the action of the inverse of the discretized parametric Helmholtz operator and use it as a nonlinear preconditioner for FGMRES.
A similar U-Net-based nonlinear preconditioner for FGMRES is proposed by~\citet{xiang2025unsupervised} for heterogeneous fluid PDEs.

\citet{kaneda2023deep} propose the deep conjugate direction method (DCDM) for solving SPD linear systems arising from discretized Poisson equations.
A CNN acts as a nonlinear preconditioner that maps the residual $\mathbf{r}_k$ to an approximation of $A^{-1}\mathbf{r}_k$.
Its output is then orthogonalized against previous search directions using Gram-Schmidt to obtain the $A$-conjugate direction $\mathbf{d}_k$.
The iterate is updated by
\[
    \mathbf{x}_{k+1} = \mathbf{x}_k + \alpha_k \mathbf{d}_k,
    \qquad
    \alpha_k=\frac{\mathbf{r}_k^T\mathbf{d}_k}{\mathbf{d}_k^TA\mathbf{d}_k},
\]
where the step size $\alpha_k$ is determined by an exact line search, analogously to the CG method.
This approach is referred to as preconditioned steepest descent with orthogonalization (PSDO) in~\cite{lan2024neuralpreconditioned}.

\citet{lan2024neuralpreconditioned} extend DCDM~\cite{kaneda2023deep} to Poisson equations with mixed Dirichlet and Neumann boundary conditions arising in incompressible fluid simulations, where the domain shape and boundary conditions evolve over time.
The resulting geometry is encoded by a three-channel tensor $\mathcal{I}$ representing fluid, solid, and air regions, and the inverse action is approximated by $A^{-1}\mathbf{r}\approx \mathcal{N}(\mathcal{I},\mathbf{r})$.
Specifically, local patches of $\mathcal{I}$ are transformed through trainable affine maps into spatially varying convolution kernels, which are then applied to the residual $\mathbf{r}$.
Unlike the nonlinear preconditioner in~\cite{kaneda2023deep}, $\mathcal{N}(\mathcal{I},\mathbf{r})$ is linear in $\mathbf{r}$ for fixed $\mathcal{I}$.

Building on DCDM~\cite{kaneda2023deep}, \citet{bai2026hybrid} propose the deep learning line search method (DLSM), in which a DeepONet predicts the search direction from the residual,
\[
    \mathbf{x}_{k+1}=\mathbf{x}_k+\alpha_k\mathbf{d}_k,
    \qquad
    \mathbf{d}_k = \mathcal{N}(\mathbf{r}_k / \|\mathbf{r}_k\|_2),
\]
with the step size $\alpha_k$ determined by the same line-search procedure as in DCDM.
Unlike DCDM, DLSM omits Gram-Schmidt orthogonalization, which is shown numerically to provide negligible benefit.
The learned iteration is further combined with a classical CG-type method in an alternating HIM; see~\Cref{sec:linear-precond-alternating}.

\citet{chen2025graph} proposes graph neural preconditioners (GNPs), which leverage GCNs~\cite{kipf2017semisupervised} to approximate the action of $A^{-1}$ within FGMRES.
For each matrix $A$, training pairs $(\mathbf{x},\mathbf{b})$, with $\mathbf{b}=A\mathbf{x}$, are generated from a mixture of Gaussian samples and Arnoldi-based samples emphasizing the bottom spectral subspace.
To reflect the scale equivariance of $A^{-1}$, \ie, $A^{-1}(\alpha\mathbf{b})=\alpha A^{-1}\mathbf{b}$ for any scalar $\alpha \neq 0$, a scaling operator and its inverse are incorporated into the network architecture.
The method is evaluated on $867$ general sparse non-SPD matrices from the SuiteSparse Matrix Collection~\cite{davis2011university}, covering 50 application areas.

\citet{dimola2026numerical} consider parametrized linear systems arising from 3D-1D mixed-dimensional elliptic PDEs, where variations in the embedded 1D geometry lead to a family of ill-conditioned systems.
They employ a U-Net to learn a nonlinear preconditioner that takes both the input vector and a geometric descriptor of the 1D manifold as inputs, enabling generalization across different geometries without retraining.
The network is trained with a residual loss, $\left\lVert A\mathcal{N}(\mathbf{r}) - \mathbf{r}\right\rVert_2$, and used within FGMRES.

\paragraph{\textbf{Hybrid losses}}

In~\cite{dimola2025neural}, the authors extend the approach in~\cite{dimola2026numerical} into a two-stage training strategy.
The first stage pretrains the neural preconditioner by minimizing a residual loss.
Let $\mathcal{K}_k(A,\mathbf r_0)$ denote the Krylov subspace and $\mathcal{W}_k := \{A\mathbf p \mid \mathbf p\in\mathcal{K}_k(A,\mathbf r_0)\}$ its image subspace.
For GMRES-type methods, the residual satisfies $\mathbf r_k \perp \mathcal{W}_k$, and its reduction is governed by the principal angle between $\mathbf{r}_{k-1}$ and $\mathcal{W}_k$~\cite{eiermann2001geometric}.
In FGMRES, letting $\mathbf v_j$ denote the $j$-th Arnoldi basis vector, the corresponding image subspace is generated by the preconditioned vectors $A\mathcal N_\theta(\mathbf v_j)$.
The pretrained model is then fine-tuned using a solver-aware loss based on this angle,
\[
    \left|\sin\left(\angle(\mathbf r_{k-1},\mathcal W_k)\right)\right|,
\]
which is obtained from the Givens rotations generated during the Arnoldi process, thereby allowing the loss to be differentiated with respect to the neural preconditioner parameters.
Minimizing this loss is designed to promote faster FGMRES convergence.

\paragraph{\textbf{Discussion}}

When neural networks are used to approximate the action of $A^{-1}$, the resulting preconditioners are generally nonlinear.
Such learned preconditioners can naturally be incorporated into iterative methods that admit nonlinear preconditioning.
For Krylov subspace methods, however, nonlinear preconditioning destroys the standard Krylov subspace structure.
In this case, flexible Krylov subspace methods, such as FCG or FGMRES, are generally required.
If one intends to employ learned preconditioners within standard Krylov methods, the learned preconditioner must still define a fixed linear operator; see, for example, the methods discussed in~\Cref{sec:linear-precond-hybrid}.

Some general remarks on \emph{learning the action of the inverse} are as follows:
\begin{enumerate}
    \item 
          The training data generation is crucial to the performance of learned preconditioners and may incur substantial computational cost.
          When paired samples are required, a common strategy to avoid explicitly solving linear systems is to first sample vectors $\mathbf{x}$ or $\mathbf{d}$ and then apply $A$ to construct the corresponding $\mathbf{b}$ or $\mathbf{r}$, respectively.
          The sampled vectors can be generated in several ways:
          \begin{enumerate}
              \item \textbf{Standard normal sampling}:
                    Training vectors are sampled from a standard normal distribution, as in~\cite{giraud2025neural, xiang2025unsupervised}.
                    This strategy is simple and inexpensive, but does not explicitly target directions that are difficult for the iterative solver.
                    Consequently, a large number of samples may be required for challenging problems.
              \item \textbf{Krylov subspace}:
                    Training vectors, such as $\mathbf{r}$ and $\mathbf{d}$, are collected during unpreconditioned Krylov iterations, as in~\cite{rudikov2024neural}.
                    Such data reflect directions naturally encountered during the iterative solution process.
              \item \textbf{Low-spectrum subspace sampling}:
                    Training vectors are sampled from an approximate eigenspace of $A$, with emphasis on eigenvectors associated with the lower end of the spectrum.
                    Compared with standard normal sampling, this strategy requires additional cost for approximating the eigenspace.
                    It has been adopted in~\cite{kaneda2023deep, lan2024neuralpreconditioned, bai2026hybrid}.
              \item \textbf{Hybrid sampling}:
                    Different sampling strategies can also be combined.
                    For example, in~\cite{chen2025graph}, half of the training vectors are sampled from a standard normal distribution and the other half from an approximate low-spectrum subspace constructed by the Arnoldi process.
          \end{enumerate}
    \item 
          For small systems, learned preconditioners may not be competitive with classical approaches in terms of either iteration counts or overall runtime, mainly because of the additional inference cost.
          For larger systems, however, they may become more advantageous, although such gains come at the expense of training-data generation and model training.
    \item 
          Since learned preconditioners are generally nonlinear, they are typically used with flexible Krylov methods rather than standard preconditioned Krylov methods, which may introduce additional computational and memory overhead.
          For SPD systems, for example, FCG is generally more expensive than PCG in both computation and memory due to the additional orthogonalization required for flexibility.
          For non-SPD systems, FGMRES is used instead of preconditioned GMRES, typically with higher memory requirements, while the additional computational overhead is often more moderate.
\end{enumerate}

Given the above discussion of existing studies, several open directions remain:
\begin{enumerate}
    \item \textbf{Training data generation:}
          Existing strategies, such as sampling from approximate eigenspaces associated with the smallest eigenvalues, can effectively emphasize low-frequency components.
          However, approximating the corresponding spectral subspaces can be computationally expensive.
          Developing effective yet low-cost strategies for training data generation therefore remains an important direction.
    \item \textbf{Multistage training with hybrid losses:}
          Multistage training strategies that combine different loss functions may improve the effectiveness of learned preconditioners.
          For example, a pretrained model can be further fine-tuned during the iterative solution process using an additional loss, as considered in~\cite{dimola2025neural, wu2026are}.
          Other combinations of training objectives remain to be explored.
    \item \textbf{Preserving preconditioner properties:}
          When nonlinear ML models are used to approximate linear preconditioning operators, important algebraic and structural properties may not be preserved automatically.
          Examples include scale equivariance, as considered in~\cite{chen2025graph}, as well as locality, symmetry, and sparsity.
          Incorporating other general or problem-specific structural constraints into ML architectures may lead to more robust and effective learned preconditioners.
    \item \textbf{Extensions to block preconditioners:}
          \emph{Preconditioner learning} may also be extended to block preconditioners, in which ML models are used to approximate or enhance preconditioners for individual blocks.
          This direction may be particularly relevant for systems with natural block structures, such as coupled multiphysics problems.
\end{enumerate}

\subsubsection{Enhancing classical preconditioners}\label{sec:linear-precond-hybrid}
Classical preconditioners, including IC, ILU, SPAI, or those based on multigrid~(MG), domain decomposition~(DD), and deflation methods, exploit properties of the coefficient matrix and the underlying problem structure to accelerate iterative solvers.
Their effectiveness, however, often relies on constructing problem-dependent components, which can incur substantial setup costs and require careful parameter selection.
Examples include sparsity patterns and matrix entries in IC, ILU, and SPAI preconditioners, transfer operators in multigrid methods, and coarse spaces in domain decomposition methods.

Unlike the approaches in~\Cref{sec:linear-precond-direct}, which directly learn the action of $A^{-1}$, the methods considered in this category incorporate ML into established preconditioning frameworks to construct or improve selected components.
The objective is to reduce setup cost, improve adaptivity or preconditioning quality, or both, while retaining the structural advantages of the underlying classical method.

For IC, ILU, and SPAI preconditioners, the matrices to be constructed are typically sparse.
Accordingly, ML methods are usually employed to predict only their nonzero entries.
In this context, GNNs are particularly suitable for learning such matrices due to the natural duality between graphs and matrices, with the sparsity pattern defining the graph connectivity; see~\Cref{fig:prelims-struct_unstruct_data}(b).
This correspondence enables GNNs to exploit the local structural information of the coefficient matrix while preserving sparsity in the learned preconditioner.
The training objectives are typically designed to enforce the algebraic properties required by the corresponding classical preconditioners, as discussed below.

Several existing studies on~\emph{enhancing classical preconditioners} are summarized in~\Cref{tab:linear-precond-hybrid-comparison}.
In the following, we organize these methods according to the underlying classical preconditioning framework.

\begin{table}[!tbp]
    \renewcommand{\arraystretch}{1.3}
    \centering
    \caption{
        Comparison of existing approaches for \emph{enhancing classical preconditioners}~(\Cref{sec:linear-precond-hybrid}).
        For IC-based preconditioners, $L_{\IC}$ denotes the factor obtained from classical IC, and $L_{\text{corr}}$ denotes a learned correction term.
        For ILU-based preconditioners, $\operatorname{S}_{L,U}$ denotes the sparsity patterns of $L$ and $U$.
        For SPAI-based preconditioners, $T$ and $D$ denote a strictly lower triangular matrix and a diagonal matrix, respectively, and $\kappa(\cdot)$ denotes the condition number.
        The vector $\mathbf{w}$ denotes a random vector used in loss functions.
    }
    \label{tab:linear-precond-hybrid-comparison}
    \resizebox{\linewidth}{!}{%
        \begin{minipage}{\linewidth}
            \centering
            \begin{tabular}{
                    llllll
                }
                \toprule
                 & \textbf{ref.}
                 & \textbf{ML model}
                 & \textbf{final prec.}
                 & \textbf{loss function}
                 & \textbf{iter. method}
                \\
                \midrule
                \multirow{1}{*}{\centering \textbf{IC}}
                 & \cite{li2023learning}
                 & $L = \GNN(A, b)$
                 & $M=L L^T$
                 & $\|LL^T\mathbf{x} - \mathbf{b}\|_2^2$
                 & PCG
                \\
                 & \cite{hausner2024neural}
                 & $L = \GNN(A)$
                 & $M=L L^T$
                 & $\|LL^T\mathbf{w} - A\mathbf{w}\|_2^2$
                 & PCG
                \\
                 & \cite{li2024deep}
                 & $L_{\IC} = \IC(A)$, $L_{\text{corr}} = \GNN(A)$
                 & $M=LL^T$, $L=L_{\IC} + L_{\text{corr}}$
                 & $\|LL^T\mathbf{w} - A\mathbf{w}\|_2^2$
                 & PCG
                \\
                 & \cite{trifonov2025learning}
                 & $L_{\IC} = \IC(A)$, $L_{\text{corr}} = \GNN(L_{\IC})$
                 & $M=LL^T$, $L=L_{\IC} + \alpha L_{\text{corr}}$
                 & $\|LL^T\mathbf{x} - \mathbf{b}\|_2^2$
                 & PCG
                \\
                \midrule
                \multirow{1}{*}{\centering \textbf{ILU}}
                 & \cite{stanaityte2020ilu}
                 & $\operatorname{S}_{L, U}$ = $\UNet(A)$
                 & $M=LU$, $(L,U)=\ILU(A; \operatorname{S}_{L, U})$
                 & binary cross-entropy
                 & prec. GMRES
                \\
                 & \cite{yusuf2024constructing}
                 & $(L, U) = \GNN(A)$
                 & $M = L U$
                 & $\|L U \mathbf{x} - \mathbf{b}\|_2$
                 & prec. LGMRES
                \\
                 & \cite{hausner2024learning}
                 & $(L, U) = \GNN(A)$
                 & $M = L U$
                 & $\|A \mathbf{w} - L U \mathbf{w}\|_2^2 + \alpha \|L U \mathbf{x}\|_2^2$
                 & prec. GMRES
                \\
                \midrule
                \multirow{1}{*}{\centering \textbf{SPAI}}
                 & \cite{ruelmann2018prospects}
                 & $R = \FNN(A)$
                 & $M^{-1}=R$
                 & $\|R-A^{-1}\|_2^2$
                 & Richardson
                \\
                 & \cite{sappl2019deep}
                 & $(T,\hat D)=\mathrm{CNN}(A)$
                 & \makecell[l]{$M^{-1}=(T+D)(T+D)^T$,\\$D=\max(\hat D,10^{-3})$}
                 & $\kappa(AM^{-1})$
                 & PCG
                \\
                 & \cite{sappl2026convolutional}
                 & $(T,\hat D)=\mathrm{CNN}(A)$
                 & \makecell[l]{$M^{-1}=(T+D)(T+D)^T$,\\$D=\log(1+\exp(\hat D))$}
                 & $\|M^{-1}A - I\|_F$
                 & PCG
                \\
                 & \cite{yang2025learning}
                 & $G=\mathrm{GNN}(A,\text{node features})$
                 & $M^{-1}=GG^\top+\varepsilon I$
                 & $\left\|\frac{AM^{-1}}{\|A\|}\textbf{w} - \textbf{w}\right\|_2^2$
                 & PCG
                \\
                \bottomrule
            \end{tabular}
        \end{minipage}
    }
\end{table}

\paragraph{\textbf{IC-based preconditioners}}

IC-based preconditioners are designed for SPD systems.
ML-based approaches learn or correct the lower-triangular factor $L$ such that $A \approx LL^T$, including its sparsity pattern and nonzero entries.
Rather than explicitly forming $M^{-1}$, the preconditioner is represented through an approximation of the factorization $M=LL^T$, so that the action of $M^{-1}$ can be applied via two triangular solves.
By preserving the IC structure, the resulting preconditioner remains SPD and can therefore be used with PCG, provided that the learned triangular factor has nonzero diagonal entries.

\citet{li2023learning} employ a GNN to learn an approximate lower-triangular factor $L$ from the SPD matrix $A$ and the right-hand side $\mathbf{b}$.
The sparsity pattern of $A$ defines the graph, with $A_{ij}$ and $\mathbf{b}_i$ used as edge and node features, respectively.
Two FNN encoders map these features to 16-dimensional latent representations, which are then updated by several message-passing layers.
The edge features are finally decoded into the entries of $L$, following the sparsity pattern of the lower-triangular part of $A$.
Ideally, one would train $L$ by minimizing the Frobenius loss $\|LL^T - A\|_F^2 = \sum_j \|LL^T\mathbf{e}_j - A\mathbf{e}_j\|_2^2$, where $\mathbf{e}_j$ denotes the $j$-th standard basis vector.
To improve training and generalization, the authors adopt the loss $\|LL^T\mathbf{x} - \mathbf{b}\|_2^2$, which incorporates the distribution of the solution vector $\mathbf{x}$.

Similarly, \citet{hausner2024neural} propose the neural incomplete factorization (NeuralIF) method, which uses a GNN to learn IC preconditioners in an unsupervised way.
The nonzero entries of $A$ are represented as edge features, and each node is associated with an eight-dimensional feature vector describing the sparsity structure and diagonal dominance of $A$.
The diagonal entries of the output $L$ are enforced to be strictly positive through an exponential transformation.
The Frobenius loss is approximated using Hutchinson's estimator~\cite{hutchinson1990stochastic}, \ie,
\[
    \|LL^T - A\|_F^2 \approx \|LL^T\mathbf{w} - A\mathbf{w}\|_2^2,
\]
where $\mathbf{w}\sim\mathcal{N}(\mathbf{0},I)$ is a standard normal random vector.

\citet{li2024deep} enhance the standard IC factor $L_{\IC}$ with a GNN-based correction $L_{\text{corr}}$,
\[
    L = L_{\IC} + L_{\text{corr}},
\]
and construct the preconditioner as $M=LL^T$.
Similar to~\cite{hausner2024neural}, Hutchinson's estimator is used to approximate the Frobenius-norm loss.
A related approach is proposed in~\cite{trifonov2025learning}, where a learnable coefficient $\alpha$ is introduced, yielding
\[
    L = L_{\IC} + \alpha L_{\text{corr}}.
\]
In contrast to~\cite{li2024deep}, the GNN takes the classical IC factor $L_{\IC}$ rather than $A$ as input, reducing the graph size by half and improving training efficiency.

\citet{millard2025pearl} formulate preconditioner learning as a contextual bandit problem and propose an approach called preconditioner enhancement through actor-critic reinforcement learning (PEARL).
The contextual bandit is a class of reinforcement learning methods~\cite{bouneffouf2020survey} in which an agent selects an action based on a given context and receives a corresponding reward, without considering long-term state transitions.
Here, the coefficient matrix $A$ serves as the context.
The actor, parameterized by an FNN, maps $A$ to a lower triangular factor $L$, yielding the preconditioner $M = L L^T + \alpha I$, where $\alpha$ stabilizes the factorization and ensures positive definiteness.
The critic, also implemented as an FNN, takes $A$ and $M$ as input and outputs a scalar reward based on solver feedback, including the iteration count, preconditioner sparsity, and final residual norm.

\paragraph{\textbf{ILU-based preconditioners}}
ILU-based preconditioners are usually used for non-SPD systems.
Analogous to IC, ML approaches can be used to learn the factors $L$ and $U$, or their sparsity patterns, such that $A \approx LU$.
The preconditioner is represented through the factorization $M=LU$ and applied via triangular solves.
The resulting preconditioner is often used with Krylov methods such as GMRES.

In standard ILU factorizations, the sparsity patterns of $L$ and $U$ are typically prescribed from the nonzero structure of $A$.
However, this can be inadequate for more complex problems, such as saddle-point systems where zero blocks in $A$ may become densely populated during the exact LU factorization.
Motivated by image segmentation, \citet{stanaityte2020ilu} employs a U-Net to predict suitable fill-in patterns of the factors $L$ and $U$ by treating $A$ as an input image and outputting a binary mask for the predicted nonzero structure.
The network is trained using binary cross-entropy, and the predicted pattern is then used to construct the ILU preconditioner.

\citet{yusuf2024constructing} propose a method called graph neural network preconditioner for advection-dominated nonsymmetric problems (GraphPAN), which uses a GNN to learn an approximate ILU factorization of sparse nonsymmetric matrices.
Given the graph representation of $A$, the GNN predicts sparse triangular factors $(L,U)$, with $\operatorname{diag}(L)=1$, yielding the preconditioner $M=L U$.
Instead of the Frobenius loss, the authors use the loss $\|LU\mathbf{x}-\mathbf{b}\|_2$, with $A\mathbf{x}=\mathbf{b}$, to penalize factorization errors on the sampled solution vectors.
The resulting preconditioner is then used with the ``loose'' GMRES (LGMRES) method~\cite{baker2005technique}.

The approach in~\cite{hausner2024neural} is extended to ILU-based preconditioners for non-SPD matrices in~\cite{hausner2024learning}, where a GNN learns the mapping $(L,U)=\GNN(A)$ and constructs the preconditioner $M=LU$.
To ensure nonsingularity, the diagonal of $U$ is fixed to one, and a specific activation function enforces nonzero diagonal entries in $L$.
Instead of directly minimizing $\kappa(AM^{-1})$, the authors design loss functions to improve the extreme singular values of the right-preconditioned matrix.
Using Hutchinson's estimator, they consider
\[
    \left\{
    \begin{aligned}
        \mathcal{L}_{\max}
         & = \|A\mathbf{w}-M\mathbf{w}\|_2^2,
        \\
        \mathcal{L}_{\min}
         & = \|MA^{-1}\mathbf{w}-\mathbf{w}\|_2^2,
    \end{aligned}
    \right.
\]
where $\mathbf w\sim\mathcal N(0,I)$.
For computational efficiency, $\mathcal L_{\min}$ is approximated using training samples satisfying $A\mathbf x=\mathbf b$ as $\widehat{\mathcal L}_{\min} = \| M \mathbf{x}-\mathbf{b}\|_2^2$.
The model is finally trained with a combined loss,
\[
    \mathcal L_{\mathrm{com}}
    = \|A\mathbf w-M\mathbf w\|_2^2
    +\alpha\|M\mathbf x\|_2^2,
\]
where $\alpha$ balances the two terms.

\paragraph{\textbf{SPAI-based preconditioners}}

SPAI-based methods explicitly construct a sparse matrix $M^{-1}$ that approximates $A^{-1}$, typically by minimizing $\|AM^{-1} - I\|_F$.
ML-based approaches can be used to predict or improve this sparse approximate inverse.
Unlike methods that learn only the action of $A^{-1}$~(\Cref{sec:linear-precond-direct}), SPAI-based approaches represent the approximate inverse explicitly in matrix form, thereby yielding a linear preconditioner.

An early attempt is presented in~\cite{ruelmann2018prospects}, where an FNN is trained in a supervised manner to map a coefficient matrix $A$ to its inverse, $A^{-1} \approx \FNN(A)$.
Both $A$ and $A^{-1}$ are vectorized as the network input and training target, respectively.
The predicted approximate inverse is then used as a preconditioner for a Richardson iteration.
Numerical results show that the learned preconditioner outperforms Jacobi and Gauss-Seidel in terms of iteration counts.
Despite the simplicity of the test problem, the results demonstrate the feasibility of learning SPAI preconditioners with neural networks.

\citet{sappl2019deep} employ a CNN to construct a preconditioner for PCG by minimizing the condition number of the preconditioned system, $\kappa(AM^{-1})$, computed via SVD.
Exploiting the symmetry of $A$, the CNN takes only its strictly lower triangular and diagonal parts as input and predicts a strictly lower triangular matrix $T$ and diagonal entries $\hat D$.
To ensure an SPD preconditioner, the diagonal is thresholded as $D=\max\{\hat D,10^{-3}\}$, and the inverse preconditioner is constructed as
\[
    M^{-1}=(T+D)(T+D)^T.
\]
A subsequent work~\cite{sappl2026convolutional} replaces the expensive condition number loss with the unsupervised Frobenius loss
\[
    \|(T+D)(T+D)^TA - I\|_F,
\]
thereby avoiding the SVD required in~\cite{sappl2019deep}.
The hard threshold on the diagonal is also replaced by a softplus function,
\[
    D=\log(1+\exp(\hat D)),
\]
which enforces positivity.
The learned preconditioner remains sparse because the sparsity of $A$ and the small convolution kernels introduce only limited local fill-in.
Unlike IC factorization, which approximates $A$, this approach directly parameterizes a factorization of $M^{-1}$ and can therefore be viewed as learning a SPAI preconditioner.

Preconditioners based on incomplete factorizations (\eg, IC and ILU) are less suitable for GPU parallelization, as their application involves triangular solves, which are inherently sequential.
To address this limitation, \citet{yang2025learning} use a GNN to construct a GPU-friendly SPAI preconditioner for PCG, whose application involves only sparse matrix-vector products.
To ensure positive definiteness, the learned preconditioner is parameterized as
\[
    M^{-1} = GG^T + \epsilon I,
\]
where the sparse matrix $G$ is predicted by the GNN and $\epsilon>0$.
This form is similar to that in~\cite{sappl2019deep,sappl2026convolutional}.

Classical SPAI methods require a predefined sparsity pattern, which can significantly affect preconditioner quality.
\citet{bankestad2024ising} propose a task-specific graph reduction (subsampling) framework based on an Ising model~\cite{cipra1987introduction}, and apply it to learn sparsity patterns for SPAI preconditioners.
By representing the coefficient matrix as a graph, the sparsity pattern selection is formulated as a graph reduction problem.
Specifically, given a prescribed sparsity pattern, the Ising model selects approximately 50\% of its entries to form the sparsity pattern of the preconditioner.
The corresponding nonzero values are then obtained by minimizing $\|AM^{-1}-I\|_F$.

\paragraph{\textbf{MG-, DD-, and deflation-based preconditioners}}
ML can also be used to enhance key components of MG-, DD-, and deflation-based preconditioners, including smoothers, subdomain solvers, transfer operators, and coarse spaces.
Classical MG and DD methods can be employed either as standalone solvers or as preconditioners, and their ML-enhanced counterparts can be used in the same way.
Given the extensive literature on these methods, we discuss them separately in~\Cref{sec:linear-mgdd} to provide a more systematic treatment and avoid redundancy.

\paragraph{\textbf{Green's function-based preconditioners}}
Green's functions (or operators) provide an explicit representation of the inverse of a linear differential operator~\cite{evans2010partial, stakgold2011greens}.
After discretization, the corresponding discrete Green's function represents, or approximates, the inverse of the resulting discrete operator and can therefore be used to construct preconditioners for the associated linear systems~\cite{benzi1999comparative, xu1997preconditioning}.

For simplicity, consider a boundary value problem with homogeneous Dirichlet boundary conditions,
\[
    \left\{
    \begin{aligned}
        L u(x) & = f(x), &  & x \in \Omega,          \\
        u(x)   & = 0,    &  & x \in \partial \Omega,
    \end{aligned}
    \right.
\]
where $u$ denotes the solution function, $L$ is a linear differential operator, $\Omega \subset \mathbb{R}^d$ is the computational domain, and $f$ is the source function.
For a source point $s\in\Omega$, the Green's function $G(\cdot,s)$ is defined by
\begin{equation}\label{equ:green}
    \left\{
    \begin{aligned}
        L_x G(x, s) & = \delta(x - s), &  & x \in \Omega,          \\
        G(x, s)     & = 0,             &  & x \in \partial \Omega,
    \end{aligned}
    \right.
\end{equation}
where $\delta(\cdot)$ denotes the Dirac delta distribution and $L_x$ indicates that the operator acts on the variable $x$.
The solution can then be represented as
\[
    u(x) = \int_{\Omega} G(x, s) f(s) \mathrm{d}s.
\]
The integral operator induced by $G$ maps the source function $f(x)$ to the solution function $u(x)$, and is therefore referred to as the Green's operator.
For the discretized linear system $A\mathbf{u}=\mathbf{b}$, the discrete Green's function accordingly provides a representation or approximation of $A^{-1}$.

In practice, the analytic form of Green's functions is often unavailable, so numerical approximation methods are preferred.
For large-scale problems, however, their numerical approximation remains challenging.
First, the Green's function is defined on $\Omega \times \Omega$, resulting in a $2d$-dimensional domain.
Second, the Dirac delta distribution $\delta(x-s)$ introduces a singularity at $x=s$, leading to highly localized and non-smooth behavior nearby.
To address these issues, ML models, particularly neural networks $\mathcal{N}$ that take $(x, s)$ as input, have been used to approximate the Green's function $G$ and thereby construct Green's function-based preconditioners.
Such learned Green's functions (or operators) are referred to as neural Green's functions (or operators) in~\cite{tang2022neural, melchers2026neural}.
In the following, we focus on their use for preconditioner construction rather than providing a comprehensive survey of Green's function learning.

\citet{ichimura2020fast} construct a preconditioner by learning local Green's functions for linear systems arising from the discretization of 3D wave equations in heterogeneous solid media.
For each spatial point, the local displacement satisfies
\[
    u_i = \sum_{j=1}^3 G_{ij} f_j, \qquad i=1,2,3,
\]
where $G_{ij}$ denotes the Green's function relating the $j$-th force component to the $i$-th displacement component.
To simplify the learning task, the Green's functions are parameterized around baseline functions as
\[
    \left\{
    \begin{aligned}
        G_{ii}^{\mathrm{NN}}
         & = G_{ii}^{\mathrm{base}}(c_1+c_2x+c_3y+c_4z),
        \\
        G_{ij}^{\mathrm{NN}}
         & = G_{ij}^{\mathrm{base}}(c_5+c_6x+c_7y+c_8z), \quad i\neq j,
    \end{aligned}
    \right.
\]
where $x$, $y$, and $z$ are spatial coordinates, and $G_{ij}^{\mathrm{base}}$ is the Green's function for a specific material configuration.
Rather than predicting the full Green's functions, an FNN maps local material properties to the low-dimensional coefficients $(c_1,\ldots,c_8)$.
The reference Green's functions used for training are obtained numerically.
The resulting local Green's-function approximation serves as a preconditioner for PCG.

\citet{sun2026learning} address the singularity of the Green's function by introducing an augmented variable $z=\phi(x,s)\in\mathbb{R}$, yielding
\[
    G(x,s)=\hat G(x,s,\phi(x,s)),
\]
where $\hat G$ is referred to as a singularity-encoded Green's function, and is approximated by an FNN with a smooth activation.
The networks are trained in an unsupervised way, combining a governing equation loss, a boundary condition loss, a normalization condition loss, and a symmetry loss.
After training, the learned Green's function is projected back to the physical domain by exploiting its symmetry:
\[
    G_\theta(x,s) = \frac{1}{2} \hat G(x,s,\phi(x,s)) + \frac{1}{2} \hat G(s,x,\phi(s,x)).
\]
Since the preconditioner constructed from the learned Green's function is dense, this approach is suitable only for small-scale systems.
For large-scale problems, the authors instead use it as the coarse-grid solver in a two-level overlapping Schwarz preconditioner; see~\Cref{sec:linear-mgdd-smoother}.
Moreover, the resulting Green's function-based preconditioner is combined with damped Jacobi or geometric multigrid to form an alternating HIM; see~\Cref{sec:linear-precond-alternating}.

Building on the multiscale neural network (MSNN) proposed in~\cite{hao2026multiscale}, \citet{xu2025neural} introduce an adaptive MSNN ($\alpha$MSNN) for learning Green's functions.
They regularize the Dirac delta in~\Cref{equ:green} using a Gaussian kernel with a small width $\varepsilon>0$,
\[
    \delta(x - s) \approx f_{\varepsilon}(x,s)
    = \left(\frac{1}{\varepsilon\sqrt{\pi}}\right)^d \exp\left( -\frac{\|x-s\|_2^2}{\varepsilon^2} \right).
\]
The Green's function is then represented as
\[
    G_{\alpha\mathrm{MSNN}}(x,s;\theta)
    = \sum_{j=1}^{L} \varepsilon^{\alpha_j(s)} \hat{f}_j
    \left(\frac{x-s}{\varepsilon^{\beta_j(s)}}, s; \theta_j \right),
\]
where $\hat f_j$ denotes the $j$-th subnetwork parameterized by $\theta_j$, and the trainable source-dependent parameters $\alpha_j(s)$ and $\beta_j(s)$ control its amplitude and spatial scale, respectively.
This decomposition is designed to capture the multiscale structure of the Green's function through localized near-field (\ie, $x\approx s$), broader middle-field, and smooth far-field components (\ie, large $|x-s|$).
In particular, the near-field branch ($j=1$) uses $\beta_1=1$ to match the width of the regularized source, and the far-field branch ($j=L$) uses $\alpha_L = 0$ and $\beta_L=0$ to represent the smooth global behavior.
For $1<j<L$, the middle-field branches use $0 < \beta_j < 1$ to capture residual structures at intermediate spatial scales.
The learned Green's function is then used to construct either an SPAI preconditioner, when it is strongly localized, or a hierarchical matrix preconditioner, when global interactions are significant.
The resulting preconditioner is employed within preconditioned GMRES.

\paragraph{\textbf{Other preconditioners}}
\citet{pasetto2017reduced} consider transient diffusion equations and employ POD to construct a ROM-based preconditioner for PCG.
At each PCG iteration, the preconditioning step approximately solves $A \mathbf{z} = \mathbf{r}$, where $\mathbf{r}$ is the current residual.
Given a reduced POD basis $P\in\mathbb{R}^{n\times m}$, the coefficient matrix $A$ is projected onto a low-dimensional subspace,
\[
    \widetilde A = P^TAP \in \mathbb{R}^{m\times m}.
\]
This yields a low-rank approximation $A^{-1} \approx P\widetilde A^{-1}P^T$.
However, since $P\widetilde A^{-1}P^T$ has rank $m<n$, using it directly as a preconditioner would restrict the Krylov subspace to at most $m$ dimensions.
To avoid this issue, the authors first apply a smoother $S$, such as Jacobi, to obtain the correction $S\mathbf{r}$.
The remaining residual, $(I-AS)\mathbf{r}$, is then projected onto the reduced space by $P^T$, corrected using $\widetilde A^{-1}$, and prolonged back to the full space by $P$.
The resulting preconditioner is
\[
    M^{-1} = S+P\widetilde A^{-1}P^T(I-AS).
\]
The POD basis $P$ is constructed from snapshots of the full-order model.
The method can also be interpreted as a two-grid algebraic multigrid scheme, in which the restriction and prolongation operators are given by the POD projections $P^T$ and $P$, respectively.
In this sense, the approach can also be viewed as a form of POD-based transfer operator learning; see~\Cref{sec:linear-mgdd-transfer_operator}.

Efficient block Jacobi preconditioners require identifying strongly connected variables; however, classical block pattern generation, such as graph-based clustering, can be computationally expensive and poorly scalable.
\citet{gotz2018machine} employ a CNN to predict the diagonal block pattern of a block-Jacobi preconditioner from the matrix sparsity pattern.
Specifically, the sparsity pattern is interpreted as an image, and the detection of diagonal blocks is formulated as a multi-label classification problem.
For a matrix $A\in\mathbb{R}^{n\times n}$, the CNN takes a band around the main diagonal as input and outputs $\hat y=(\hat y_i)_{i=1}^n\in[0,1]^n$, where $\hat y_i$ is the predicted probability that row $i$ starts a diagonal block.
The corresponding labels are $y=(y_i)_{i=1}^n\in\{0,1\}^n$, with $y_i=1$ for block starts and $y_i=0$ otherwise.
The network is trained by minimizing the binary cross-entropy
\[
    \mathcal{L}(y,\hat{y})
    = -\sum_{i=1}^{n} \left[ y_{i}\log(\hat{y}_i) + (1-y_i)\log(1-\hat{y}_i) \right].
\]
Given the predicted block partition, $D=\operatorname{blockdiag}(D_1,\ldots,D_N)$, the block-Jacobi preconditioner is constructed as
\[
    M^{-1}=D^{-1}
    =\operatorname{blockdiag} \left(D_1^{-1},\ldots,D_N^{-1} \right).
\]

\citet{trappler2025statedependent} focus on variational data assimilation, where the analysis step typically leads to a high-dimensional nonlinear least-squares problem that is often solved using a Gauss-Newton method.
At each iteration, the resulting linearized system is SPD and depends on the current state $s$, \ie, $A(s) \mathbf{x}= \mathbf{b}(s)$, which can be solved by PCG with a state-dependent preconditioner $M(s)^{-1} = L(s)L(s)^T$.
Given the leading eigenpairs of the eigendecomposition $A(s) \approx U(s)\Lambda(s)U(s)^T$, spectral preconditioners~\cite{diouane2024efficient} construct
\[
    L(s) = I + U(s) \left( \mu \Lambda(s)^{-1/2}-I \right) U(s)^T,
\]
where $U(s)\in\mathbb{R}^{n\times r}$, $\Lambda(s)\in\mathbb{R}^{r\times r}$, and $\mu >0$ is a parameter.
The authors employ a U-Net to learn the mapping from the state $s$ to approximations of the leading eigenvectors and eigenvalues, denoted by $U_\theta(s)$ and $\Lambda_\theta(s)$, respectively.
A QR decomposition is applied to the predicted eigenvectors to enforce orthonormality, and a scaled sigmoid function bounds the predicted eigenvalues between zero and a prescribed upper bound.
The network is trained by minimizing the low-rank approximation error
\[
    \mathcal{L}_\theta(s) = \|A_\theta(s) - A(s)\|_F^2 \approx \frac{1}{J_\text{max}} \sum_j^{J_\text{max}} \|A_\theta(s)\mathbf{w}_j - A(s)\mathbf{w}_j\|_2^2,
\]
where $A_\theta(s) = U_\theta(s)\Lambda_\theta(s)U_\theta(s)^T$, and $\mathbf{w}_j$ is a Gaussian random vector.

\paragraph{\textbf{Discussion}}

As we mentioned in~\Cref{sec:linear-core_components-param}, parameters in preconditioners can be selected within the framework of \emph{automated solver selection}; see~\Cref{sec:linear-solver_selection}.
Examples include learning parameters for ILU-type preconditioners using classification models, as explored in~\cite{xu2005study,xu2005data,holloway2007neural} and~\cite[Section~4]{moore2022machine}, or using reinforcement learning~\cite{kuefler2008using}.

Approaches for \emph{enhancing classical preconditioners} combine ML models with classical preconditioning frameworks while often preserving the linear structure of the resulting operator.
In such approaches, ML models are employed to predict decomposition factors, sparsity patterns, or scalar parameters in classical preconditioners.
As a result, the learned preconditioners can be applied within a wide range of iterative methods, including standard Krylov subspace methods.

General remarks on existing methods and future directions for \emph{enhancing classical preconditioners} are summarized as follows:
\begin{enumerate}
    \item \textbf{Comparison with \emph{learning the action of the inverse}~(\Cref{sec:linear-precond-direct}):}
          \begin{itemize}
              \item Approaches based on \emph{learning the action of the inverse} generally yield nonlinear preconditioners unless a linear model architecture is used, whereas most methods for \emph{enhancing classical preconditioners} retain a linear preconditioning operator.
              \item Consequently, the former typically require flexible Krylov methods, while the latter can be combined with standard preconditioned Krylov methods.
                    This can be an advantage of \emph{enhancing classical preconditioners}.
                    As discussed in~\Cref{sec:linear-precond-direct}, PCG has lower per-iteration computational and memory costs than FCG, whereas preconditioned GMRES and FGMRES have similar per-iteration costs, with preconditioned GMRES requiring slightly less storage.
              \item Retaining a linear preconditioner also facilitates theoretical analysis, for example, through the spectrum or condition number of the preconditioned system.
                    Such analysis is generally more difficult for nonlinear learned preconditioners and remains relatively limited.
          \end{itemize}
    \item \textbf{Parallel construction and application:}
          The factors of IC- or ILU-based preconditioners can be predicted by ML models such as GNNs in a highly parallel manner~\cite{li2023learning,hausner2024neural}, avoiding some of the sequential dependencies arising in classical incomplete factorization.
          This is attractive for large-scale computations on modern parallel hardware.
          However, even when their construction is parallelized, the resulting IC- or ILU-based preconditioners still require triangular solves during application, which are themselves difficult to parallelize efficiently.
          In this respect, learned SPAI-based preconditioners are well suited to parallel architectures because they can be applied through sparse matrix-vector products.
    \item \textbf{Flexible sparsity patterns:}
          Most existing approaches for learning IC- or ILU-based factors adopt a prescribed sparsity pattern, often inherited from the coefficient matrix or a classical incomplete factorization.
          Although this controls memory and computational costs, it may restrict the preconditioner quality.
          Learning or adaptively selecting sparsity patterns therefore represents a promising direction, with the challenge of balancing approximation quality, fill-in, and computational cost.
    \item \textbf{Learning correction terms:}
          In~\cite{li2024deep,trifonov2025learning}, ML models are used to learn corrections to standard IC factors rather than constructing the preconditioner entirely from scratch.
          Such hybridization retains a reliable classical baseline while reducing the burden on the learning model.
          This idea could be extended to other classical preconditioners, such as ILU and SPAI, by learning corrections to their factors, entries, or sparsity patterns.
    \item \textbf{Offline cost and generalization:}
          Methods of \emph{enhancing classical preconditioners} also require an offline training stage and their performance may depend on how closely test matrices resemble the training distribution.
          Their computational benefit is therefore most compelling for sequences of related linear systems, where the training cost can be amortized over many solves when the systems are sufficiently similar.
          Achieving effective training and improving generalization across matrix sizes, discretizations, and problem parameters remain important challenges.
    \item \textbf{Loss functions:}
          Designing effective loss functions remains a key challenge in \emph{preconditioner learning}.
          For CG-type methods applied to SPD systems, one natural objective is to reduce the spectral condition number of the preconditioned system, which enters a standard convergence bound; spectral clustering also matters.
          In practice, however, such spectral quantities are often expensive to evaluate, so many works rely on proxy objectives such as $\|I - AM^{-1}\|_F^2$.
          Such proxies may not accurately reflect the convergence behavior of the iterative solver, motivating loss functions that more directly target spectral properties.
\end{enumerate}

\subsubsection{Alternating hybrid iterative methods}\label{sec:linear-precond-alternating}

\emph{Alternating hybrid iterative methods} combine a classical iterative scheme with an ML-based iterative scheme in an alternating manner.
Similar to \emph{initial guess learning}~(\Cref{sec:linear-core_components-initial}), their motivation can also be considered from two perspectives:
\begin{enumerate}
    \item From the view of classical iterative methods, stationary (relaxation) methods often converge slowly because they are inefficient at reducing low-frequency errors.
          Standard (unpreconditioned) Krylov methods may also converge slowly for ill-conditioned systems.
    \item From the view of ML, the spectral bias of neural network training makes it more difficult to learn local high-frequency features accurately.
          Stationary methods, in contrast, are effective at reducing the remaining errors and can therefore refine the learned solution.
\end{enumerate}
Within this framework, the classical component may consist of a stationary iterative method or a Krylov subspace method, while the ML-based component is typically an ML-preconditioned Richardson or Krylov method.
\Cref{fig:precond}(d) shows one example of an alternating HIM, where $n_c$ iterations of a classical relaxation method and $n_m$ iterations of an ML-preconditioned Richardson method are applied alternately.
We regard \emph{alternating HIMs} as a usage strategy for \emph{preconditioner learning}, since preconditioners learned using the approaches in~\Cref{sec:linear-precond-direct,sec:linear-precond-hybrid} can be incorporated with classical iterative methods in an alternating manner.
Several existing approaches for \emph{alternating hybrid iterative methods} are summarized in~\Cref{tab:linear-precond-alternating-comparison}.
To organize the existing literature, we discuss \emph{alternating HIMs} in three groups: two representative families and other approaches.
\textbf{(1)} the hybrid, iterative, numerical, and transferable solver (HINTS)~\cite{zhang2024blending} and related variants;
\textbf{(2)} the Fourier neural solver (FNS)~\cite{cui2022fourier, cui2025hybrid} and related variants; and
\textbf{(3)} other approaches.

\begin{table}[!ht]
    \renewcommand{\arraystretch}{1.3}
    \centering
    \caption{
        Comparison of existing approaches for \emph{alternating hybrid iterative methods}~(\Cref{sec:linear-precond-alternating}).
        The last column shows the type of learned preconditioner: \emph{learning the action of $A^{-1}$}~(\Cref{sec:linear-precond-direct}) or \emph{enhancing classical preconditioners}~(\Cref{sec:linear-precond-hybrid}).
        For the latter, the underlying classical preconditioner is given in parentheses.
        In~\cite{wu2026are}, in addition to Jacobi and ML-preconditioned iterations, a physics-aware Anderson acceleration (PA-AA) strategy is introduced.
    }
    \label{tab:linear-precond-alternating-comparison}
    \resizebox{\linewidth}{!}{%
        \begin{minipage}{\linewidth}
            \centering
            \begin{tabular}{
                    p{1.0cm}llll
                }
                \toprule
                \textbf{ref.}
                 & \textbf{tested problem}
                 & \textbf{classical iter.}
                 & \textbf{ML-preconditioned iter.}
                 & \textbf{learned prec. type}
                \\
                \midrule
                \cite{zhang2024blending}
                 & Poisson and Helmholtz eq.
                 & Jacobi, Gauss-Seidel
                 & DeepONet-prec. Richardson
                 & learning the action of $A^{-1}$
                \\
                \cite{kahana2023geometry}
                 & Darcy flow and linear elasticity
                 & Gauss-Seidel
                 & DeepONet-prec. Richardson
                 & learning the action of $A^{-1}$
                \\
                \cite{versano2024attentionbased}
                 & Helmholtz eq.
                 & Gauss-Seidel, GMRES
                 & DeepONet-prec. Richardson
                 & learning the action of $A^{-1}$
                \\
                \cite{bai2026hybrid}
                 & Poisson eq. from incompressible flow
                 & CG or PCG
                 & DeepONet-prec. Richardson (DLSM)
                 & learning the action of $A^{-1}$
                \\
                \cite{kopanicakova2025deeponet}
                 & diffusion and Helmholtz eq.
                 & Jacobi
                 & DeepONet-prec. Richardson/Krylov
                 & enhancing class. (subspace corr.)
                \\
                \multirow{2}{*}{\cite{lee2025fast}}
                 & \multirow{2}{*}{Helmholtz scattering}
                 & Stationary methods
                 & DeepONet-prec. Richardson
                 & \multirow{2}{*}{enhancing class. (subspace corr.)}
                \\
                 &
                 & Krylov
                 & DeepONet-prec. Krylov
                 &
                \\
                \cite{wang2024labelfree}
                 & Poisson, Helmholtz, and diffusion eq.
                 & Jacobi
                 & FNO-prec. Richardson
                 & learning the action of $A^{-1}$
                \\
                \cite{hu2025hybrid}
                 & Poisson eq.
                 & Richardson, GS, multigrid
                 & MIONet-prec. Richardson
                 & learning the action of $A^{-1}$
                \\
                \cite{cui2022fourier,cui2025hybrid}
                 & diffusion, convection-diffusion, Helmholtz
                 & Jacobi, Cheby-semi, Krylov
                 & FFT \& FNN-prec. Richardson
                 & enhancing class. (spectral)
                \\
                \cite{rizzuti2019learned}
                 & Helmholtz equation
                 & CG on normal eq.
                 & CNN-based correction
                 & learning the correction
                \\
                \cite{wu2026are}
                 & diffusion and Helmholtz eq.
                 & Jacobi
                 & same as~\cite{zhang2024blending} or \cite{cui2022fourier,cui2025hybrid}
                 & same as~\cite{zhang2024blending} or \cite{cui2022fourier,cui2025hybrid}
                \\
                \cite{sun2026learning}
                 & Poisson and Helmholtz eq.
                 & Jacobi
                 & learned Green's func.-prec. Richardson
                 & enhancing class. (Green's fun.)
                \\
                \bottomrule
            \end{tabular}
        \end{minipage}
    }
\end{table}

\paragraph{\textbf{HINTS and related variants}}

\citet{zhang2024blending} propose the HINTS method, which combines a DeepONet-preconditioned Richardson iteration with classical relaxation methods such as Jacobi and Gauss-Seidel for parametric linear PDE systems.
For a linear PDE $L_x(\mathbf{u}(x); \mathbf{c}(x))=\mathbf{f}(x)$, $x\in\Omega$, DeepONet is trained offline to approximate the solution operator $(\mathbf{c},\mathbf{f})\to\mathbf{u}$, where $\mathbf{c}$ denotes the PDE parameters.
Within HINTS, the trained network $\mathcal{N}$ acts as a preconditioner by mapping the current residual $\mathbf{r}_k=\mathbf{b}-A\mathbf{x}_k$ and $\mathbf{c}$ to an update direction,
\[
    \mathbf{x}_{k+1} = \mathbf{x}_k + \mathbf{d}_k,
    \qquad
    \mathbf{d}_k = \mathcal{N}(\mathbf{r}_k, \mathbf{c}).
\]
By alternating the learned and classical iterations, HINTS may achieve nearly uniform convergence across eigenmodes.
With the transferability of DeepONet, a model trained on coarse grids can be applied to much larger systems, although its performance may deteriorate under distribution shifts.
HINTS can also be incorporated into multigrid as a smoother or used as a preconditioner for FGMRES.

Moreover, several variants of HINTS have been proposed to further enhance its performance.
\citet{kahana2023geometry} investigate and improve the geometry transferability of HINTS by incorporating operator-level transfer learning~\cite{goswami2022deep}.
In this setting, a DeepONet is first trained on a source geometry with sufficient data and is then adapted to a related target geometry using only a small amount of target-domain data.
Specifically, the target DeepONet is initialized with the parameters learned on the source geometry, and only selected layers are fine-tuned for the target geometry.
Analogously, \citet{versano2024attentionbased} extend HINTS by incorporating a geometry-aware DeepONet with masked self-attention layers~\cite{vaswani2017attention}.
The masked self-attention layers filter out points outside the computational domain, avoiding artificial extensions of the input beyond the physical domain.
The model is trained on simple reference geometries and can be applied to different geometries without additional fine-tuning, thereby enabling geometry transferability.

In~\cite{bai2026hybrid}, the authors propose an ML-based iterative method called DLSM that combines DeepONet with a line-search procedure, where DeepONet predicts the search direction; see~\Cref{sec:linear-precond-direct} for more details.
Building on DLSM, they further propose an alternating HIM, HyDEA (Hybrid Deep lEarning line-search directions and iterative methods for Accelerated solutions), which alternates DLSM with a classical CG-type method to solve the pressure Poisson equation arising in viscous incompressible flow simulations.

While HINTS and its variants primarily use DeepONet to directly approximate a preconditioner, \citet{kopanicakova2025deeponet} instead construct preconditioners within a subspace correction framework for parametric linear systems arising from PDE discretizations.
They propose two approaches: direct preconditioning (DP) and trunk basis (TB); see~\Cref{sec:linear-mgdd-transfer_operator} for more details.
Preconditioners based on these two approaches are combined with classical preconditioners to form multiplicative or additive preconditioners for Krylov subspace methods.

Based on the DP and TB preconditioning strategies~\cite{kopanicakova2025deeponet}, \citet{lee2025fast} develop relaxation-based and Krylov-based meta-solvers for 3D Helmholtz scattering problems, for which standard Jacobi iterations may diverge. The relaxation-based meta-solver combines DP with a classical relaxation scheme, while the Krylov-based meta-solver combines TB with a Krylov subspace method. The DeepONet is trained on non-scattering problems and then used without retraining for scattering problems with different scatterer geometries.

Training the DeepONet in HINTS typically requires numerical solutions as labels for different parameter settings.
To avoid this cost, \citet{wang2024labelfree} propose a label-free variant of HINTS, in which the DeepONet is replaced by an FNO trained using a physics-informed loss~\cite{li2024physicsinformed}.
The loss penalizes the PDE residual and can therefore be viewed as an unsupervised, physics-based training objective.
\citet{hu2025hybrid} develop an alternating HIM that replaces the DeepONet in HINTS with a multiple-input operator network (MIONet)~\cite{jin2022mioneta}.
MIONet can separately model linear and nonlinear dependencies on multiple PDE inputs; for example, for the Poisson equation, the solution operator is nonlinear in the diffusion coefficient but linear in the source term.
The authors also provide a theoretical analysis of the resulting alternating HIM, including its convergence conditions, spectral behavior, and convergence rate.

\paragraph{\textbf{FNS and related variants}}

\citet{cui2022fourier} propose the FNS method, which combines stationary iterations with a learned correction in Fourier space that can be interpreted as an ML-preconditioned Richardson iteration.
From a spectral perspective, consider an eigendecomposition
\[
    A=Q\Lambda Q^{-1},
    \qquad
    A^{-1}=Q\Lambda^{-1}Q^{-1}.
\]
The learned correction in FNS can be interpreted as an approximation of $A^{-1}$, where the eigenvector matrix $Q$ is approximated by a Fourier basis via fast Fourier transform (FFT), and the inverse eigenvalues $\Lambda^{-1}$ are learned from the PDE parameters $\mathbf c$ using a network $\mathcal{N}$, \ie,
\[
    Q\approx\mathcal F,
    \qquad
    \operatorname{diag}(\Lambda^{-1})\approx\mathcal N(\mathbf c).
\]
The resulting learned preconditioner remains linear, in contrast to nonlinear operators such as DeepONet used in HINTS~\cite{zhang2024blending}.

Subsequently, \citet{cui2025hybrid} propose an improved FNS.
Rather than directly approximating the eigenvector matrix $Q$ by the Fourier basis $\mathcal{F}$ as in~\cite{cui2022fourier}, they introduce a transition operator $T$ such that
\[
    Q=\mathcal{F}T,
\]
where $T$ is a learnable transition matrix and is approximated by a linear convolutional operator.
Furthermore, FNS has been extended to solve PDEs on unstructured grids in~\cite{liu2026nodal} through the use of GNNs.

\paragraph{\textbf{Other approaches}}

\citet{rizzuti2019learned} propose an iterative solver for the Helmholtz equation that alternates classical Krylov iterations with CNN-based corrections.
In this approach, a few steps of the Krylov method are first performed, followed by a CNN that refines the intermediate solution.
The CNN takes the Krylov output as input and produces an updated solution, which is then passed to the subsequent Krylov iterations for further refinement.
The CNN architecture incorporates a sequence of ResNet modules to effectively correct the iterate, which can be interpreted as an ML-preconditioned Richardson iteration method.

The reliability and robustness of alternating HIMs are investigated in~\cite{wu2026are}, with a focus on HINTS~\cite{zhang2024blending} and FNS~\cite{cui2022fourier}.
The study compares loss strategies based on the update $\mathbf{d}_k = \mathbf{x}_{k+1} - \mathbf{x}_k$ and on the residual $\mathbf{r}_k = \mathbf{b} - A\mathbf{x}_k$, and also considers unrolling in the training process of the alternating iteration.
An iteration based on Anderson acceleration (AA)~\cite{walker2011anderson} is further introduced to accelerate the alternating fixed-point scheme.
Overall, the results indicate that the loss choice is crucial and approach-dependent. Since the updates may become small even when the residual remains large, physics-based residual losses can improve performance.
In~\cite{sun2026learning}, the neural Green's function-based preconditioner, discussed in~\Cref{sec:linear-precond-hybrid}, is combined with classical iterative methods, such as damped Jacobi and geometric multigrid, to form an alternating HIM.

\paragraph{\textbf{Discussion}}

Since \Cref{sec:linear-precond-alternating} (\emph{alternating HIMs}) discusses a framework for combining classical and ML-based iterative schemes, the approaches described in~\Cref{sec:linear-precond-direct} (\emph{learning the action of the inverse}) and \Cref{sec:linear-precond-hybrid} (\emph{enhancing classical preconditioners}) can naturally be incorporated into this framework, as demonstrated in~\cite{bai2026hybrid, kopanicakova2025deeponet, sun2026learning}.
We have the following remarks on existing studies and future directions for \emph{alternating HIMs}:

\begin{enumerate}
    \item \textbf{Training data generation and training cost}:
          Most ML models in \emph{alternating HIMs} are trained using supervised data, whose generation typically requires classical PDE solves; an exception is the physics-informed label-free approach in~\cite{wang2024labelfree}.
          Due to the offline data generation and training costs, \emph{alternating HIMs} are particularly suitable for multi-query settings, where these costs can be amortized over repeated solves.
          Their computational benefits also tend to be more obvious for large-scale systems, since some models can be trained on smaller discretizations and transferred to larger ones.
          For small systems, the ML inference overhead may offset the benefits.
    \item \textbf{Effect of $n_c$ and $n_m$}:
          $n_c$ and $n_m$ denote the numbers of classical and ML-preconditioned iterations, respectively, within each alternating cycle; see also~\Cref{fig:precond}(d).
          Most existing studies adopt $n_c>1$ and $n_m=1$.
          In~\cite{hu2025hybrid}, the authors theoretically analyze the convergence rates of methods that alternate Gauss-Seidel iterations with neural operator-preconditioned Richardson iterations.
          The results show that the choice of $n_c$ can strongly affect convergence: too small an $n_c$ may lead to divergence, whereas increasing $n_c$ initially improves the acceleration but eventually diminishes it.
    \item \textbf{Stability}:
          In \emph{alternating HIMs}, learned preconditioners can be nonlinear, as in HINTS, or linear, as in FNS.
          Nonlinear preconditioners may yield rapid residual reduction and substantial speedup in the early iterations, but can exhibit stagnation or even residual growth at later stages~\cite{wu2026are}.
          Compared with nonlinear learned mappings, linear preconditioners are generally more amenable to stability and convergence analysis and can be designed to satisfy structural or spectral properties required by classical iterative schemes.
          Developing strategies to mitigate such stagnation, including appropriately designed learned operators with improved stability properties, remains an important direction.
    \item \textbf{Extension}:
          Since \emph{alternating HIMs} can be regarded as a usage strategy of learned preconditioners, most ML-based preconditioners can, in principle, be incorporated into this framework.
          Moreover, existing studies focus mainly on linear systems, and the framework may also be extended to nonlinear systems.
\end{enumerate}

\subsection{Enhancements of subspace correction methods}\label{sec:linear-mgdd}
Based on the principle of \emph{subspace correction}, multigrid (MG), domain decomposition (DD), deflation, and augmentation methods can be used to solve large-scale linear systems efficiently, in which the error is reduced through successive corrections on appropriately chosen subspaces.
From a high-level perspective, these methods share a common conceptual framework for solving linear systems, since they iteratively reduce the error and improve the solution through subspace corrections.
However, their effectiveness depends critically on the quality of their components, which are usually designed or selected based on prior knowledge and problem-specific expertise.
Recent studies have explored the use of ML techniques to automatically learn or optimize these components from data.
The structure of this subsection and the topics covered are summarized in~\Cref{tab:mgdd-structure}.
Before reviewing these ML-enhanced approaches, we first provide the necessary preliminaries on classical multigrid, domain decomposition, deflation, and augmentation methods in~\Cref{sec:linear-mgdd-prelims}.

\begin{table}[!ht]
    \caption{
        Outline of topics in~\Cref{sec:linear-mgdd} on \emph{enhancements of subspace correction methods}.
        The symbol ``--'' indicates that the component is not applicable to the corresponding method, and ``\xmark'' indicates that we found no related work on ML enhancement in the reviewed literature.
    }
    \label{tab:mgdd-structure}
    \begin{tabular}{lllll}
        \toprule
                                                           & \textbf{multigrid}
                                                           & \textbf{domain decomposition}
                                                           & \textbf{multiscale}
                                                           & \textbf{deflation/augmentation}
        \\ \midrule
        \multirow{2}{*}{\secref{sec:linear-mgdd-smoother}} & smoothers                                                                                      & subdomain solvers & \multirow{2}{*}{--} & \multirow{2}{*}{--} \\ \cmidrule(lr){2-3} & coarse-grid solvers & coarse-grid solvers & & \\ \midrule
        \secref{sec:linear-mgdd-coarsening}                & coarsening in AMG                                                                              & --                & --                  & --                  \\ \midrule
        \secref{sec:linear-mgdd-transfer_operator}         & \multicolumn{4}{l}{transfer operators / coarse spaces / deflation (or augmentation) subspaces}                                                                 \\ \midrule
        \secref{sec:linear-mgdd-parameters}                & parameters                                                                                     & parameters        & \xmark              & \xmark              \\ \bottomrule
    \end{tabular}
\end{table}

Since these methods can be used as preconditioners, most approaches reviewed in this subsection can also be interpreted as \emph{preconditioner learning}; see~\Cref{sec:linear-precond}.
More specifically, ML is used here to construct the components of these methods, so they are instances of \emph{enhancing classical preconditioners}; see~\Cref{sec:linear-precond-hybrid}.

\subsubsection{Preliminaries}\label{sec:linear-mgdd-prelims}
When solving large-scale linear systems arising from PDE discretizations, classical iterative methods applied directly to the whole system may converge slowly.
For instance, many stationary iterative methods efficiently reduce high-frequency error components but perform poorly on low-frequency errors, often requiring many iterations to achieve convergence.
A common remedy is provided by \emph{subspace correction methods}~\cite{xu1992iterative,tang2009comparison}, including multigrid~\cite{hackbusch_multi-grid_1985,wesseling_introduction_1995,trottenberg2001multigrid,brandt_multigrid_2011}, domain decomposition~\cite{smith_domain_2004,toselli_domain_2005,dolean_introduction_2015,quarteroni_domain_1999}, deflation~\cite{nicolaides1987deflation, morgan1995restarted}, and augmentation methods~\cite{frank2001construction, parks2006recycling}.
The main idea is to accelerate convergence by performing corrections associated with suitably chosen subspaces, where error components that are difficult to eliminate by global iterations alone can be treated more efficiently.

Instead of solving only the original linear system, these methods combine global updates with solutions of smaller problems defined on suitably chosen subspaces.
The resulting subspace correction is then incorporated into the global solution.
Within this framework, MG, DD, deflation, and augmentation methods can significantly improve convergence compared with classical iterative solvers.
For simplicity, we focus on two-grid MG methods and two-level DD methods, although the underlying ideas can be naturally extended to multilevel variants.
We discuss them below.

\paragraph{\textbf{MG methods}}
In MG methods, corrections are performed on a hierarchy of coarse-grid spaces, where errors are reduced through successive corrections on multiple levels.
Smoothing methods (or smoothers) applied on fine grids efficiently damp high-frequency error components, while coarse-grid corrections eliminate low-frequency errors that cannot be effectively reduced on the fine grid alone.
The solvers for coarse-grid systems are called coarse-grid solvers.

MG methods are commonly classified into geometric multigrid (GMG) and algebraic multigrid (AMG).
GMG methods construct coarse grids using geometric information from the underlying problem, whereas AMG methods construct coarse spaces algebraically from the system matrix.
The selection of coarse grids is typically achieved through a coarsening strategy, also referred to as coarse/fine (C/F) splitting~\cite{ruge1987algebraic}.
Transfer operators, namely restriction and prolongation, are used to construct the coarse-grid system and transfer information between grid levels.

\begin{figure}[!htbp]
    \centering

    \begin{subfigure}[t]{0.48\textwidth}
        \vspace{0pt}
        \centering
        \includegraphics[width=0.7\textwidth]{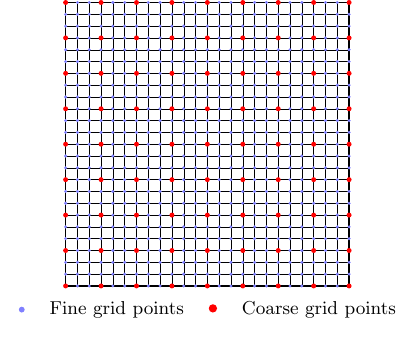}
        \caption{
            Fine- and coarse-grid points in a two-level MG hierarchy.
            Blue and red points represent fine- and coarse-grid points, respectively.
            Restriction and prolongation operators transfer information between the fine and coarse grids, where smoothing and coarse-grid correction act on different error components.
        }
        \label{fig:prelims-mgdd_grid-mg}
    \end{subfigure}
    \hfill
    \begin{subfigure}[t]{0.48\textwidth}
        \vspace{0pt}
        \centering
        \includegraphics[width=0.7\textwidth]{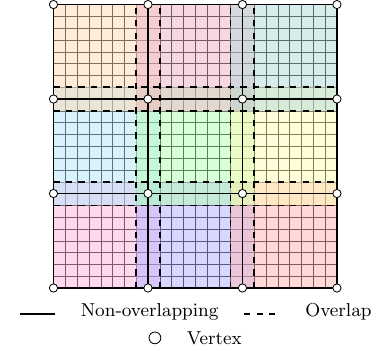}
        \caption{
            Overlapping subdomains in DD.
            Nine overlapping subdomains (colored regions) define local correction subspaces; solid lines mark non-overlapping interfaces and dashed lines indicate overlap.
            An exemplary coarse space is associated with the vertices of the non-overlapping partition.
        }
        \label{fig:prelims-mgdd_grid-dd}
    \end{subfigure}
    \caption{
        Subspace correction structures in \textbf{(a)} MG and \textbf{(b)} DD methods.
        MG methods combine fine-grid smoothing with coarse-grid correction, whereas DD methods combine local subdomain corrections with a coarse space correction.
        Both approaches accelerate convergence through corrections on suitably constructed subspaces.
    }
    \label{fig:prelims-mgdd_grid}
\end{figure}

We consider a two-grid V-cycle method~\cite{trottenberg2001multigrid} as an example.
\Cref{fig:prelims-mgdd_grid-mg} illustrates the fine and coarse grids used in this setting.
Given the $k$-th iterate $\mathbf{x}_k$ and a prolongation operator $P\in\mathbb{R}^{n\times n_c}$, where $n_c$ denotes the dimension of the coarse-grid system, a two-grid iteration is summarized in~\Cref{alg:mg}.
\begin{center}
    \begin{minipage}{0.75\textwidth}
        \begin{algorithm}[H]
            \caption{Two-grid V-cycle}\label{alg:mg}
            \begin{algorithmic}[1]
                \State $\mathbf{x}_{k}' = \mathbf{x}_{k} + M^{-1}(\mathbf{b} - A\mathbf{x}_{k})$ \Comment{Pre-smoothing (possibly multiple iterations)}
                \State $\mathbf{r}_{k}' = \mathbf{b} - A\mathbf{x}_{k}'$ \Comment{Compute the residual}
                \State $A_c = P^T A P$, \quad $\mathbf{r}_c = P^T \mathbf{r}_{k}'$ \Comment{Restriction to the coarse grid}
                \State Solve $A_c \mathbf{e}_c = \mathbf{r}_c$ \Comment{Solve the coarse-grid system}
                \State $\mathbf{x}_{k+1}' = \mathbf{x}_k' + P \mathbf{e}_c$ \Comment{Prolongation and correction update}
                \State $\mathbf{x}_{k+1} = \mathbf{x}_{k+1}' + M^{-1}(\mathbf{b} - A\mathbf{x}_{k+1}')$ \Comment{Post-smoothing}
            \end{algorithmic}
        \end{algorithm}
    \end{minipage}
\end{center}

A multilevel V-cycle can be obtained by recursively replacing the coarse-grid correction step with another two-grid V-cycle until the coarsest level is reached.
In addition to the standard V-cycle, other cycling strategies have also been developed, including the W-cycle, F-cycle, and full multigrid (FMG) method~\cite{brandt1977multi}.
We refer readers to~\cite{trottenberg2001multigrid,brandt_multigrid_2011} for further details on multigrid methods.

\paragraph{\textbf{DD methods}}

DD methods can be categorized into one-level and two-level methods.
One-level DD methods partition the global problem into smaller subproblems, apply local solvers on individual subdomains, and combine the resulting local corrections to update the global solution.
These methods are particularly effective at reducing localized or high-frequency error components that are difficult to eliminate using global iterations alone.
A limitation of one-level DD methods is the lack of scalability, since the number of iterations generally increases with the problem size.
To address this issue, two-level DD methods introduce a coarse problem that efficiently reduces low-frequency error components, thereby improving both convergence and scalability.
We refer readers to~\cite{smith_domain_2004,toselli_domain_2005,dolean_introduction_2015,quarteroni_domain_1999} for further details on DD methods.

DD methods can also be classified into overlapping and non-overlapping approaches.
In overlapping methods, neighboring subdomains share common regions through which information is exchanged.
Classical examples include additive and multiplicative Schwarz methods.
In non-overlapping methods, subdomains are disjoint and interact only through their interfaces.
Representative examples include Schur complement methods and substructuring approaches such as finite element tearing and interconnecting (FETI)~\cite{farhat1991method} and balancing domain decomposition by constraints (BDDC)~\cite{dohrmann2003preconditioner}.

Within the \emph{subspace correction} framework, DD methods involve two types of subspaces.
First, the subdomains themselves define local subspaces.
Solving local problems on these subspaces primarily reduces localized or high-frequency error components.
Second, two-level DD methods introduce a coarse subspace, similar to the coarse-grid space in MG methods, to capture global low-frequency errors and improve scalability.

For example, the two-level additive Schwarz preconditioner is given by
\begin{equation}\label{equ:dd_additive_Schwarz}
    M^{-1}_{\mathrm{two\text{-}level}} = \underbrace{ \sum_{i=1}^{N_{\mathrm{dom}}} R_i^T A_i^{-1} R_i }_{\text{local subspace corrections}} + \underbrace{ P_0 A_0^{-1} P_0^T }_{\text{coarse-space correction}},
\end{equation}
where $N_{\mathrm{dom}}$ denotes the number of subdomains and $R_i$ is the restriction operator associated with the $i$-th subdomain.
The local matrix $A_i = R_i A R_i^T$ represents the system restricted to the $i$-th subdomain.
Moreover, $P_0$ denotes the coarse-space basis (or prolongation operator), and
$A_0 = P_0^T A P_0$ is the coefficient matrix of the coarse problem.

\paragraph{\textbf{Deflation and augmentation methods}}

The convergence of (preconditioned) Krylov subspace methods is strongly influenced by the spectral properties of the matrix.
For SPD problems, small eigenvalues may slow down convergence by increasing the condition number of the system.
More generally, slowly converging error components are often associated with invariant subspaces corresponding to unfavorable parts of the spectrum.
Deflation and augmentation techniques~\cite{nicolaides1987deflation, morgan1995restarted, frank2001construction, parks2006recycling} accelerate Krylov subspace methods by incorporating spectral information into the iterative process, leading to deflated and augmented Krylov methods~\cite{gaul2013framework, coulaud2013deflation}.
Both approaches typically exploit eigenvectors associated with slowly converging modes, but they incorporate these subspaces in different ways.

Deflation methods improve convergence by removing undesirable spectral components from the linear system through projection operators.
For simplicity, we consider the deflated conjugate gradient (DCG) method~\cite{frank2001construction}.
Define the projection operator $P = I - AZ(Z^TAZ)^{-1}Z^T$, where the columns of $Z\in \mathbb{R}^{N\times r}$, with $r \ll N$, span the deflation subspace, usually constructed from approximate eigenvectors of $A$ associated with small eigenvalues.
Assuming $A$ is SPD and $Z^TAZ$ is nonsingular, the solution can be decomposed as $\mathbf{x} = (I - P^T)\mathbf{x} + P^T\mathbf{x}$, and these two components are computed separately.
The first component corresponds to the ``coarse-space'' correction and is computed explicitly as $(I - P^T)\mathbf{x} = Z(Z^TAZ)^{-1}Z^T\mathbf{b}$, which only requires solving a small $r\times r$ linear system.
Using the relation $AP^T = PA$, the second component is obtained by solving the deflated system $PA\tilde{\mathbf{x}}=P\mathbf{b}$, followed by $P^T\mathbf{x} = P^T\tilde{\mathbf{x}}$.
By projecting out slowly converging components, the spectrum of the deflated operator is improved, leading to faster convergence.

Augmentation methods take a complementary approach: instead of explicitly modifying the operator, they enlarge the Krylov search space with an additional subspace containing useful spectral information~\cite{parks2006recycling}.
Given an augmentation subspace $\mathcal{W} = \mathrm{span}\{w_1, w_2, \cdots, w_r\}$, the approximate solution is computed in the enlarged space $\mathbf{x}_0 + \mathcal{K}_m(A, \mathbf{r}_0) + \mathcal{W}$, where $\mathcal{K}_m(A, \mathbf{r}_0)$ denotes the $m$-dimensional Krylov subspace generated by $A$ and the initial residual $\mathbf{r}_0$.
The augmentation subspace is typically chosen to approximate invariant subspaces associated with slowly converging error components.
Although conceptually different, deflation and augmentation methods are closely related and can be interpreted within a common framework~\cite{gaul2013framework}.
We refer readers to~\cite{gaul2013framework, soodhalter2020survey} for further details on deflation and augmentation methods.

\subsubsection{Smoothers, subdomain solvers, and coarse-grid solvers}\label{sec:linear-mgdd-smoother}

In MG and DD methods, smoothers, subdomain solvers, and coarse-grid solvers are key components that play a similar role in solving local or coarse-grid subproblems.
These components can be enhanced or replaced by ML models and are discussed separately below.

\paragraph{\textbf{Smoothers (MG)}}
Smoothers in MG methods play a crucial role in damping high-frequency errors on fine grids and are essential for achieving fast convergence.
Relaxation methods are commonly used as smoothers in MG; however, in practice, the optimal choice may vary across grid levels.

\citet{huang2023learning} propose to learn multigrid smoothers using CNNs.
At multigrid level $l$ ($l=0, \ldots, L-1$), the smoother for $A^{(l)}\mathbf{x}^{(l)}=\mathbf{b}^{(l)}$ is written as
\[
    \mathbf{x}^{(l)}_{k+1}
    = \mathbf{x}^{(l)}_k + (M^{(l)})^{-1}
    \bigl(\mathbf{b}^{(l)}-A^{(l)}\mathbf{x}^{(l)}_k\bigr),
    \qquad
    k=0, 1, \ldots, K^{(l)}_{\text{smo}} - 1,
\]
where $K^{(l)}_{\text{smo}}$ denotes the number of smoothing steps.
In their method, a CNN, denoted by $\mathcal{N}^{(l)}$, is used to approximate the action of $(M^{(l)})^{-1}$, which can be interpreted as learning a level-dependent preconditioner; see~\Cref{sec:linear-precond}.
Rather than training smoothers independently at different levels, they propose a recursive training strategy to promote effective smoothing behavior across all grid levels.

\paragraph{\textbf{Subdomain solvers (DD)}}
In DD methods, the global linear system is decomposed into multiple smaller systems defined on subdomains, which are then solved by local subdomain solvers to reduce errors within each subdomain.
Here, we focus on ML-based approaches for solving the linear systems arising on subdomains, rather than on using ML techniques to directly discretize PDEs on subdomains.
We interpret the latter as DD-inspired ML methods and will briefly discuss them in~\Cref{sec:related-inspired_nn-architectures}.
For further details, we also refer to~\cite[Section~5.2]{klawonn2024machine}.

The deep statistical solvers (DSS) method~\cite{donon2020deep} is a GNN-based approach for solving PDEs on unstructured meshes.
For a discretized system $A \mathbf{x} = \mathbf{b}$ of size $N$, DSS predicts an approximate solution $\hat{\mathbf{x}}=\operatorname{DSS}(N,A,\mathbf{b})$, and is trained using a residual loss
\[
    \mathcal{L} = \frac{1}{N}\|\mathbf{b} - A \hat{\mathbf{x}}\|_2^2.
\]
However, the number of message-passing layers required by DSS grows with the system size, which limits its applicability to large-scale problems.
To address this issue, \citet{nastorg2024multilevela} propose DDM-GNN, a two-level additive Schwarz preconditioner (\Cref{equ:dd_additive_Schwarz}) in which the local subdomain solves are replaced by a DSS model.
For the $i$-th subdomain and a global residual $\mathbf r$, DSS computes
\[
    \hat{\mathbf z}_i
    =
    \operatorname{DSS}
    \left(
    N_i, A_i, \frac{R_i\mathbf r}{\|R_i\mathbf r\|}
    \right),
\]
where $N_i$ denotes the system size on the subdomain, $A_i = R_iAR_i^T$, and the local residual is normalized to prevent stagnation.
The resulting DDM-GNN preconditioner is applied to the global residual as
\[
    \mathbf z
    =
    M_{\mathrm{DDM\text{-}GNN}}^{-1}(\mathbf r)
    =
    \sum_{i=1}^{N_{\mathrm{dom}}}
    R_i^T \|R_i\mathbf r\|\,\hat{\mathbf z}_i
    +
    \mathbf r_c,
\]
where $\mathbf r_c:=P_0A_0^{-1}P_0^T\mathbf r$ is the coarse-grid correction computed by a direct solver.
The learned preconditioner is used within PCG.

\paragraph{\textbf{Coarse-grid solvers (MG and DD)}}
Coarse-grid solvers in MG and DD methods are responsible for reducing low-frequency errors and are therefore essential for achieving fast convergence.

\citet{azulay2023multigridaugmented} combine classical iterative methods with a U-Net model to construct preconditioners for solving heterogeneous Helmholtz equations at high wavenumbers within an FGMRES framework.
The U-Net $\mathcal{N}$ approximates an inverse operator by taking Helmholtz parameters $\mathbf{c}$ and the residual $\mathbf{r}$ as input and outputting the error correction $\mathbf{e}$, \ie,
\[
    \mathbf{e} = \mathcal{N}(\mathbf{c}, \mathbf{r}).
\]
Two preconditioners are proposed.
The first mimics a two-grid V-cycle, consisting of Jacobi pre-smoothing, U-Net correction, and Jacobi post-smoothing.
The U-Net acts as a learned coarse-grid correction operator: its encoder-decoder structure, analogous to restriction and prolongation, captures global low-frequency components on a coarser representation and maps them back to the fine grid.
In the second approach, the U-Net first processes the residual to capture the low-frequency component, followed by refinement using a shifted Laplacian multigrid V-cycle, which provides additional smoothing~\cite{umetani2009multigrid, erlangga2006novel}.
Furthermore, the authors propose an encoder-solver framework that separates the U-Net inputs to improve generalization.
The ``encoder'' extracts hierarchical features from the Helmholtz parameter, while the ``solver'' uses these features with the residual to produce updates during the iterative solution.

The encoder-solver framework in~\cite{azulay2023multigridaugmented} is extended in~\cite{lerer2024multigridaugmented} by replacing the standard convolutional layer at the coarsest level of the solver U-Net with an implicit layer.
This implicit layer is designed to mimic a direct coarse-grid solver in multigrid methods, where a convolution kernel is inverted via the FFT rather than applied directly, enabling more effective global information propagation.
In~\cite{cui2025neural}, the Chebyshev semi-iteration method is employed as the coarse-grid solver for Helmholtz equations with high wavenumber and heterogeneous media.
At different grid levels, the parameters of the Chebyshev semi-iteration method are predicted by a CNN, as described in the context of parameter learning in~\Cref{sec:linear-core_components-param}.
\citet{li2025neural} propose the neural algebraic multigrid (NAMG) operator, which integrates attention mechanisms into the classical AMG framework.
In particular, the coarse-grid correction is approximated using a self-attention model.
In addition, the singularity-encoded Green's function learning method is used in~\cite{sun2026learning} as the coarse-grid solver within two-level overlapping Schwarz preconditioners; see~\Cref{sec:linear-precond-hybrid}.

\paragraph{\textbf{Discussion}}

Despite these developments, relatively little work has considered complex problem settings.
More broadly, many existing ML-based preconditioning techniques could in principle be employed as inexact local or coarse-grid solvers within MG and DD frameworks, suggesting substantial opportunities for future research.

In PDE settings, smoothers are typically used to reduce high-frequency error components.
Since neural networks generally exhibit a spectral bias toward low-frequency components, designing ML-based smoothers that effectively damp high-frequency errors remains challenging.
Nevertheless, several approaches have been proposed, including local solver-based methods~\cite{moseley2023finite, dolean2024multilevel} and methods based on random Fourier features~\cite{wu2026iterative}; see also~\cite{xu2025understanding}.

Coarse-grid correction, in contrast, aims to eliminate low-frequency errors, which aligns naturally with the spectral bias of neural networks.
This observation suggests that a wide range of ML approaches, for example, those in~\cite{meng2025when, cuomo2022scientific}, could be leveraged to construct effective coarse-grid solvers.

\subsubsection{Coarsening in algebraic multigrid}\label{sec:linear-mgdd-coarsening}
In this subsection, we focus on coarsening strategies in AMG, while transfer operators are discussed in the next subsection.
The coarse grid is typically selected from the fine grid and is constructed to capture low-frequency residual components at reduced computational cost.
In GMG, the coarse grid can be obtained directly from the geometric structure of the problem.
In contrast, in AMG, the coarse grid is constructed algebraically from the coefficient matrix $A$.
At each level, AMG selects an appropriate coarse grid and corresponding transfer operators.
The choice of the coarse/fine splitting is crucial, as it balances the convergence speed and computational cost.
However, identifying an optimal coarse grid is intractable; in fact, the problem is NP-complete~\cite{maclachlan2007greedy}.
To address this challenge, ML-based methods have been proposed to construct coarse grids in AMG, either by \emph{directly selecting the coarse grid} or by \emph{learning the strong threshold}, a parameter in classical coarsening strategies.

\paragraph{\textbf{Directly selecting the coarse grid}}

Coarse-grid generation is commonly formulated as a graph coarsening problem~\cite{chen2022graph}.
In this view, the coarse grid corresponds to a subgraph of the fine grid, where coarse nodes are selected based on the strength of connections.
\citet{taghibakhshi2021optimizationbased} employ a GNN-based reinforcement learning agent to generate coarse graphs.
The agent is trained on small graphs and then applied to larger, unstructured graphs.
The proposed approach outperforms the greedy coarsening strategy~\cite{maclachlan2007greedya} in terms of the F-fraction, defined as the ratio of fine nodes to the total number of nodes.

In~\cite{gottschalk2021coarsening}, the authors propose a kriging-based approach to construct AMG components, including both coarse-fine splitting and prolongation.
They assume that the smooth error in AMG can be modeled algebraically as a Gaussian process.
Given data on the coarse grid (treated as observations in the kriging framework) and a prescribed covariance structure, the kriging model is used to infer the conditional distribution of values at fine-grid nodes, including their means and variances.
Building on the fact that a Gaussian process exhibits low conditional variance near observed points and higher variance elsewhere, the coarse grid is constructed in a greedy manner: at each step, the point with the largest conditional variance is selected and added to the coarse grid.

\paragraph{\textbf{Learning the strong threshold}}

AMG coarse-fine splitting relies on the notion of strength of connection, which is typically defined via a threshold parameter (often referred to as the strong threshold)~\cite{ruge1987algebraic}.
Classical approaches include the Ruge-St\"uben method~\cite{ruge1987algebraic}, Cleary-Luby-Jones-Plassmann (CLJP)~\cite{cleary1998coarse}, parallel modified independent set (PMIS)~\cite{de2006reducing}, and hybrid modified independent set (HMIS)~\cite{de2006reducing}.
Given a matrix $A \in \mathbb{R}^{N \times N}$ and a threshold parameter $0 < \theta \leq 1$, the variable $i$ is said to strongly depend on variable $j$ ($j \neq i$) if
\[
    |A_{ij}|
    \geq
    \theta \max_{l \neq i}\left\{|A_{il}|\right\},
    \qquad i,j=1,\ldots,N.
\]
The choice of this threshold is problem-dependent and usually requires careful tuning, often based on heuristic strategies.
Recently, several works have exploited ML methods to automatically find appropriate threshold values.

\citet{antonietti2023accelerating} propose AMG-ANN, which uses (artificial) neural networks to predict the optimal strong threshold parameter $\theta$ for AMG.
In this framework, the coefficient matrix is interpreted as a grayscale image, and CNNs are used to extract features $\mathbf{y}_{\text{feat}}$ (a vector), \ie,
\[
    \mathbf{y}_{\text{feat}} = \CNN(A).
\]
These features, together with the mesh size $h$ and a candidate threshold $\theta$, are fed into an FNN that outputs the corresponding approximate convergence factor,
\[
    \rho = \FNN(\mathbf{y}_{\text{feat}},-\log_2 h,\theta),
\]
The approximate convergence factor is defined as
\[
    \rho=
    \left(
    \frac{\|\mathbf r_{N_{\mathrm{it}}}\|_2}
    {\|\mathbf r_0\|_2}
    \right)^{1/N_{\mathrm{it}}},
\]
where $N_{\mathrm{it}}$ is the number of iterations required to reach the prescribed residual tolerance.
The network is trained in a supervised manner.
Since the network output is differentiable with respect to $\theta$, the optimal threshold can be obtained by minimizing the predicted convergence factor via gradient descent.
The method is demonstrated on 2D elliptic and stationary Stokes problems.

The AMG-ANN method~\cite{antonietti2023accelerating} is extended to 3D elliptic PDEs in~\cite{caldana2024deep}, where the FNN additionally takes the polynomial degree of the finite element discretization as input.
Instead of predicting convergence factors, the model estimates computational time (including setup and solve time) and a corresponding uncertainty.
Furthermore, in~\cite{antonietti2025deep}, the authors extend AMG-ANN to simultaneously select both smoothers and the strong threshold for linear systems arising from polytopal discretizations.

\citet{zou2024autoamgsettled} propose an auto-tuning framework, called AutoAMG($\theta$), for learning strong threshold parameters.
Different from prior work that uses CNNs for feature extraction~\cite{antonietti2023accelerating, caldana2024deep}, this method introduces a graph convolutional isomorphism network (GCIN), which combines GCNs~\cite{kipf2017semisupervised} and GINs~\cite{xu2019how}.
This hybrid architecture improves feature extraction compared to using either model alone.
The extracted features are then passed to an FNN that predicts the strong threshold $\hat \theta$.

Building on~\cite{jiang2022general}, where Gaussian process regression (GPR) is used to predict the splitting parameter for alternating-direction implicit (ADI) methods, \citet{luo2026kernellearning} employ GPR to learn the strong threshold for AMG.
Specifically, they approximate a mapping from the dimension of the linear system to the optimal strong threshold.
This is done using a kernel learning strategy, in which the kernel is selected from linear combinations of multiplicative combinations of basic kernels, such as radial basis and exponential kernels.
The method is evaluated using several metrics, including mean squared error, prediction interval coverage, and the Bayesian information criterion (BIC).

\paragraph{\textbf{Discussion}}

\emph{Learning the strong threshold} can be regarded as an instance of \emph{parameter learning}, which was discussed in~\Cref{sec:linear-core_components-param}.
The adaptive selection of strong thresholds, for example based on local matrix information, remains an open question.
Beyond the strong threshold, alternative measures for defining connection strength may also be explored, such as Gaussian process regression~\cite{gottschalk2021coarsening} and reinforcement learning.

\subsubsection{Coarse spaces and transfer operators}\label{sec:linear-mgdd-transfer_operator}

Global corrections are essential for scalable and robust subspace correction methods. For efficiency, they are usually built from low-dimensional subspaces, often associated with a coarse grid, that capture the components most important for convergence.
Constructing the corresponding transfer operators between the original space and these subspaces can therefore be costly, which motivates the use of ML models to learn effective global correction components.
The same idea appears under different terminology across classical frameworks: multigrid often refers to transfer operators, domain decomposition to coarse spaces, and deflation (or augmentation) methods to deflation (or augmentation) subspaces.
In many settings, the maps to and from the correction subspace are determined by a single operator, so learning transfer operators often amounts to learning this map.
Despite these differences in terminology, the resulting learning tasks are conceptually closely related, and many approaches are applicable across different settings.
One example is coarse spaces based on multiscale methods, which have been used in various types of subspace correction methods.
In the following, we review these approaches, grouped by the type of subspace correction method on which they primarily focus.

\paragraph{\textbf{Multigrid}}

In~\cite{katrutsa2017deep}, the two-grid GMG method is interpreted as a DNN, where each multigrid operation corresponds to a network layer.
Within this framework, the restriction and prolongation operators are treated as trainable parameters associated with these layers, together with the Jacobi damping parameter in the smoother.
Consequently, the iteration matrix can be viewed as a function of these parameters, and the learning problem is formulated as minimizing its spectral radius with respect to them.
Since direct computation of the spectral radius is expensive and non-differentiable, the authors instead optimize a differentiable approximation based on Gelfand's formula~\cite{kozyakin2009accuracy}, where the Frobenius norm of matrix powers is used to approximate the spectral radius.
The later published version of this work~\cite{katrutsa2020blackbox} considers a more general formulation in which the restriction and prolongation operators, as well as the Jacobi damping parameter, are treated directly as trainable variables.
Instead of relying on an explicit neural network interpretation, the multigrid method is implemented within a differentiable programming framework, where gradients are computed via automatic differentiation.

\citet{greenfeld2019learning} employ an FNN to learn prolongation operators within a GMG method for 2D diffusion PDEs with spatially varying coefficients on structured grids.
The method constructs the prolongation matrix using local information, where the FNN maps local discretization stencils to the nonzero entries of the prolongation operator under a predefined sparsity pattern.
The loss function is defined based on the spectral radius of the two-grid error propagation matrix, which characterizes the convergence rate.
In practice, due to the high computational cost, this loss is replaced by the Frobenius norm.
This also avoids backpropagation through eigendecompositions, which has been reported to be prone to numerical instability in certain settings~\cite{wang2019backpropagationfriendly, li2023learning}.
In~\cite{luz2020learning}, this approach is extended to PDEs on unstructured grids through the use of GNNs.
The GNN takes a graph representation of the coefficient matrix as input and predicts a prolongation operator under a prescribed sparsity pattern.
This sparsity pattern is determined using classical AMG coarsening methods such as Ruge-St\"uben~\cite{ruge1987algebraic} or CLJP~\cite{cleary1998coarse}.
Building on this idea, in~\cite{moore2021learning} and~\cite[Section~5]{moore2022machine}, the authors propose employing a GNN to learn a correction to a classical prolongation operator, rather than constructing it from scratch.
This correction formulation is inspired by residual neural networks (ResNets)~\cite{he2016deep}.
Similar correction-based strategies have also been explored for learning incomplete Cholesky factors in~\cite{trifonov2025learning, li2024deep}; see~\Cref{sec:linear-precond-hybrid} for further discussion.

In bootstrap AMG (BAMG)~\cite{brandt2011bootstrap}, local weighted least squares (LWLS) are used to construct prolongation operators, which are then refined through a bootstrap procedure.
However, BAMG employs spatially shared weights in the LWLS formulation, and relies on an expensive bootstrap process involving eigenvalue computations.
To address these limitations, \citet{wang2023learningbased} propose a learning-based LWLS (L-LWLS) method for constructing prolongation operators.
Specifically, an MLP is used to learn spatially varying weights in the LWLS model, and the resulting LWLS problem is solved via gradient descent.
In addition, a GNN is used to extract local features from the coefficient matrix and to learn a correction to the LWLS solution, effectively replacing the costly bootstrap procedure.
\citet{holguin2024accelerating} reinterpret the prolongation operator in GMG, which maps coarse-grid solutions to fine-grid representations, as performing a super-resolution task.
In image analysis, super-resolution models aim to reconstruct high-resolution images from low-resolution inputs.
Based on this analogy, they employ a super-resolution generative adversarial network (SR-GAN)~\cite{ledig2017photo} to approximate the prolongation operator.
\citet{tomasi2022construction} consider linear systems arising from finite element discretizations and propose using DNNs to construct $L^2$-projection-based transfer operators~\cite{krause2016parallel}.
The prolongation operator is given by
\[
    P = M_h^{-1}B_h,
\]
where $M_h$ is the fine-grid mass matrix, and $B_h$ is the coupling matrix between the fine and coarse meshes.
Instead of directly learning $P$, the networks predict local entries of $B_h$ from the corresponding local information in $M_h$, after which $P$ is assembled and used in the multigrid hierarchy.

\paragraph{\textbf{Domain decomposition}}

Coarse spaces are important for ensuring the robustness and scalability of DD methods.
For challenging problems, such as elliptic PDEs with highly heterogeneous coefficients, standard coarse spaces may be insufficient.
Adaptive coarse spaces address this issue by adding components that are difficult to eliminate by local solvers.
When these components are identified through local eigenvalue problems, the resulting spaces are often referred to as spectral coarse spaces~\cite{galvis2010domain, galvis2010domainreduced, mandel2007adaptive, klawonn2015feti, spillane2014abstract}.
Such methods are generally more robust, but their setup can be expensive, since it may require solving local generalized eigenvalue problems, for instance on subdomain interfaces.
ML techniques have therefore been explored to reduce this cost.
We refer to~\cite{heinlein2021combining, klawonn2024machine} for comprehensive overviews of these topics and of combining ML with DD methods more generally.

In adaptive finite element tearing and interconnecting-dual primal (FETI-DP)~\cite{farhat2001feti, klawonn2012deflation}, many interface eigenvalue problems do not contribute to the coarse space, but identifying them a priori is challenging.
In~\cite{heinlein2019machine}, this selection is formulated as a classification problem, leading to ML-FETI-DP, which combines ML with adaptive FETI-DP.
FNN classifiers use samples of the coefficient function from adjacent subdomain components to predict whether an interface eigenvalue problem is required; because this input is discretization-independent, trained models can generalize across meshes.
Further developments consider random training datasets~\cite{heinlein2020machine}, local interface-neighborhood inputs~\cite{heinlein2021machine}, 3D elliptic problems~\cite{heinlein2021combininga}, and critical-interface identification for adaptive generalized Dryja-Smith-Widlund (AGDSW) coarse spaces~\cite{heinlein2018adaptive, heinlein2021combining, heinlein2022predicting}.

Following up on this classification approach, \citet{klawonn2024learninga} then treat coarse space construction as a regression task, using FNNs to learn interface constraints directly and thereby avoid local eigenvalue problems.
Developed first for regular domain decompositions with straight edges, this approach is extended to irregular decompositions in~\cite{klawonn2024learning} and to nonlinear FETI-DP methods in~\cite{klawonn2025learning}; see~\Cref{sec:nonlinear}.
Furthermore, \citet{ciaramella2022spectral} propose two alternative approaches for constructing coarse spaces based on PCA and DNNs, rather than solving generalized eigenvalue problems explicitly.
The coarse space is designed to capture dominant eigenvectors of the smoother on the fine grid, which correspond to low-frequency error components.
The PCA-based approach approximates these eigenvectors using smoothed random vectors, while the DNN-based approach interprets a two-level DD method as a DNN and treats the coarse basis functions as trainable parameters, similarly to~\cite{katrutsa2017deep}.

In~\cite{chung2021learning}, FNNs are integrated with an adaptive BDDC method~\cite{kim2015bddc} for stochastic elliptic problems.
An FNN maps the stochastic coefficient, represented by a Karhunen-Lo\`eve expansion, to eigenvectors of local generalized eigenvalue problems, which are then used to construct the coarse space.
\citet{kopanicakova2025deeponet} employ DeepONet to construct preconditioners within a subspace correction framework using direct preconditioning (DP) and trunk basis (TB) approaches.
In the DP approach, the DeepONet takes the PDE parameters $\mathbf c$ and the residual $\mathbf r_k$ as inputs and approximates the corresponding coarse-space correction,
\[
    \mathcal{N}(\mathbf{c}, \mathbf{r}_k)
    \approx
    P_0 A_0(\mathbf{c})^{-1} R_0 \mathbf{r}_k,
\]
where $A_0(\mathbf{c})=R_0A(\mathbf{c})P_0$ and $R_0=P_0^T$.
In the TB approach, basis functions extracted from the DeepONet trunk network are used to construct the prolongation operator $P_0$, which is orthogonalized via QR factorization to ensure full rank and improve conditioning.
The DP-based preconditioner is nonlinear and is therefore used within FGMRES, whereas the TB-based preconditioner is still linear and can thus be used with standard Krylov methods.
The resulting DP- and TB-based preconditioners can also be combined with classical stationary iterations or Krylov subspace methods, yielding alternating HIMs; see~\Cref{sec:linear-precond-alternating}.
Compared with the classical Nicolaides coarse space~\cite{nicolaides1987deflation}, the TB-based coarse space is more expensive but remains effective for problems with jumping coefficients or indefiniteness.
Its construction is cheaper than more sophisticated coarse spaces based on the Dirichlet-to-Neumann map~\cite{dolean2012analysis} or GenEO (Generalized Eigenproblems in the Overlap)~\cite{spillane2014abstract}.
Finally, in a GNN-based approach, Robin-type interface conditions and coarse spaces are derived in~\cite{taghibakhshi2023mggnn}; see~\Cref{sec:linear-mgdd-parameters} for more details.

\paragraph{\textbf{Multiscale methods}}
Multiscale methods, such as multiscale finite volume methods (MsFVMs)~\cite{hajibeygi2008iterative} and multiscale finite element methods (MsFEMs)~\cite{efendiev2013generalized}, can be used to construct two-level preconditioners.
Specifically, the multiscale basis functions can be used as coarse spaces in DD and as transfer operators in MG~\cite{heinlein2025algebraic, alves2024computational, fu2024efficient}.
However, constructing such multiscale basis functions typically requires solving local (spectral) problems, which can be computationally expensive.
The following works employ ML models to learn prolongation operators and reduce this cost.

The pore-level multiscale method (PLMM)~\cite{mehmani2021multiscale} is designed to solve elliptic equations arising from linear-elastic mechanics at the pore scale, where constructing the prolongation matrix is costly.
In~\cite{li2024machine}, the authors propose a CNN-based framework, combining U-Net and ResNet, to accelerate the construction of the prolongation matrix within two-level PLMM preconditioners.
The network takes as input cropped images of decomposed pore-scale subdomains, encoding both local geometry and boundary condition information, and predicts local basis functions that are subsequently assembled to form the prolongation matrix.

\citet{liu2025learning} employ a U-Net to reduce the construction cost of two-level multiscale preconditioners used within PCG for solving Darcy flow problems with random permeability coefficients.
For each realization of the permeability field, constructing the multiscale prolongation operator typically requires solving multiple local spectral problems, which can be computationally expensive.
To address this issue, a U-Net is trained to learn the mapping from local permeability fields to the corresponding generalized multiscale basis functions, which are then used to efficiently construct the prolongation operator.

\paragraph{\textbf{Deflation and augmentation methods}}

Deflation and augmentation methods are commonly used to accelerate the convergence of Krylov subspace methods.
They both introduce a subspace (called a deflation or augmentation subspace) to remove components of the error or residual associated with ``difficult'' eigenvalues, typically small or near-zero ones, which are often linked to slowly converging (low-frequency) error components.
However, in practice, constructing such subspaces may be computationally costly, which motivates the use of ML methods to learn deflation or augmentation subspaces.

When solving a sequence of systems of linear equations, the Krylov subspace recycling technique can be employed to accelerate the convergence of Krylov subspace methods by reusing the invariant subspace from previous systems.
In~\cite{carlberg2016krylovsubspace}, the authors employ a POD method to construct a subspace from the previous solutions to serve as the recycling subspace.
\citet{trivedi2019datadriven} accelerate GMRES for electromagnetic simulations by using PCA and CNN models to construct a low-dimensional subspace that captures the dominant components of the solution, which is then used to augment the Krylov subspace.

Inspired by Krylov subspace recycling methods~\cite{parks2006recycling}, \citet{luo2024neural} propose a neural Krylov iteration (NeurKItt) method, in which an FNO predicts an approximate invariant subspace associated with the smallest eigenvalues of the corresponding linear system from the PDE parameters.
    A QR decomposition is applied to the network output to obtain an orthonormal basis.
    The learned subspace is then incorporated into a deflated Krylov method to accelerate convergence.

In~\cite{arisaka2025learning}, the authors employ the gradient-based meta-solving (GBMS) framework~\cite{arisaka2023principled} to train a GNN model for learning the deflation subspace.
The GNN takes the coefficient matrix and the right-hand side vector as input and predicts the deflation vectors.
\citet{kopanicakova2025leveraging} consider parametric linear systems and propose two approaches for learning the deflation operator using DeepONet.
The first is the TB approach~\cite{kopanicakova2025deeponet}, in which the outputs of the trunk network are used to construct the deflation operator.
Inspired by classical recycling techniques~\cite{parks2006recycling}, the second approach constructs the deflation operator using solution approximations generated by DeepONet.

\paragraph{\textbf{Discussion}}

For multigrid methods, the coarsening strategy and transfer operators should ideally be designed in a consistent manner so that the coarse-grid correction can effectively represent and eliminate low-frequency error components.
Most existing works treat these two components separately: one is learned by ML models, while the other is constructed through classical approaches.
For example, \emph{transfer operator learning} methods typically assume a prescribed sparsity pattern for the transfer operators, which implicitly fixes the coarse grid.
A promising research direction is therefore to develop adaptive approaches that jointly learn the coarsening and the associated transfer operators.

\subsubsection{Parameters in multigrid and domain decomposition}\label{sec:linear-mgdd-parameters}
MG and DD methods are sophisticated iterative techniques involving numerous parameters, which are often chosen based on experience or prior knowledge.
ML can be employed to automatically select or learn these parameters, thereby improving overall performance.

\paragraph{\textbf{MG parameters}}

In the previous subsections, we discussed ML-based approaches for learning individual MG components, including smoothers, coarse-grid solvers, and transfer operators.
Here, we turn to other aspects of MG design, such as selecting effective combinations of these components.
This task can be formulated as a discrete optimization problem, for which evolutionary algorithms, such as genetic algorithms (GAs) and genetic programming (GP), are particularly suitable.
We refer to~\cite[Section 8]{schmitt2024automating} for a comprehensive review of the automated design of MG methods.

As an early example, \citet{oosterlee2003genetic} introduce a GA-based approach for automatically selecting multigrid components, including pre- and post-smoothers, relaxation parameters, and transfer operators, from a predefined search space.
Fourier analysis is used to evaluate the performance of candidate configurations.
In~\cite{thekale2010optimizing}, a branch-and-bound method is proposed to determine the optimal number of multigrid cycles within a full MG scheme.
The method efficiently explores and prunes the search space, with the objective of minimizing the computational cost required to achieve a prescribed accuracy.
During the evolutionary search process, many candidate parameter configurations are generated.
Since directly evaluating the computational cost of each candidate is expensive, \citet{petrushov2023automated} propose a hybrid evolutionary strategy that combines an evolutionary algorithm with a neural network model for optimizing AMG configurations.
Specifically, the neural network is trained to predict the solution time and thus serves as a surrogate model.
During the mutation process, this surrogate is used to discard candidates with high predicted computational cost.

Instead of restricting to standard V- or W-cycles, \citet{schmitt2020constructing, schmitt2021evostencils} aim to automatically discover flexible GMG cycle structures.
By representing a GMG solver as a tree of mathematical expressions, they formulate the design process as a program optimization problem.
A context-free grammar is introduced, and grammar-guided genetic programming (G3P) is used to explore the space of candidate solvers.
The quality of each candidate is evaluated in terms of both convergence rate and computational performance.
This framework has been applied to Helmholtz problems in~\cite{schmitt2022evolving} and extended to AMG in~\cite{parthasarathy2024evolving}.
Moreover, it has also been generalized to FAS, the nonlinear multigrid method; see~\Cref{sec:nonlinear} for further details.

\paragraph{\textbf{DD parameters}}

In overlapping additive Schwarz methods, the overlap size plays a critical role in overall performance.
A larger overlap improves information exchange between subdomains and typically reduces the number of iterations required for convergence.
However, it also increases the size of the local subproblems, leading to higher computational cost.
Therefore, there exists a trade-off between convergence rate and computational efficiency.
In~\cite{burrows2013learning}, the authors employ a regression model to determine the optimal overlap width that minimizes the total number of floating-point operations (FPOs).
They consider a 2D Poisson problem with discontinuous coefficients and construct a model that takes local neighborhood features (\eg, maximum and minimum coefficient values) as input and predicts the corresponding number of FPOs.
The trained model is then used to estimate the FPO count for previously unseen local configurations.

\citet{taghibakhshi2022learning} exploit GCNs~\cite{kipf2017semisupervised} to learn optimal subdomain interface conditions for a one-level optimized Schwarz method applied to problems defined on unstructured grids.
This work is further extended in~\cite{taghibakhshi2023mggnn}, where the authors propose a multigrid graph neural network (MG-GNN) framework for learning optimized parameters in two-level DD methods in an unsupervised manner.
In particular, they focus on learning Robin-type interface conditions and coarse spaces.
\citet{knoke2023domain, knoke2024neural} study the solution of time-harmonic Maxwell's equations using optimized Schwarz methods~\cite{el2015quasi}, whose convergence behavior is highly sensitive to the choice of interface operators.
To address this issue, they employ FNNs to approximate optimal interface operators between subdomains.

\subsection{Automated solver selection}\label{sec:linear-solver_selection}

For a given linear system, various combinations of iterative methods and preconditioners generally exhibit different convergence behaviors.
When solving nonlinear systems, one may also need to choose an appropriate linear solver at each nonlinear iteration.
In practice, selecting a solver is typically based on experience or knowledge of the problem.
Moreover, when different iterative schemes, preconditioners, and their associated parameters are taken into account, the number of possible combinations becomes extremely large.
\emph{Automated solver selection} uses ML methods to identify an optimal iterative method, a preconditioner, or their combination from a given search space based on properties of the problem of interest.
Note that \emph{automated solver selection} is one of the earliest topics within the field of ML-enhanced iterative methods, for example, \cite{xu2005study, bhowmick2006application}; here, we only review some representative works and do not aim to be exhaustive.

We first introduce, in~\Cref{sec:linear-solver_selection-prelims}, a general workflow for approaches based on \emph{numerical features} and \emph{structural features}.
We then review representative methods following this paradigm in~\Cref{sec:linear-solver_selection-numerical_structural_features}.
Finally, in~\Cref{sec:linear-solver_selection-other}, we discuss two alternative approaches based on online learning.

\subsubsection{Preliminaries}\label{sec:linear-solver_selection-prelims}

\begin{figure}[!htbp]
    \centering
    \includegraphics[width=1.0\textwidth]{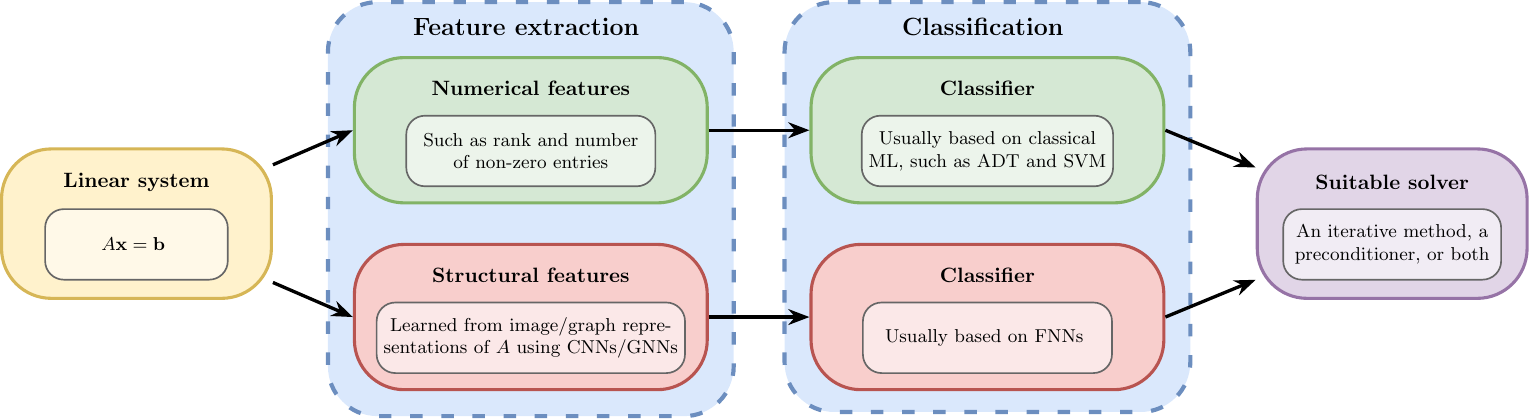}
    \caption{
        Schematic of \emph{automated solver selection} based on numerical and structural features.
        The general workflow consists of two stages: feature extraction and classification.
        The upper and lower paths can be used independently or in combination.
        In the upper path, selected numerical features, such as the rank and the number of nonzero entries, are fed into a classifier, typically implemented using classical ML methods such as ADT and SVM.
        In the lower path, CNNs or GNNs are used to extract structural features from image- or graph-based representations of the input matrix, and an FNN performs the classification.
        When the two paths are combined, the extracted numerical and structural features are fused into a unified representation and fed into a classifier for solver selection.
    }
    \label{fig:automated-selection}
\end{figure}

\Cref{fig:automated-selection} provides a schematic illustration of \emph{automated solver selection} based on numerical and structural features.
Given a linear system $A\mathbf{x} = \mathbf{b}$, these approaches typically involve two steps~\cite{bhowmick2006application}:
\begin{enumerate}
    \item \textbf{Feature extraction:}
          Various features are extracted from $A$ to characterize the target linear system and can generally be classified into two types~\cite{tang2022graph,xiong2025mmautosolver}:
          \begin{itemize}
              \item \textbf{Numerical features:}
                    Numerical features refer to global matrix properties, such as the rank, condition number, and proportion of nonzero entries of the coefficient matrix $A$.
              \item \textbf{Structural features:}
                    By representing the coefficient matrix $A$ as an image or a graph, CNNs or GNNs can be employed to extract \emph{structural features}, which typically capture local properties and sparsity patterns of the matrix; see also~\Cref{fig:prelims-struct_unstruct_data}(b) and (d).
          \end{itemize}
    \item \textbf{Classification:}
          An ML-based classifier takes the extracted features as input and outputs a suitable solver configuration, including an iterative method, a preconditioner, or a combination of both.
          Classical classification algorithms, such as alternating decision trees (ADTs) and support vector machines (SVMs), can be used.
          For structural features extracted by CNNs or GNNs, FNNs are often used as classifiers, as in~\cite{yamada2018preconditioner,tang2022graph,xiong2025mmautosolver}.
\end{enumerate}

When both numerical and structural features are used, they are typically fused before being passed to the classifier, as in~\cite{tang2022graph,xiong2025mmautosolver}.
Most existing studies focus on the design or improvement of feature extraction and classifier construction, while solver performance is commonly evaluated using metrics such as solution time, memory cost, iteration count, and solution accuracy.

\subsubsection{Based on numerical and structural features}\label{sec:linear-solver_selection-numerical_structural_features}

According to the types of features considered, existing approaches can be categorized into three groups: (1) \textbf{numerical features only}, (2) \textbf{structural features only}, and (3) \textbf{both numerical and structural features}.
These methods are summarized in~\Cref{tab:linear-solver_selection-comparison} and described below.

\begin{table}[!ht]
    \renewcommand{\arraystretch}{1.3}
    \centering
    \caption{
        Comparison of existing approaches for \emph{automated solver selection} based on numerical and structural features~(\Cref{sec:linear-solver_selection-numerical_structural_features}).
    }
    \label{tab:linear-solver_selection-comparison}
    \begin{tabular}{
            lllllp{6cm}
        }
        \toprule
        \multirow{2}{*}{\textbf{ref.}}
         & \multicolumn{2}{c}{\textbf{features}}
         & \multicolumn{2}{c}{\textbf{classifier}}
         & \multirow{2}{*}{\makecell{\textbf{solver search space}\\(\# of candidates used for different problems)}}
        \\
         & \textbf{numerical} (\#)                                                                                  & \textbf{structural}
         & \textbf{numerical}                                                                                       & \textbf{structural}
         &
        \\
        \midrule
        \cite{bhowmick2006application}
         & not specified
         & \xmark
         & ADT, SVM
         & \xmark
         & Krylov methods + preconditioners (48, 242)
        \\
        \cite{bhowmick2009lowcost}
         & $57 \xrightarrow{\text{reduced}} 27$
         & \xmark
         & ADT, SVM
         & \xmark
         & Krylov methods + preconditioners (80)
        \\
        \cite{sun2025new}
         & $66 \xrightarrow{\text{reduced}} 8$
         & \xmark
         & RF, XGBoost, GBDT
         & \xmark
         & Krylov solver + restart parameter (6)
        \\
        \cite{holloway2007neural}
         & 35
         & \xmark
         & FNNs
         & \xmark
         & ILUTP\_Mem parameters (72)
        \\
        \cite{kuefler2008using}
         & 32
         & \xmark
         & reinforcement learning
         & \xmark
         & preconditioned solvers (576)
        \\
        \cite{peairs2011using}
         & 1 (residual norm)
         & \xmark
         & reinforcement learning
         & \xmark
         & GMRES restart parameters (6)
        \\
        \cite{yeom2016datadriven}
         & 35
         & \xmark
         & kNN + PVS
         & \xmark
         & 9 iterative solvers + 13 preconditioners (117)
        \\
        \cite{yamada2018preconditioner}
         & \xmark
         & CNN
         & \xmark
         & FNN
         & preconditioners (6)
        \\
        \cite{tang2022graph}
         & 10
         & GNN
         & FNN
         & FNN
         & 5 Krylov solvers + 9 preconditioners (33 compatible)
        \\
        \cite{xiong2025mmautosolver}
         & 17
         & CNN
         & FNN
         & FNN
         & Krylov solver + restart parameter (19)
        \\
        \bottomrule
    \end{tabular}
\end{table}

\paragraph{\textbf{Numerical features only}}

In~\cite{bhowmick2006application}, the authors propose a framework that employs ML techniques to learn mappings from matrix features to an efficient solver.
Many features are considered, including the rank, sparsity pattern, and spectral properties of the coefficient matrix.
The method adopts ADT~\cite{freund1999alternating} based on AdaBoost~\cite{freund1997decisiontheoretic} as the classifier.
AdaBoost is an ensemble learning algorithm that helps mitigate overfitting, particularly in scenarios where the training dataset is limited or dominated by a specific class of problems.
In such cases, a model may perform well on the training data but generalize poorly to unseen instances.

However, some features, such as the condition number of the coefficient matrix, are expensive to compute.
To mitigate this issue, \citet{bhowmick2009lowcost} propose a low-cost classifier construction method that removes low-information and redundant features, and orders the remaining ones by computational cost so that inexpensive features are evaluated first.
Features are then incrementally added until the desired prediction accuracy is achieved, thereby reducing the overall feature extraction cost.
Here, low-information features refer to those that contribute little to predictive performance, typically characterized by having variance below a predefined threshold.
Redundant features, on the other hand, are those with similar statistical properties (\eg, identical mean and variance), for which only the feature with the lowest computational cost is retained while the others are eliminated.
Building on this work, \citet{sun2025new} propose an extended method that incorporates additional feature selection techniques.
Specifically, genetic algorithms and information gain~\cite{larose2014discovering} are used to further reduce the feature set.
Random forest (RF)~\cite{breiman2001random}, XGBoost~\cite{chen2016xgboost}, and gradient boosting decision tree (GBDT)~\cite{friedman2001greedy} are tested as classifiers.

The authors of~\cite{xu2005study, xu2005data} propose employing an SVM model to automate preconditioner selection by predicting the solvability of a given linear system, \ie, whether a linear system can be successfully solved under a given solver configuration.
Their framework considers multiple preconditioners, including ILU(0), ILU(k), and ILUT~\cite{saad1994ilut}, together with their associated parameters.
The SVM model takes matrix features and solver parameters as input and predicts whether preconditioned GMRES converges within a prescribed iteration limit and residual tolerance.
This idea is later extended to parameter selection for SPAI-type preconditioners in~\cite{kim2010data}.
Moreover, different classification models were further investigated in~\cite{xu2016using}.

\citet{holloway2007neural} use an NN to predict the convergence behavior of preconditioned GMRES(50) for sparse linear systems.
The study focuses on the ILU threshold pivoting memory-efficient (ILUTP\_Mem) preconditioner~\cite{chen2004ilutp_mem}, parameterized by three tuning variables.
The network takes matrix features and preconditioner parameters as input and predicts whether the resulting preconditioned system converges.
The authors further show that prediction accuracy improves when the training and test matrices are more similar.

Supervised learning methods for solver selection typically require a large amount of training data.
Alternatively, \citet{kuefler2008using} propose a reinforcement learning framework for constructing preconditioned iterative solvers for sparse linear systems.
Their approach formulates solver selection as a sequential decision-making process.
The state is characterized by matrix features and the current stage of the solution process, the actions include operations such as matrix scaling, reordering, and preconditioning, and the reward is based on the negative runtime.
Analogous to~\cite{holloway2007neural}, the authors consider GMRES(50) preconditioned by ILUTP\_Mem, while exploring different combinations of ILUTP\_Mem parameters.
In~\cite{peairs2011using}, this framework is further extended to select the restart parameter in GMRES($m$).
In~\cite[Section~4]{moore2022machine}, the author employs FNNs to predict suitable parameters for ILUT preconditioners based on the coefficient matrix $A$.
Methods in~\cite{xu2005study, xu2005data, holloway2007neural, kuefler2008using, peairs2011using} and~\cite[Section~4]{moore2022machine} can also be regarded as learning parameters in ILU-type preconditioners and GMRES($m$); see~\Cref{sec:linear-core_components-param}.

\citet{yeom2016datadriven} propose a data-driven method to predict the fastest combination of preconditioner and solver for sparse linear systems.
They introduce a performance vector space (PVS), inspired by the \emph{word2vec} algorithm in natural language processing~\cite{mikolov2013distributed}.
Each matrix is treated as an ``input word'', while well-performing solver configurations are treated as its ``context words''.
A neural network with one hidden layer is used to learn embeddings such that matrices with similar solver performance are close in the PVS.
For inference, a new matrix is approximated as a sparse linear combination of training matrices in the feature space.
Assuming that matrices with similar features admit similar solver performance, the same coefficients are used to compute its embedding in the PVS.
Finally, a $k$-nearest neighbors (kNN) classifier is applied in the PVS to select the solver configuration.

\paragraph{\textbf{Structural features only}}

In~\cite{yamada2018preconditioner}, CNNs are employed to select a preconditioner for GMRES by reformulating the task as an image classification problem.
Prior to applying the CNN, a feature image is constructed from the coefficient matrix.
To handle matrices of varying sizes, the original matrix is partitioned into 64 $\times$ 64 regions, with each region mapped to one pixel in the feature image.
The sparse coefficient matrix is thereby represented as a three-channel image.
The blue channel encodes the normalized size of the original matrix, the green channel represents the local density of nonzero entries in each region, and the red channel represents the average magnitude of nonzero values in each region.
The resulting feature image is then fed into a CNN classifier, which predicts the optimal preconditioner for the given system.

\paragraph{\textbf{Numerical and structural features}}

\citet{tang2022graph} propose a GNN-based model to select the optimal combination of solver and preconditioner, with a particular focus on Krylov subspace methods.
The problem is formulated as a multi-label classification task, where the model takes a coefficient matrix as input and predicts a binary vector.
The architecture contains two parallel branches.
One branch represents the sparse matrix as a graph and applies graph convolution and pooling layers to node and edge features to capture structural information, such as diagonal dominance.
In the other branch, global matrix features, such as the number of nonzero entries, are encoded as a feature vector and processed by an FNN.
The resulting two embeddings are then fused and passed to a classifier (another FNN) for final prediction.
During training, the performance of each combination of solver and preconditioner is evaluated using a score function that jointly considers runtime and residual error.

In~\cite{xiong2025mmautosolver}, a multimodal machine learning framework~\cite{morency2022tutorial} is introduced for automatically selecting the optimal combination of solvers and preconditioners for sparse linear systems.
Multimodal learning was originally designed to extract information from multiple data sources, such as images and audio.
Similar to~\cite{tang2022graph}, the method considers both numerical and structural features, treating them as two distinct modalities.
Specifically, numerical features and structural features are learned by an FNN and a CNN, respectively.
The learned representations are then fused into a classifier (another FNN) to predict the optimal combination of solver and preconditioner.

\subsubsection{Other approaches}\label{sec:linear-solver_selection-other}

Similar to~\cite{kuefler2008using}, \citet{khodak2025pcgbandit} propose the so-called ``PCGBandit'' method, which employs online learning, more specifically an adversarial bandit algorithm, to adaptively select preconditioner configurations for PCG.
This approach builds on~\cite{khodak2024learning}, where a similar bandit framework is used to tune the SOR relaxation parameter; see~\Cref{sec:linear-core_components-param}.
The PCGBandit method is designed for sequences of linear systems arising in time-dependent fluid and magnetohydrodynamic simulations.
It reformulates preconditioner selection along this sequence as a bandit problem: at each time step, an action corresponds to choosing a preconditioner configuration, and the feedback is the solution time of the resulting linear system.
The objective is to minimize the cumulative solution time over the whole simulation.

\citet{zabegaev2024automated} focus on the selection of linear solvers in time-dependent multiphysics simulations.
At each time step, the method aims to select an optimal solver configuration, including splitting schemes, Krylov subspace methods, and preconditioners.
This task can be formulated as a black-box optimization problem, where the performance of each candidate configuration can only be evaluated through expensive solver executions.
To address this, the authors adopt an online learning strategy: at each simulation time step, the performance of the selected solver configuration is used to update the selection model.
The resulting solver-configuration selection problem is then addressed by Bayesian optimization with Gaussian processes~\cite{liu2021gptune} or by a heuristic approach.

\paragraph{\textbf{Discussion of automated solver selection}}

\emph{Automated solver selection} has been studied for several decades and is particularly attractive in practice because it is non-intrusive and often does not require modifications to existing code.
Future work may focus on extending preconditioner selection methods to more complex settings, such as block preconditioners, where different blocks may employ different preconditioning strategies.
For example, in coupled systems, block-wise matrix features or the coupling information could be used to determine suitable preconditioners for different blocks.
Moreover, meta-learning frameworks have been applied to algorithm selection~\cite{smith2009cross}, and iterative solver selection can naturally be viewed as a specialized instance of algorithm selection.

\section{Machine learning-enhanced nonlinear iterative methods}\label{sec:nonlinear}

This section reviews ML-enhanced iterative methods for nonlinear systems, following an organization similar to that used for linear systems, as shown in~\Cref{fig:overview-HIMs}(b).
\Cref{sec:nonlinear-prelims} introduces the relevant classical nonlinear solvers, \Cref{sec:nonlinear-core_components} discusses learning their core components, and \Cref{sec:nonlinear-precond} focuses on ML-based nonlinear preconditioning, including nonlinear extensions of subspace correction methods.

\subsection{Preliminaries}\label{sec:nonlinear-prelims}

Consider a system of nonlinear equations
\begin{equation*}
    F(\mathbf{x}) = 0,
\end{equation*}
where $F: \mathbb{R}^N \rightarrow \mathbb{R}^N$ is a nonlinear function and $\mathbf{x} \in \mathbb{R}^{N}$ is the unknown vector.
Nonlinear systems are typically solved numerically by iterative methods that aim to compute an approximate solution $\mathbf{x}^* \in \mathbb{R}^N$ such that $F(\mathbf{x}^*) \approx 0$.
Following the notation for linear systems, $\mathbf{x}_k$ denotes the $k$-th iterate, and the corresponding residual is given by $F(\mathbf{x}_k)$.

Nonlinear systems can be solved by various iterative methods, including fixed-point iterations and Newton-type methods.
Gradient-based methods can also be applied to the residual merit function $\psi(\mathbf{x})=\frac{1}{2}\|F(\mathbf{x})\|_2^2$, whose gradient is $J(\mathbf{x})^{\mathsf T}F(\mathbf{x})$, where $J$ is the Jacobian of $F$. However, a stationary point of $\psi$ need not solve $F(\mathbf{x})=0$.
Fixed-point iterations are obtained by reformulating the original system into the form $\mathbf{x} = \phi(\mathbf{x})$, leading to the iteration
\begin{equation*}
    \mathbf{x}_{k+1} = \phi(\mathbf{x}_k),
\end{equation*}
starting from an initial guess $\mathbf{x}_0$.
Under a suitable contraction condition, this fixed-point iteration converges to a fixed point.
A vector $\hat{\mathbf{x}}$ is called a fixed point of $\phi$ if $\hat{\mathbf{x}} = \phi(\hat{\mathbf{x}})$.
A simple example is given by $\phi(\mathbf{x}) = \mathbf{x} + F(\mathbf{x})$.
In nonlinear solvers, fixed-point iterations are often called Picard iterations~\cite{picard1890memoire}, a term also used for the classical method of successive approximations for ordinary differential equation (ODE) initial-value problems and, in numerical PDEs, for linearizations obtained by lagging nonlinear terms. Although simple, Picard iterations may converge slowly.

Newton's method, also known as the Newton-Raphson method, is a widely used iterative method for solving systems of nonlinear equations and can also be regarded as a fixed-point iteration.
It is based on a first-order Taylor expansion of the nonlinear function $F(\mathbf{x})$ at the current iterate $\mathbf{x}_k$, which leads to a linearized approximation of the nonlinear problem.
At each iteration, Newton's method is given by
\begin{equation*}
    \mathbf{x}_{k+1} = \mathbf{x}_k - J(\mathbf{x}_k)^{-1}F(\mathbf{x}_k),
\end{equation*}
where $J(\mathbf{x}_k)$ denotes the Jacobian matrix of $F$ evaluated at $\mathbf{x}_k$, which is often large and sparse.
Here, a linearized system $J(\mathbf{x}_k) \mathbf{d}_k = -F(\mathbf{x}_k)$ must be solved.
Thus, each Newton step requires a linear solver, and linear iterative methods can be employed for this task.
Solving these linearized systems is typically the major part of the cost of Newton's method.
To reduce the cost, various variants of Newton's method have been developed.
For example, inexact Newton's method computes an approximate Newton direction $\mathbf{d}_k$ satisfying
\begin{equation*}
    \|J(\mathbf{x}_k) \mathbf{d}_k + F(\mathbf{x}_k)\| \leq \eta_k \|F(\mathbf{x}_k)\|,
\end{equation*}
where $\eta_k \in [0, 1)$ controls the accuracy of the linear solve.
Inexact Newton's method is often preferred in practice, especially for large-scale problems, as it allows the linearized systems to be solved approximately and reduces computational cost.

Newton's method is sensitive to the initial guess, and its quadratic convergence is typically guaranteed only when the initial guess is sufficiently close to the exact solution.
Globalization strategies are therefore used to improve robustness by enlarging the set of initial guesses from which Newton's method converges.
Common globalization techniques~\cite{MR2278446} include backtracking line searches and trust-region methods; cf.~\cite{dennis_numerical_1996,eisenstat_globally_1994}.
Other robustness-enhancing approaches include homotopy and continuation methods, pseudo-time stepping, mesh sequencing, and nonlinear preconditioning; cf.~\cite{allgower_continuation_1993,cliffe_numerical_2000,cresta_nonlinear_2007,hwang_parallel_2005,kelley1998convergence,knoll_jacobian-free_2004,kubicek_computational_1983,liu2022nonlinear,watson_globally_1990,watson_algorithm_1997}.
Nonlinear preconditioning~\cite{cai2002nonlinearly} improves robustness by transforming the nonlinear system before applying Newton-type iterations.
Nonlinear domain decomposition methods, including nonlinear additive Schwarz methods~\cite{hwang_parallel_2005,dolean2016nonlinear}, nonlinear elimination techniques~\cite{rose1992use,cai2011inexact}, and nonlinear FETI-DP and BDDC~\cite{klawonn2014nonlinear}, provide important examples of this perspective.
Multilevel nonlinear solvers and preconditioners, such as nonlinear multilevel domain decomposition preconditioners~\cite{klawonn2014nonlinear,heinlein2020additive,heinlein2023adaptive} and the full approximation scheme (FAS)~\cite{brandt1977multi}, a nonlinear multigrid method, exploit coarse-level nonlinear corrections to accelerate convergence and improve scalability.

\subsection{Core components learning}\label{sec:nonlinear-core_components}
We review methods for learning \emph{core components} of nonlinear iterative solvers: the \emph{initial guess}, \emph{update function}, and \emph{parameters}, discussed in~\Cref{sec:nonlinear-core_components-initial,sec:nonlinear-core_components-update,sec:nonlinear-core_components-param}, respectively.

\subsubsection{Initial guess learning}\label{sec:nonlinear-core_components-initial}

An appropriate initial guess is crucial for nonlinear solvers.
For example, Newton's method is highly sensitive to the choice of the initial guess.
A poor initial guess may result in slow convergence, stagnation, or even divergence.
ML techniques can therefore be employed to improve nonlinear solvers by generating more suitable initial guesses.
When combined with Newton's method, such approaches are often referred to as \emph{hybrid Newton's methods}~\cite{lechevallier2025hybrid}.
Similar to the linear case, neural network-based initial guess learning approaches often exploit the spectral bias of neural networks; \ie, networks trained for a relatively small number of iterations tend to capture the low-frequency components of the exact solution.
As a result, they can provide effective initial guesses for classical nonlinear solvers.

\citet{huang2020intdeep} propose a hybrid solver, called the deep learning initialized iterative (Int-Deep) method, which combines deep learning with classical nonlinear solvers, such as Newton's method and Picard iteration, for solving nonlinear PDEs.
The method consists of two phases.
In the first phase, the PDE is reformulated as an expectation minimization problem and solved using a mesh-free DNN.
Since the network is trained for only a limited number of iterations, the resulting solution is not sufficiently accurate.
In the second phase, this approximate solution is used as the initial guess for a classical nonlinear solver applied to the finite element-discretized system.
With such an initialization, the iterative solver converges rapidly and refines the solution to high accuracy.
The Int-Deep framework has subsequently been applied to 3D semilinear elliptic equations~\cite{jianguohuang2023hybrid}, elliptic variational inequalities of the second kind~\cite{cao2024hybrid}, and Navier--Stokes Darcy models~\cite{huang2025newtons}.

In~\cite{yahia2021mixing}, the authors propose a hybrid approach that combines Newton's method with neural networks for solving nonlinear systems arising from the inverse and direct kinematics (IK and DK) of cable-driven parallel robots (CDPRs).
In such problems, Newton's method is highly sensitive to the choice of initial guess, and random initializations often lead to divergence.
To address this issue, neural networks are trained in a supervised manner to learn the mapping from problem parameters to approximate solutions.
These predictions are then used as initial guesses for Newton's method, improving its convergence.

\citet{aghili2025accelerating} consider nonlinear parametric elliptic PDEs and propose a hybrid method in which an FNO is used to predict an approximate solution as the initial guess for Newton's method.
Specifically, the FNO is trained to learn the mapping from PDE parameters to solutions.
The training loss consists of three terms: a data-driven loss measuring the discrepancy between predicted and reference solutions, a gradient loss, and a discretization-informed loss based on the residual of the discretized PDE.

In~\cite{lechevallier2025hybrid, lechevallier2024physics}, the authors propose a hybrid Newton's method for $\text{CO}_2$ storage simulations, which involve a time-dependent problem.
Due to the use of an implicit time integration scheme, a system of nonlinear equations has to be solved at every time step to compute the saturation, which is computationally expensive.
To alleviate this cost, an FNO is employed to learn the mapping from the problem parameters and the saturation at the previous time step to the saturation at the current time step.
The FNO prediction is then used as the initial guess for Newton's method.

\citet{jin2025fast} propose a neural hybrid Newton's method for solving stiff nonlinear time-evolution PDEs discretized by implicit time integration schemes, where a system of nonlinear equations must be solved at each time step.
Their method employs a CNN to approximate the mapping from the solution at the previous time step to the solution at the current time step.
The CNN is trained in an unsupervised manner by minimizing the residual of the implicit time-stepping scheme.
The CNN prediction is then used as the initialization for Newton's method.

\citet{taghikhani2025neuralinitialized} propose a hybrid approach, called neural-initialized Newton (NiN), for parametric nonlinear systems arising in computational solid mechanics.
The method employs a neural operator to learn the mapping from problem parameters to the solution function and uses the predicted solution as an initialization for Newton's method.
Specifically, the authors consider implicit finite operator learning (iFOL)~\cite{najianasl2026physicsinformed}, a physics-informed neural operator framework that incorporates finite element residuals directly into the loss function.
The resulting prediction is then refined through Newton's method applied to nonlinear systems arising from finite element discretizations.

Moreover, hybrid Newton methods have been developed and applied across a broad range of problems in computational science and engineering, including simulation of flow in porous media~\cite{abbasi2022improved}, real-time simulation of elastic structures~\cite{odot2022deepphysics}, heterogeneous elastic localization in computational mechanics~\cite{kelly2025thermodynamicallyinformed}, fluid-structure interaction problems~\cite{tiba2025machinelearning}, and parametric transient problems~\cite{papadopoulos2024fusing}.

\paragraph{\textbf{Discussion}}

In the nonlinear setting, leveraging ML models to improve initial guesses has been widely studied, as it is relatively simple to implement and initialization often has a stronger influence on convergence than in linear iterative methods.
Similar to the linear case, any ML model that provides approximate solutions can be used for \emph{initial guess learning}.
Moreover, for nonlinear problems, obtaining an initial guess within the basin of attraction of the nonlinear solver may be more important than producing a highly accurate approximation, although most existing methods primarily aim to improve the accuracy of the initial guess.
Indeed, nonlinear iterations may stagnate even when the initial guess is close to the exact solution, while rapid convergence may still be achieved from a less accurate initial guess lying within the basin of attraction~\cite{susanto2009newton}.
Another challenge is the selection of informative input features, which need to capture the nonlinear behavior of \(F\) and are therefore often strongly problem dependent.

\subsubsection{Update function learning}\label{sec:nonlinear-core_components-update}

In standard Newton's method, each iteration requires solving a linearized system to compute the update, which can be computationally expensive.
Moreover, for nonlinear PDEs with multiple solutions, Newton's method may suffer from instability near bifurcation points.
To address these challenges, \citet{hao2024newton} propose the Newton-informed neural operator method, where a neural operator $\mathcal{N}$ (implemented via DeepONet) is trained to approximate the Newton update mapping, \ie, $\Delta u = \mathcal{N}(u)$, with $u$ denoting the current iterate and $\Delta u$ the corresponding update.
The training objective consists of two parts. The first is a supervised loss, which matches the neural operator output with the Newton updates in discrete form.
The second is the Newton loss (a physics-informed loss), which enforces the linearized equation at the operator level, \ie, $J(u)\mathcal{N}(u) = -\mathcal{F}(u(x))$, where $\mathcal{F}$ denotes the nonlinear differential operator and $J(u)$ is its Fr\'echet derivative.
Once trained, the model can generate multiple solutions from different initial guesses, and is capable of discovering solutions not included in the training dataset.

\paragraph{\textbf{Discussion}}
Existing work on \emph{update function learning} for nonlinear iterative methods remains very limited.
Beyond the ``Newton-informed neural operator method'' in~\cite{hao2024newton}, one can consider learning more general update formulas, for example by incorporating information from previous iterates or residuals.
Moreover, \emph{update function learning} is closely related to the idea of ``learning to optimize''~\cite{chen2022learning}, where ML models are trained to learn iterative update rules for optimization algorithms.

\subsubsection{Parameter learning}\label{sec:nonlinear-core_components-param}

Numerical relaxation techniques can be employed to accelerate fixed-point iterations by introducing a relaxation parameter at each iteration.
\citet{silva2021machine} propose two general frameworks that employ random forest regression~\cite{breiman2001random} to predict optimal relaxation parameters for nonlinear iterative methods applied to multiphase porous media flow.
These frameworks take dimensionless numbers and simulation properties as input and select the relaxation parameter to minimize the number of nonlinear iterations.
In particular, they are applied to the Picard iteration with numerical relaxation~\cite{salinas2017improving}.
This approach is further extended to an online learning setting in~\cite{silva2023adaptive}, where the model is continuously refined during the simulation.

The pseudo-transient continuation method~\cite{kelley1998convergence} is a globalization technique for Newton's method, in which each iteration involves solving an associated time-dependent problem.
A local pseudo-time stepping strategy can be used to accelerate convergence and improve the convergence region, especially for nonlinear systems arising from finite element discretizations.
Here, locality is defined with respect to mesh-dependent quantities.
Classical approaches typically rely on a global Courant-Friedrichs-Lewy (CFL) number~\cite{CMRM} to select pseudo-time steps.
To improve this, neural networks are used in~\cite{zandbergen2025improving} to predict spatially varying pseudo-time steps based on local information.
Since the learned parameter is defined at the element level and directly modifies the solver, this approach can be regarded as intrusive.
Moreover, it can also be interpreted as a nonlinear preconditioner.

\paragraph{\textbf{Discussion}}
Learning parameters in other nonlinear iterative methods could also be explored, such as step sizes in line search methods or parameters in trust-region methods.

\subsection{Nonlinear preconditioner learning}\label{sec:nonlinear-precond}

The classical Newton's method often suffers from stagnation or even failure due to unbalanced nonlinearities, where a few dominant nonlinear components slow down the overall convergence~\cite{cai2002nonlinearly}.
Nonlinear preconditioning techniques aim to identify and balance such nonlinearities, thereby improving the convergence and robustness.
Existing nonlinear preconditioning approaches often rely on problem-specific knowledge, careful parameter tuning, or classical nonlinear domain decomposition methods, such as nonlinear additive Schwarz methods~\cite{hwang_parallel_2005,luo2023pin^la-settled}.
ML techniques have therefore been leveraged to construct or enhance nonlinear preconditioning strategies, including approaches related to nonlinear domain decomposition and subspace correction methods discussed below.

\citet{luo2023pin^la-settled} propose a nonlinear preconditioning method, referred to as the preconditioned inexact Newton with learning capability ($\text{PIN}^\mathcal{L}$), to efficiently identify and mitigate unbalanced nonlinearities, where a few dominant nonlinear components severely hinder the convergence of Newton iterations.
The $\text{PIN}^\mathcal{L}$ method is an unsupervised learning approach based on POD.
During the training stage, a suitably chosen training problem is solved using a standard inexact Newton's method to collect residual and solution snapshots, which are then used to construct projection operators for the residual and solution spaces, respectively.
Then, a nonlinear preconditioning stage uses these projection operators to construct an approximate system of nonlinear equations that captures the dominant low-frequency nonlinear structures.
Several low-dimensional subspace Newton iterations with projected Jacobian systems are performed to generate an improved initial guess.
Finally, starting from this improved initial guess, the original system of nonlinear equations is solved by a global inexact Newton's method.
Therefore, this method can also be interpreted as a form of \emph{initial guess learning}~(\Cref{sec:nonlinear-core_components-initial}).
Building on $\text{PIN}^\mathcal{L}$, \citet{gong2026stagnationshortening} propose a stagnation-shortening inexact Newton (SSIN) method, in which the nonlinear preconditioner is constructed online from residual and solution snapshots collected directly during the stagnation phase of the target problem, rather than from a separate training problem as in~\cite{luo2023pin^la-settled}.

\citet{lee2025neuraloperator} propose a neural preconditioned Newton (NP-Newton) method for solving parametric nonlinear systems arising from PDE discretizations.
The approach introduces a nonlinear right preconditioner based on a neural operator model, called the fixed-point neural operator (FPNO).
Given a system of nonlinear equations $F(\mathbf{x})=0$, the method considers the right-preconditioned system $F(G(\mathbf{x}))=0$, where $G$ denotes the nonlinear right preconditioner implemented by the FPNO.
The FPNO is designed to emulate a fixed-point iteration and is defined as $G(\mathbf{x}) = \mathbf{x} + \alpha G_B(\mathbf{x})$, where $G_B$ is a backbone neural operator and $\alpha$ is an adaptive step size predicted from the nonlinear residual through an FNN and a hyperbolic tangent function.
At each nonlinear iteration, the method first applies neural preconditioning to obtain an improved iterate, followed by a standard Newton correction.
This method can also be regarded as a nonlinear extension of \emph{alternating HIMs}~(\Cref{sec:linear-precond-alternating}), \ie, a HIM that alternates between standard Newton's method and a learned fixed-point iteration.

As in the linear case, subspace correction methods for nonlinear systems can also be viewed as a nonlinear preconditioning technique.
The nonlinear FETI-DP method~\cite{klawonn2014nonlinear} is a nonlinear non-overlapping DD method that has been successfully applied to various problems.
Adaptive coarse spaces developed for the linear FETI-DP method can accelerate linear convergence, and the same coarse spaces can also be incorporated into nonlinear FETI-DP to improve both linear and nonlinear convergence~\cite{heinlein2023adaptive}.
\citet{klawonn2024learninga} propose learning coarse constraints for the linear FETI-DP method using FNNs; see~\Cref{sec:linear-mgdd-transfer_operator} for more details.
This idea is further extended to the nonlinear FETI-DP framework in~\cite{klawonn2025learning}.
Moreover, the automated design of GMG cycle structures~\cite{schmitt2020constructing, schmitt2021evostencils} is discussed in~\Cref{sec:linear-mgdd-parameters} and has been extended to nonlinear problems via the full approximation scheme (FAS), as demonstrated in~\cite{parthasarathy2023evolving}.

\paragraph{\textbf{Discussion}}

Although many methods for \emph{initial guess learning} for nonlinear iterative methods have been proposed, research on \emph{ML-enhanced nonlinear iterative methods} remains relatively limited, particularly for \emph{update function learning} and \emph{nonlinear preconditioning techniques}.
This is mainly because nonlinear iterative methods are generally more difficult to stabilize and analyze than linear solvers.
For example, learning effective nonlinear preconditioners is typically more difficult than learning linear preconditioners.
Another challenge lies in identifying suitable features and generating effective training data, since nonlinear iterative behavior often depends strongly on the underlying problem structure and initialization.
Nevertheless, many techniques developed for linear problems (discussed in~\Cref{sec:linear}) may also be applicable to nonlinear settings, although their effectiveness is problem dependent.
For instance, when Newton's method is employed, the resulting linearized systems can be solved using \emph{ML-enhanced linear iterative methods}.

Existing literature primarily focuses on employing ML for nonlinear right-preconditioning, which preserves the original nonlinear function; however, nonlinear left-preconditioning remains relatively underexplored.
\emph{Alternating hybrid iterative methods}, as discussed in~\Cref{sec:linear-precond-alternating}, may also be extended to nonlinear solvers.
In addition, other subspace correction methods developed for linear solvers, as discussed in~\Cref{sec:linear-mgdd}, may also be extended in a similar manner.
If nonlinear convergence can be improved at relatively low additional cost, the resulting computational savings can be substantial, especially for problems requiring many nonlinear iterations.
Moreover, since neural networks are inherently nonlinear models, they may be particularly well suited for learning nonlinear corrections or preconditioning operators.

\section{Related topics}\label{sec:related}

The integration of ML techniques with classical iterative methods, not necessarily limited to methods for solving systems of equations, has attracted increasing attention in recent years.
This paper focuses specifically on ML-enhanced iterative methods for solving systems of linear and nonlinear equations.
In addition, several other active research topics are closely related to the theme of this paper.
In this section, we briefly review some closely related topics.

In~\Cref{sec:related-time_inte}, we review ML-enhanced time integration methods.
Time integration methods advance the solution iteratively over successive time steps, and their performance can therefore be improved using ML models.
In~\Cref{sec:related-inspired_nn}, we briefly discuss papers on iterative method-inspired training data, neural network architectures, and training algorithms.
Concepts and techniques from classical iterative methods can be incorporated into ML methods through these ingredients.
\Cref{sec:related-direct} reviews ML-enhanced direct solvers for linear systems.
Direct solvers are also widely used in applications, and ML can be employed to reduce their computational cost or improve their performance, for example, by automatically selecting matrix reorderings to reduce fill-in.
The purpose of this section is to provide a brief review of several recent works on these related topics rather than an exhaustive survey.

\subsection{Machine learning-enhanced time integration methods}\label{sec:related-time_inte}
Time integration methods are widely used for solving time-dependent problems arising from ODEs and PDEs.
These methods are generally classified as explicit or implicit schemes.
Explicit methods advance the solution to the next time step using only previously computed solution values and do not require solving a system of equations at each time step.
They are usually simple and computationally efficient, but their stability generally imposes step-size restrictions, including CFL-type restrictions for suitable PDE discretizations.
In contrast, implicit methods generally provide better stability and allow larger time steps, but they require solving a linear or nonlinear system at each time step, which can be computationally expensive.

The performance of time integration methods depends strongly on the choice of time-step size, stability, accuracy, and computational cost.
To improve these aspects, ML techniques have recently been introduced into time integration methods.
Here, we focus only on forward problems, \ie, combining ML and time integration methods for solving time-dependent problems.
Related works on inverse problems and model discovery can be found in~\cite{qin2019data, qin2021datadriven, rudy2019deep, raissi2018multistep, keller2021discovery}.

\citet{um2020solverintheloop} propose a solver-in-the-loop framework for iterative PDE solvers, with time integration as a specific example.
Their method uses neural networks to learn correction terms for discretization errors in numerical update rules, where the correction model is trained through interaction with the solver during time integration.
\citet{liu2022hierarchical} propose a hierarchical neural time-stepping framework (HiTS) that learns flow maps (or update rules) at multiple temporal scales using a ResNet-based model.
The method combines short-step and long-step neural integrators so that small-scale models provide local accuracy while large-scale models mitigate long-term error accumulation.

\citet{yan2024time} propose learning time integration schemes for solving time-dependent PDEs on coarse spatial grids.
In particular, neural networks are employed to learn the coefficients of linear multistep methods so as to compensate for the numerical errors introduced by coarse spatial discretization.
For explicit time integration methods, small time steps are accurate but slow, while large time steps are fast but inaccurate or even unstable.
Based on this, \citet{huang2023fasta} propose a deep learning-based corrector called neural vector (NeurVec).
In this method, a classical numerical solver is performed with a large time step, followed by a neural network correction term that learns the integration error caused by using a coarse time step.

\citet{sheth2021intelligent} employ an ML method to select the time-step size in reservoir simulation problems that are large but still efficient to solve.
Specifically, at each time step, random forests or neural networks are used to learn the mapping from simulation features to the appropriate next time-step size.
These features include throughput values (related to CFL conditions), throughput gradients, previous Newton iterations, and previous time-step sizes.
Training data is generated using an ``aggressive'' time-stepping strategy that explores large time steps and retains only successfully converged steps for training.
Closely related are ML-enhanced pseudo-time-stepping techniques, which can be used as globalization strategies to improve nonlinear convergence.
In such methods, ML can be used to predict local pseudo-time steps or related parameters, thereby accelerating convergence and enlarging the effective convergence region.
One such approach is proposed in~\cite{zandbergen2025improving}; see also~\Cref{sec:nonlinear}.

\subsection{Machine learning inspired by classical iterative methods}\label{sec:related-inspired_nn}

As discussed in~\Cref{sec:prelims-ml}, an ML system typically consists of three key ingredients: \emph{data}, \emph{models}, and \emph{learning}~\cite{deisenroth2020mathematics}.
In this subsection, we review how ideas from iterative methods can be used to improve each of these ingredients.
First, in~\Cref{sec:related-inspired_nn-data}, we consider enhancements to the \emph{data} component, focusing on \emph{training data sampling} strategies inspired by multigrid (or multilevel) methods.
Second, in~\Cref{sec:related-inspired_nn-architectures}, we consider the \emph{model} component, which, in the context of neural networks, corresponds to the \emph{architecture}.
In particular, classical iterative methods such as multigrid and domain decomposition have inspired novel neural network architectures through their multilevel structures and divide-and-conquer principles, respectively.
This topic is closely related to ML-enhanced MG and DD methods discussed in~\Cref{sec:linear-mgdd}; however, those methods explicitly construct and solve linear systems using MG and DD methods, whereas the present discussion focuses on neural network-based discretization and architectures.
Third, in~\Cref{sec:related-inspired_nn-training_algorithms}, we briefly discuss the \emph{learning} component, where training algorithms are improved using preconditioning techniques.

\subsubsection{Training data sampling}\label{sec:related-inspired_nn-data}

Neural networks are trained on data, and their performance is strongly influenced by the quality and distribution of the training samples.
Data sampling strategies inspired by multigrid (or multilevel) methods can therefore be employed to improve both training efficiency and generalization performance.

To mitigate the spectral bias in PINNs, \citet{tsai2025mldpinn} propose a multilevel dataset training strategy, called MLD-PINNs.
By switching between datasets with different numbers of collocation points in a multigrid-inspired training sequence (\eg, a V-cycle), the model can effectively capture both low- and high-frequency components.
Neural operator learning typically requires high-resolution training data, which can be computationally expensive.
To address this issue, \citet{rowbottom2026multilevel} propose a multilevel Monte Carlo (MLMC) approach to accelerate the training.
Specifically, MLMC is used to estimate the expected gradients of the loss function across multiple grid levels: coarse grids capture most of the variance, while a few fine-grid samples are used to refine the estimate.

\subsubsection{Neural network architectures}\label{sec:related-inspired_nn-architectures}

Classical iterative method-inspired neural network architectures have received increasing attention in recent years.
The design of many neural network architectures is guided by underlying mathematical principles, which makes them more interpretable than purely black-box models.
In particular, classical iterative methods such as MG and DD have been introduced into neural network design.
Incorporating these principles into neural networks may improve their robustness and enhance their ability to capture complicated features of the underlying problem.

Ideas from iterative methods have also inspired the design of neural network architectures.
For example, in residual neural networks (ResNet)~\cite{he2016deep}, the network is designed to learn residual mappings, where the input $x$ is added to the output of the residual mapping, yielding $y = x + \mathcal{N}(x)$.
This discrete iterative perspective is further extended by neural ODEs~\cite{chen2018neural}, which introduce a continuous-depth formulation in which the evolution of the hidden state is governed by a learnable differential equation.
Similarly, in~\cite{boudec2024learning, boudec2025learning}, the solution of a PDE is represented in an iterative form, and a neural network is used to learn update terms.
This is closely related to \emph{update function learning} (see~\Cref{sec:linear-core_components-update}), but without relying on linear systems derived from classical discretization methods.

Drawing inspiration from classical Krylov subspace methods and Runge-Kutta methods, \citet{doncevic2024recursively} propose a generic framework for solving systems of equations and time-integration problems, called the recursively recurrent neural network (R2N2) model, which consists of an inner recurrence and an outer recurrence.
The inner recurrence is used to generate an information subspace, such as the Krylov subspace in Krylov methods or the multiple integration stages in Runge-Kutta methods.
The outer recurrence then determines an update within the generated subspace.
The corresponding coefficients are learned through meta-learning.

The multigrid idea has also been incorporated into neural network architectures.
\citet{he2019mgnet} propose a unified framework called MgNet that connects CNNs with multigrid methods, where feature extraction and pooling operations correspond to smoothing and restriction steps, respectively.
MgNet has been applied to image processing tasks~\cite{he2019mgnet} and operator learning for PDEs~\cite{zhu2024enhanced}.
In~\cite{wang2026learnable}, graph convolutions are embedded into the smoothing and transfer operations of a GMG V-cycle, yielding a learnable multigrid architecture.

\citet{riccietti2021multilevel} introduce the multilevel idea into PINNs by representing the PDE solution as the sum of a coarse term and a fine term.
These terms are modeled by two neural networks with different numbers of parameters.
The networks are trained alternately using training data at different resolutions.
\citet{wang2024multistage} propose multi-stage neural networks (MSNNs) to approximate a function, where a sequence of networks is trained iteratively to learn residuals.
This results in a decomposition of the form $\mathcal{N}(x)=\sum_{i=1}^{\#_\text{net}} \mathcal{N}_i(x)$, where $\#_\text{net}$ denotes the number of networks.

Similarly, \citet{yang2026multilevel} propose a multilevel deep learning framework for solving PDEs, where the final solution is expressed as $\mathcal{N}(x)=\sum_{i=1}^{\#_\text{net}} \mathcal{N}_i(x)$.
Their method decomposes the training process into multiple levels, where each level aims to learn the residual errors left by the previous ones.
At each level, an adaptive sampling strategy is employed to concentrate training points in regions with large residuals and sharp gradients.
\citet{guo2025mgfno} integrate multigrid concepts into FNO to better capture high-frequency components, where three different networks are trained to represent the solution at different resolutions.
\citet{le2021surrogate} introduce multigrid ideas into the U-Net architecture for surrogate modeling in fluid dynamics.
\citet{he2024mgno} propose MgNO, a neural operator that incorporates linear operators into its formulation.
A multigrid structure is used to parameterize these operators.

DD principles, which reflect a divide-and-conquer strategy, have also inspired neural network architectures.
The central idea is to divide the computational domain into subdomains, train local neural solvers on each subdomain, and then combine them using overlapping or non-overlapping DD strategies.
We refer to~\cite[Section~4]{klawonn2024machine} for a comprehensive review of domain decomposition-inspired neural network architectures.

In~\cite{li2020d3m, li2020deepa}, the Deep Ritz method~\cite{e2018deep} and PINNs~\cite{raissi2019physicsinformed} are employed as subdomain solvers within overlapping Schwarz methods, respectively.
Inspired by two-level DD methods, \citet{mercier2021coarse} introduce a coarse network to handle the coarse problem, thereby improving scalability, \ie, the ability to handle an increasing number of subdomains.
\citet{sun2024domaina} employ Deep Ritz networks as subdomain solvers within a non-overlapping DD framework.
\citet{klaes2025flowrateconserving} combine alternating Schwarz methods with CNN-based surrogate models for macroscopic blood flow simulations, where the CNNs act as subdomain solvers.

\citet{jagtap2020conservative} propose conservative PINNs (cPINNs), which introduce non-overlapping DD into PINNs by training separate PINNs on local subdomains and coupling them through conservation conditions at the interfaces.
Building on cPINNs, the authors of~\cite{jagtap2020extended} propose extended PINNs (XPINNs), which generalize this framework to general PDEs and arbitrary space-time domain decompositions.
Motivated by the finite element method, \citet{moseley2023finite} propose the finite basis PINNs (FBPINNs) method, in which the PDE solution is expressed as $u(x)=\sum_i w_i(x)\mathcal{N}_i(x)$.
Here, the functions $w_i(x)$ are window functions based on overlapping domain decomposition that localize the domain, and the neural networks $\mathcal{N}_i(x)$ learn local basis functions.
FBPINNs have been further extended to multilevel variants~\cite{dolean2024multilevel}, time-dependent problems and operator learning~\cite{heinlein2025multifidelity}, adaptive subdomains~\cite{botvinick-greenhouse2025abpinns}, and randomized neural networks~\cite{anderson2024elmfbpinn, shang2025overlappinga}.

Moreover, DD ideas have also been applied to image processing tasks.
For image classification tasks, \citet{gu2022decomposition} propose decomposing a global CNN into several subnetworks by partitioning its width, \ie, along the channel dimension.
These subnetworks are trained independently and in parallel.
A transfer learning strategy is then used to compose the trained subnetworks and initialize the global CNN, which is subsequently fine-tuned to further improve its performance.
\citet{klawonn2024domain} incorporate DD concepts into CNNs for image classification by partitioning an image into multiple subimages.
Each subimage is processed by a local CNN to produce a local decision, and these decisions are then aggregated by a global DNN to obtain the final prediction.
The domain decomposition-based U-Net (DDU-Net) for image segmentation proposed in~\cite{verburg2025ddunet} integrates DD with the U-Net architecture.
In this approach, the input image is partitioned into several overlapping subimages, and the U-Net is applied to each subimage to extract local features.

\subsubsection{Training algorithms}\label{sec:related-inspired_nn-training_algorithms}
For a neural network $\mathcal{N}_\theta(x)$ with parameters $\theta$, training involves minimizing a loss function $\mathcal{L}(\theta)$ over a dataset $\mathcal{D}$.
    Many first- and second-order gradient-based optimization methods can be interpreted through preconditioning or curvature-based scaling of the update direction.
    A generic preconditioned gradient iteration takes the form
\[
    \theta^{(k+1)}
    =
    \theta^{(k)}
    -
    \eta M_k^{-1}
    \nabla_\theta \mathcal{L}(\theta^{(k)};\mathcal{D}),
\]
where $\eta$ denotes the step size (or learning rate) and $M_k^{-1}$ scales or preconditions the gradient.
    When $M_k=I$, the iteration is standard gradient descent.
    Adaptive first-order methods such as AdaGrad~\cite{duchi2011adaptive} use diagonal, iteration-dependent scaling based on past gradients, while Adam~\cite{kingma2017adam} combines such adaptive scaling with a momentum term based on first-moment estimates.
    When $M_k=H(\theta^{(k)})$, where $H(\theta^{(k)})$ denotes the Hessian matrix, the iteration corresponds to Newton's method, which incorporates curvature information.
    For nonlinear least-squares problems, Gauss-Newton replaces the Hessian by the approximation $J_k^TJ_k$, where $J_k$ is the Jacobian of the residual.
    Quasi-Newton methods, such as L-BFGS~\cite{liu1989limited}, construct approximations to the Hessian or its inverse.
    More recently, Shampoo~\cite{gupta2018shampoo} and SOAP~\cite{vyas2025soap} construct structured preconditioners that exploit the matrix structure of gradients to capture richer curvature information.
    We refer to~\cite{kopanicakova2026introduction} for further discussion on optimization methods for training SciML models.
    Ideas from classical iterative methods, including preconditioning, can be incorporated into optimization algorithms for neural network training, as discussed below.

\citet{liu2024preconditioning} focus on employing neural networks to learn the solution of linear systems arising from the finite difference method.
In particular, they introduce an ILU preconditioner to improve the condition number of the system.
Similarly, \citet{song2025matrix} consider PINNs based on CNNs (PICNNs) for solving PDEs discretized by the finite difference method.
For linear PDEs, the resulting system is often ill-conditioned, which leads to slow convergence during training.
To address this issue, they propose a matrix preconditioning framework based on an ILU factorization of the Jacobian matrix, which is computed via automatic differentiation and matrix coloring.
\citet{wen2026preconditioned} introduce a preconditioned physics-informed loss for neural operator training to mitigate the ill-conditioning.

Multigrid (or multilevel) ideas have been introduced into neural network training.
\citet{gunther2020layerparallel} propose a layer-parallel training strategy for ResNets~\cite{he2016deep} by interpreting the network as a forward Euler discretization of a dynamical system.
From this perspective, the training process can be reformulated as an optimal control problem, where the network parameters act as control variables.
To solve this problem efficiently, they employ the multigrid-reduction-in-time (MGRIT) method~\cite{falgout2014parallel}, where the time domain corresponds to the layer domain.
This allows for parallel computation across layers and replaces the conventional sequential propagation with a nonlinear multigrid method.
Moreover, this layer-parallel training paradigm has been extended to other architectures, including gated recurrent units (GRUs)~\cite{cho2014properties} in~\cite{moon2022parallel} and transformers~\cite{vaswani2017attention} in~\cite{jiang2026layerparallel}.

Several papers introduce domain decomposition-based nonlinear preconditioning into the training of neural networks.
These approaches can be broadly divided into two classes: splitting the training parameters $\theta$ and splitting the training dataset $\mathcal{D}$.
First, based on nonlinear preconditioning techniques~\cite{cai2002nonlinearly}, \citet{kopanicakova2024enhancing} propose a nonlinear preconditioning framework for L-BFGS to improve the training of PINNs in stiff and ill-conditioned problems.
In particular, they introduce one-level additive and multiplicative Schwarz methods by partitioning the network parameters in a layerwise manner.
Moreover, \citet{lee2026twolevel} extend it to the two-level case, where a coarse-level training step is introduced to accelerate the convergence.
Similarly, the multi-preconditioned L-BFGS method for training FBPINNs is proposed in~\cite{salvado-benasco2026multipreconditioned}; see also~\cite{moseley2023finite,dolean2024multilevel}.
\citet{trotti2025parallel} introduce the additively preconditioned trust-region strategy (APTS)~\cite{gross2009unifying} to neural network training, where the trust-region method is a globalization technique for optimization problems.
Second, \citet{a.cruzalegria2025dataparallel} build on APTS to propose a data-parallel training method for neural networks.

\subsection{Machine learning-enhanced direct solvers for linear systems}\label{sec:related-direct}

Direct methods, such as LU, Cholesky, and QR factorizations, solve linear systems through matrix factorization.
For large sparse systems, however, factorization may introduce substantial fill-in, increasing both computational and memory costs~\cite{davis2006direct}.
Matrix reordering is therefore commonly applied before factorization to reduce fill-in, although finding an optimal ordering is NP-hard~\cite{yannakakis1981computing}.
Recently, ML has been explored for selecting or constructing suitable reorderings; below, we briefly review several representative studies.

\citet{cui2019machinelearningdriven} consider DD-based power grid analysis, where linear systems on subdomains are solved using sparse direct solvers.
To reduce the cost of direct solves, they employ SVMs or neural networks to select an appropriate matrix ordering strategy for reducing fill-in from three candidates, based on matrix features such as the number of nonzeros.
This idea is similar to \emph{automated solver selection} based on matrix features~(\Cref{sec:linear-solver_selection-numerical_structural_features}).
Similar ML-based selection strategies are also considered in~\cite{tang2025selection,chen2023optimized}.

\citet{gatti2022graph} propose using reinforcement learning and GNNs for graph partitioning.
Specifically, a GraphSAGE-based network is trained using the advantage actor-critic (A2C) algorithm~\cite{williams1992simple}.
Rather than directly learning a matrix ordering, the method recursively refines vertex separators, which are then used to construct a nested dissection ordering to reduce fill-in in sparse LU factorization.
Analogously, reinforcement learning-based reordering methods for reducing fill-in in sparse LU factorization have also been proposed in~\cite{dasgupta2023alpha, harb2025refill}.

\subsection{Discussion of related topics}\label{sec:other-discussion}

A central open question for the related topics discussed in this section is how learned components can be used without losing the safeguards of the underlying numerical methods.
For ML-enhanced time integration, this means understanding which predictions can be safely embedded into classical schemes and under which conditions they preserve stability, error control, and long-time reliability.
Such methods should therefore be assessed not only in terms of short-term accuracy or speedup, but also in terms of their robustness outside the training regime and their interaction with classical adaptivity mechanisms.

For iterative-method-inspired ML, the central challenge is to turn structural analogies into reliable design principles.
Architectures inspired by unrolling, residual or recurrent updates, neural ODEs, DD, or multigrid encode useful structural ideas such as repeated correction, locality, hierarchy, and parallelism.
However, many assumptions underlying the analysis of classical iterative methods are not satisfied by neural network models, so such connections are often heuristic rather than theoretically justified.
    Nevertheless, some of these ideas have been shown to improve properties such as training stability or the mitigation of spectral bias, whereas in other cases the gains are relatively modest.
    Although approximation results are available for many architectures, they do not yet explain when and why these structural choices improve training, generalization, or convergence in scientific computing tasks.
Future work should identify and theoretically analyze which structural ingredients are responsible for the observed gains and whether they remain effective in large-scale implementations, where communication, memory use, and scalability become decisive.

As for ML-enhanced direct solvers, existing studies have mainly focused on matrix reordering for reducing fill-in.
A possible research direction is to evaluate ML-based reordering methods not only in terms of fill-in reduction, but also in terms of practical factorization time and memory consumption, as demonstrated in~\cite{li2026factorizationinloop}.
Moreover, a systematic comparison between ML-enhanced iterative and direct methods remains lacking, particularly in terms of computational efficiency, memory consumption, robustness, and scalability across different problem classes.

\section{Conclusions and future directions}\label{sec:conclusions}

This review has examined \emph{hybrid iterative methods (HIMs)}, that is, ML-enhanced iterative methods for systems of linear and nonlinear equations.
Most existing work concerns linear systems, whereas nonlinear HIMs are less extensively studied and remain an important direction for future research.
We have also briefly discussed related techniques, including ML-enhanced time-stepping methods, ML approaches inspired by classical iterative algorithms, and ML-enhanced direct solvers for linear systems.

\paragraph{\textbf{Key observations from existing studies}}
Two central observations emerge from the existing literature on HIMs:
\begin{enumerate}
    \item \textbf{HIMs augment rather than replace classical iterative methods.}
          HIMs are not purely ML-based replacements for classical solvers.
          Instead, they integrate ML-based components into established iterative frameworks.
          Classical iterative methods offer robustness, mathematical structure, interpretability, and convergence guarantees, while ML can provide efficient data-driven construction or adaptation of key solver components, such as initial guesses and preconditioners; see, for example, \Cref{sec:linear-core_components,sec:linear-precond}.
          In this way, HIMs aim to improve computational efficiency while preserving, whenever possible, desirable properties of classical solvers that are often unavailable in black-box ML solvers.
    \item \textbf{Classical iterative methods and ML models can provide complementary numerical capabilities.}
          Many classical methods for systems arising from differential equations have strong smoothing properties: high-frequency error components are reduced rapidly, whereas low-frequency components often converge slowly.
          This behavior is closely related to the locality of the underlying discretizations, which are typically based on local basis functions and sparse interactions.
          In contrast, neural networks provide nonlinear and nonlocal representations that can capture global solution features.
          They are therefore natural candidates for approximating global or low-frequency components that are difficult for local relaxation schemes.
          This complementary behavior is reflected in several HIMs, including \emph{initial guess learning} and \emph{alternating HIMs}; see \Cref{sec:linear-core_components-initial,sec:linear-precond-alternating}, respectively.
          Other ML techniques, such as PCA- or POD-based approaches, also play an important role, although neural network-based methods currently dominate recent activities.
\end{enumerate}

\paragraph{\textbf{Main future directions}}
Based on these observations, we identify two main directions for future research on HIMs:
\begin{enumerate}
    \item \textbf{ML can be used to reduce the computational cost of classical solvers.}
          HIMs can provide computational benefits when savings in repeated or parametric solves outweigh the costs of data generation, training, and inference (see~\Cref{sec:linear-core_components-initial}).
          Trained models can approximate or select costly solver components, such as preconditioners, coarse-space corrections, transfer operators, initial guesses, update directions, and parameter settings.
          In this way, HIMs can exploit the approximation capabilities of ML models and, in particular, the ability of neural networks to approximate nonlinear mappings, to represent complex solver operations with reduced online cost.
          The usefulness of such predictions, however, depends on how they are embedded into the surrounding solver.
          When combined with classical correction, relaxation, projection, or Krylov-type mechanisms, ML-based approximations can improve efficiency while retaining important numerical safeguards.
    \item \textbf{ML may be particularly useful for problems for which classical solver theory provides limited guidance.}
          A particularly promising direction for HIMs is to target problems for which classical theory provides only limited guidance, or for which the most effective solver or solver configuration is not known a priori.
          For model problems such as the Poisson equation, whose standard discretizations lead to symmetric positive definite systems, classical iterative methods are comparatively mature and well understood.
          In contrast, the design of effective solvers can be considerably more challenging for nonsymmetric, indefinite, saddle-point, or highly nonlinear systems with complex and problem-dependent properties.
          In such settings, data-driven approaches to constructing or selecting solvers, preconditioners, solver components, and algorithmic parameters may be especially valuable, in particular when existing theory is difficult to extend.
          The current development and limitations of nonlinear HIMs are discussed in~\Cref{sec:nonlinear}.
\end{enumerate}

\paragraph{\textbf{Challenges and practical applicability}}
Despite the progress reviewed in this paper, several challenges remain:
\begin{enumerate}
    \item \textbf{Generalization, systematic evaluation, and appropriate integration with classical solvers remain central challenges.}
          Learned components must generalize beyond the training distribution, their offline training cost must be justified by repeated or large-scale use, and their integration should respect the mathematical structure of the underlying problem.
          Although many existing studies on HIMs have reported improvements over classical iterative methods, these gains are often problem-dependent and may not persist across different problem settings.
          Direct comparisons across existing studies should therefore be interpreted with caution, as they often use different test problems, ML implementations, and evaluation metrics.
          Systematic benchmarking under common settings remains an important direction for future research.
          Choosing ML models that match the structure of the underlying problem, such as GNNs for unstructured grids or sparse matrix-valued data, also remains essential.
    \item \textbf{Real-world validation and large-scale applicability require substantially more investigation.}
          When validating HIMs on real-world applications and large-scale HPC platforms, scalability, communication overhead, end-to-end efficiency, robustness, generalization, and reliability are critical considerations.
          Existing studies still primarily focus on relatively simple problems, and the performance of HIMs on more complex and application-relevant problems, involving complicated geometries, heterogeneous or multiscale coefficients, large-scale discretizations, or strongly coupled multiphysics systems, remains to be systematically investigated.
          Moreover, \emph{initial guess learning}, \emph{alternating HIMs}, and \emph{automated solver selection} (see~\Cref{sec:linear-core_components-initial,sec:linear-precond-alternating,sec:linear-solver_selection}) may have particular practical applicability to real-world problems, as they can often be integrated into established classical solvers with relatively minor modifications to existing implementations.
          Representative tests include nonlinear hyperelasticity problems~\cite{wang2026pretrain}, incompressible fluid simulations~\cite{lan2024neuralpreconditioned}, mixed-dimensional elliptic problems~\cite{dimola2026numerical}, Helmholtz scattering~\cite{lee2025fast}, pore-scale mechanics~\cite{li2024machine}, $\text{CO}_2$ storage~\cite{lechevallier2025hybrid}, and fluid-structure interaction~\cite{tiba2025machinelearning}.
          More broadly, ML-enhanced finite element procedures have also been explored for practical applications~\cite{bathe2026enhancing}.
          Finally, while this paper focuses primarily on square linear and nonlinear systems of equations, HIMs may also be extended to other classes of problems, such as overdetermined linear systems and more general nonlinear optimization problems.
\end{enumerate}

\section*{Acknowledgements}

The author Yuhuang Meng would like to acknowledge the financial support from the China Scholarship Council (CSC).

\section*{Declarations}

\begin{itemize}
    \item \textbf{Competing interests}: The authors declare no competing interests.
    \item \textbf{Data availability}: No datasets were generated or analyzed during the current study.
    \item \textbf{Author contributions}: Y.M. conducted the primary literature search, reviewed and categorized the relevant references, and prepared the original draft of the manuscript. A.H. conceived the original idea for the article and identified additional relevant references. Y.M., J.Z., and A.H. contributed to the conceptualization of the article, discussed its content, and reviewed and revised the manuscript.
\end{itemize}

\appendix
\section{Open-source code and datasets}\label{sec:open_code}

\subsection{Classical iterative solvers and matrix datasets}\label{sec:open_code-classical}
Several open-source packages for classical sparse iterative solvers and sparse matrix datasets are listed in~\Cref{tab:open_code-classical}.

\begin{table}[!ht]
    \footnotesize
    \centering
    \caption{
        Packages on classical sparse iterative solvers and sparse matrix datasets.
    }
    \label{tab:open_code-classical}
    \begin{tabular}{
            lp{8cm}l
        }
        \toprule
        \textbf{package name}
            & \textbf{description}
            & \textbf{download link}
        \\
        \midrule
        PETSc
            & the portable, extensible toolkit for scientific computation
            & \url{https://petsc.org/release/}
        \\
        Trilinos
            & collection of open-source packages for scientific application development
            & \url{https://trilinos.github.io/}
        \\
        Hypre
            & high-performance preconditioners
            & \url{https://www.llnl.gov/casc/hypre/}
        \\
        Ginkgo
            & high-performance numerical linear algebra library
            & \url{https://gko-project.org/}
        \\
        SuiteSparse Matrix Collection
            & formerly known as the University of Florida Sparse Matrix Collection
            & \url{https://sparse.tamu.edu/}
        \\
        \bottomrule
    \end{tabular}
\end{table}

\subsection{ML-enhanced linear iterative methods}

Open-source code and datasets for \emph{core components learning}~(\Cref{sec:linear-core_components}), \emph{preconditioner learning}~(\Cref{sec:linear-precond}), \emph{enhancements of subspace correction methods}~(\Cref{sec:linear-mgdd}), and \emph{automated solver selection}~(\Cref{sec:linear-solver_selection}) are provided in~\Cref{tab:open_code-core_components,tab:open_code-precond,tab:open_code-mgdd,tab:open_code-solver_selection}, respectively.

\begin{table}[!ht]
    \footnotesize
    \centering
    \caption{
        Open-source code for \emph{core components learning}~(\Cref{sec:linear-core_components}), including \emph{initial guess learning}, \emph{update function learning}, and \emph{parameter learning}.
    }
    \label{tab:open_code-core_components}
    \begin{tabular}{
            lp{8.8cm}l
        }
        \toprule
        \textbf{ref.}
            & \textbf{description}
            & \textbf{download link}
        \\
        \midrule
        \cite{ozbay2021poisson}
            & Poisson CNN, learning initial guesses for multigrid \& BiCGSTAB
            & \url{https://github.com/aligirayhanozbay/poisson_CNN}
        \\
        \cite{luna2021accelerating}
            & learning initial guesses for GMRES
            & \url{https://github.com/ML4FnP/GMRES-Learning}
        \\
        \cite{eshaghi2026nows}
            & NOWS, neural operator warm starts
            & \url{https://github.com/eshaghi-ms/NOWS}
        \\
        \cite{wang2026pretrainingfinetuning}
            & a pretraining-finetuning computational framework for material homogenization
            & \url{https://github.com/yizheng-wang/HomoGenius}
        \\
        \cite{wang2026pretrain}
            & a pretraining and warm-start framework for PDEs
            & \url{https://github.com/yizheng-wang/PFEM}
        \\
        \cite{hsieh2019learning}
            & enhancing classical iterative methods for PDEs
            & \url{https://github.com/ermongroup/Neural-PDE-Solver}
        \\
        \cite{stanziola2021helmholtz}
            & learning an update function for Helmholtz equations
            & \url{https://github.com/ucl-bug/helmnet}
        \\
        \cite{khodak2024learning}
            & learning to relax: learning SOR parameters using reinforcement learning
            & \url{https://github.com/mkhodak/learning-to-relax}
        \\
        \cite{arisaka2023principled, arisaka2024accelerating}
            & GBMS and NI-GBMS, gradient-based meta-solving for parameter learning
            & \url{https://github.com/arisakaso/nigbms}
        \\
        \bottomrule
    \end{tabular}
\end{table}

\begin{table}[!ht]
    \footnotesize
    \centering
    \caption{
        Open-source code for \emph{preconditioner learning}~(\Cref{sec:linear-precond}).
    }
    \label{tab:open_code-precond}
    \begin{tabular}{
            llp{7cm}
        }
        \toprule
        \textbf{ref.}
            & \textbf{description}
            & \textbf{download link}
        \\
        \midrule
        \cite{ackmann2020machinelearned}
            & machine-learned preconditioners for GCR
            & \url{https://github.com/JanAckmann/MLPrecon}
        \\
        \cite{rudikov2024neural}
            & FCG-NO, neural operator-based preconditioner for FCG
            & \url{https://github.com/arudikov/FCG-NO}
        \\
        \cite{chen2025graph}
            & GNP, Graph neural preconditioners for FGMRES
            & \url{https://github.com/jiechenjiechen/GNP}
        \\
        \cite{kaneda2023deep}
            & DCDM, deep conjugate direction method
            & \url{https://github.com/ayano721/2023_DCDM}
        \\
        \cite{lan2024neuralpreconditioned}
            & neural-preconditioned steepest descent with orthogonalization
            & \url{https://github.com/kai-lan/MLPCG}
        \\
        \cite{bai2026hybrid}
            & DLSM and HyDEA, alternating hybrid iterative method
            & \url{https://github.com/HMB9666/HyDEA}
        \\
        \cite{hausner2024neural, hausner2024learning}
            & neural incomplete Cholesky and LU factorizations
            & \url{https://github.com/paulhausner/neural-incomplete-factorization}
        \\
        \cite{sappl2026convolutional}
            & sparse approximate inverse preconditioner for PCG
            & \url{https://github.com/jsappl/DeepPreconditioning}
        \\
        \cite{yang2025learning}
            & SPAI preconditioner learning for PCG solvers on GPUs
            & \url{https://github.com/Adversarr/LearningSparsePreconditioner4GPU}
        \\
        \cite{millard2025pearl}
            & PEARL, preconditioner learning using reinforcement learning
            & \url{https://github.com/djm3622/precondition-discovery-contextual-bandit}
        \\
        \cite{bankestad2024ising}
            & sparsity pattern learning for SPAI preconditioners
            & \url{https://github.com/mariabankestad/IsingOnGraphs}
        \\
        \cite{sun2026learning}
            & singularity-encoded Green's function learning
            & \url{https://github.com/q1sun/SgEncd_Green_Function}
        \\
        \cite{gotz2018machine}
            & sparsity pattern learning for block-Jacobi preconditioners
            & \url{https://github.com/Markus-Goetz/block-prediction/}
        \\
        \cite{trappler2025statedependent}
            & state-dependent preconditioners for variational data assimilation
            & \url{https://github.com/VTrappler/ML_preconditioners}
        \\
        \cite{zhang2024blending}
            & HINTS, combining relaxation methods and ML-precon. iterations
            & \url{https://github.com/kopanicakova/HINTS_precond}
        \\
        \cite{versano2024attentionbased}
            & a HINTS variant based on a geometry-aware DeepONet
            & \url{https://github.com/idanv87/convnet}
        \\
        \cite{cui2022fourier, cui2025hybrid}
            & FNS, Fourier neural solver
            & \url{https://github.com/cuichen1996/FourierNeuralSolver}
        \\
        \cite{wu2026are}
            & physics-aware Anderson acceleration
            & \url{https://github.com/HNU-WYH/DL-HIM}
        \\
        \bottomrule
    \end{tabular}
\end{table}

\begin{table}[!ht]
    \footnotesize
    \centering
    \caption{
        Open-source code for \emph{enhancements of subspace correction methods}~(\Cref{sec:linear-mgdd}).
    }
    \label{tab:open_code-mgdd}
    \begin{tabular}{
            lp{8cm}p{7cm}
        }
        \toprule
        \textbf{ref.}
            & \textbf{description}
            & \textbf{download link}
        \\
        \midrule
        \cite{huang2023learning}
            & learning multigrid smoothers
            & \url{https://github.com/jerryhuangru/Learning-optimal-multigrid-smoothers-via-neural-networks}
        \\
        \cite{azulay2023multigridaugmented}
            & multigrid-augmented deep learning preconditioners
            & \url{https://github.com/BGUCompSci/HelmholtzAccCNN.jl}
        \\
        \cite{lerer2024multigridaugmented}
            & multigrid-augmented deep learning preconditioners using compact implicit layers
            & \url{https://github.com/BGUCompSci/CompactImplicitHelmholtz}
        \\
        \cite{taghibakhshi2021optimizationbased}
            & algebraic multigrid coarsening using reinforcement learning
            & \url{https://github.com/compdyn/rl_grid_coarsen}
        \\
        \cite{caldana2024deep}
            & AMG-ANN, learning the strong threshold parameter using CNNs
            & \url{https://github.com/MatteoCaldana/AMG-ANN}
        \\
        \cite{zou2024autoamgsettled}
            & AutoAMG($\theta$), learning the strong threshold parameter using GNNs
            & \url{https://github.com/zhf-0/autoamg}
        \\
        \cite{katrutsa2017deep,katrutsa2020blackbox}
            & deep multigrid, learning transfer operators and damping factor in the Jacobi smoother
            & \url{https://github.com/amkatrutsa/dmg}
        \\
        \cite{greenfeld2019learning}
            & learning prolongation operators
            & \url{https://github.com/danielgreenfeld3/Learning-to-optimize-multigrid-solvers}
        \\
        \cite{luz2020learning}
            & learning prolongation operators for PDEs on unstructured grids
            & \url{https://github.com/ilayluz/learning-amg}
        \\
        \cite{taghibakhshi2023mggnn}
            & MG-GNN, learning Robin-type interface conditions and coarse spaces
            & \url{https://github.com/JRD971000/Code-Multilevel-MLORAS}
        \\
        \cite{luo2024neural}
            & NeurKItt, neural Krylov iteration for accelerating linear system solving
            & \url{https://github.com/smart-JLuo/NeurKItt}
        \\
        \cite{schmitt2021evostencils}
            & EvoStencils, a grammar-based genetic programming approach for constructing efficient multigrid methods
            & \url{https://github.com/jonas-schmitt/evostencils}
        \\
        \cite{taghibakhshi2022learning}
            & learning interface conditions in domain decomposition solvers
            & \url{https://github.com/compdyn/learning-oras}
        \\
        \bottomrule
    \end{tabular}
\end{table}

\begin{table}[!ht]
    \footnotesize
    \centering
    \caption{
        Open-source code for \emph{automated solver selection}~(\Cref{sec:linear-solver_selection}).
    }
    \label{tab:open_code-solver_selection}
    \begin{tabular}{
            lp{8.5cm}l
        }
        \toprule
        \textbf{ref.}
            & \textbf{description}
            & \textbf{download link}
        \\
        \midrule
        \cite{khodak2025pcgbandit}
            & PCGBandit, learning preconditioners for PCG using online learning
            & \url{https://github.com/mkhodak/PCGBandit}
        \\
        \cite{zabegaev2024automated}
            & automated solver selection for simulation of multiphysics processes in porous media
            & \url{https://github.com/Yuriyzabegaev/solver_selector}
        \\
        \cite{weng2026embeddingbased}
            & embedding-based methods for linear solver performance prediction
            & \url{https://gitlab.com/h.liu/embedding-based-predictor}
        \\
        \bottomrule
    \end{tabular}
\end{table}

\subsection{ML-enhanced nonlinear iterative methods}

Open-source code and datasets for \emph{ML-enhanced nonlinear iterative methods}~(\Cref{sec:nonlinear}) are provided in~\Cref{tab:open_code-nonlinear}.

\begin{table}[!ht]
    \footnotesize
    \centering
    \caption{
        Open-source code for \emph{ML-enhanced nonlinear iterative methods}~(\Cref{sec:nonlinear}).
    }
    \label{tab:open_code-nonlinear}
    \begin{tabular}{
            lp{7.3cm}p{6.8cm}
        }
        \toprule
        \textbf{ref.}
            & \textbf{description}
            & \textbf{download link}
        \\
        \midrule
        \cite{lechevallier2025hybrid, lechevallier2024physics}
            & hybrid Newton's method
            & \url{https://github.com/PINN-Well-opening-and-closing-events/Yads}
        \\
        \cite{odot2022deepphysics}
            & DeepPhysics: a physics aware deep learning framework for real-time simulation
            & \url{https://github.com/AlbanOdot/DeepPhysics-article}
        \\
        \cite{kelly2025thermodynamicallyinformed}
            & TherINO: Thermodynamically-informed iterative neural operators
            & \url{https://github.com/conlain-k/therino}
        \\
        \cite{hao2024newton}
            & Newton informed neural operator for computing multiple solutions of nonlinear PDEs
            & \url{https://github.com/xlliu2017/Newton-Informed-Neural-Operator}
        \\
        \cite{silva2023adaptive}
            & adaptive learning acceleration for nonlinear PDE solvers
            & \url{https://github.com/vlssanonymous/mlsolveracc}
        \\
        \cite{zandbergen2025improving}
            & improving pseudo-time stepping convergence for CFD simulations with neural networks
            & \url{https://github.com/searhein/nn-pseudo-time-stepping-data}
        \\
        \bottomrule
    \end{tabular}
\end{table}

\bibliography{references}

\end{document}